\pdfoutput=1

\RequirePackage{fix-cm}

\documentclass[leqno, fleqn, 12pt]{article}

\usepackage{latexsym}
\usepackage{amsmath}
\usepackage{amscd}
\usepackage{amssymb}
\usepackage{verbatim}

\usepackage[T1]{fontenc}

\usepackage[upright]{fourier}

\usepackage{fullpage}

\usepackage{tikz-cd}

\tikzcdset{row sep/special/.initial=12ex}
\tikzcdset{row sep/spe/.initial=10.5ex}
\tikzcdset{column sep/spec/.initial=5em}
\tikzcdset{column sep/specc/.initial=8em}
\tikzcdset{column sep/speccc/.initial=10em}

\newskip\nineskipamount \nineskipamount=9pt plus 0pt minus 0pt
\newskip\zeroskipamount \zeroskipamount=0pt plus 0pt minus 0pt
\usepackage[noorphans,vskip=\nineskipamount]{quoting}

\makeatletter
\renewcommand{\@makefntext}[1]{\vspace*{0.5ex}\parindent=0em
\hspace*{-0.4em}
\hbox to 0.4em{\hss\@makefnmark}\hspace*{0.4em}{#1}
}
\makeatother

\newcounter{mysectionnumber}
\newcommand{\mysection}[2]{\setcounter{footnote}{0}
\setcounter{myparnum}{0}
\refstepcounter{mysectionnumber}
\vspace{21pt}{\Large {\themysectionnumber.} {#1}}\label{#2}\vspace*{15pt}}

\newcommand{\mynonumbersection}[1]{
\vspace{21pt}
{\Large{}\hspace*{0.00em}$\phantom{1.}$ {#1}\vspace*{15pt}}
}

\newcommand{\myit}[1]{\textbf{\textit{#1}}\hspace{0.0em}}

\newcounter{myparnum}
\newcounter{mylemmanum}[myparnum]
\newcommand{\myuppar}[1]{\vspace{\medskipamount}\textbf{#1}\hspace*{0.5em}}

\newcommand{\myitpar}[1]{\vspace{\medskipamount}\textit{\textbf{#1}}\hspace*{0.5em}}

\newcommand{\myepar}[1]{\vspace{\medskipamount}\textit{#1}\hspace*{0.5em}}

\newcounter{myappendnumber}
\newcounter{myaparnum}[myappendnumber]
\newcommand{\myappend}[2]{\setcounter{footnote}{0}
\setcounter{myaparnum}{0}
\refstepcounter{myappendnumber}
\vspace{21pt}{\Large A\dff.{\themyappendnumber.}\oss {#1}}\label{#2}\vspace*{15pt}}

\newcommand{\proof}{\vspace{\medskipamount}{\textbf{{\emph{Proof}.}}\hspace*{1em}}}
\newcommand{\prooftitle}[1]{\vspace{\medskipamount}{\textbf{{\emph{#1}.}}\hspace*{1em}}}

\newcommand{\eproof}{ $\blacksquare$}

\newcommand{\dis}{\displaystyle}

\def\sss{\hspace{0.05em}\ }
\def\dss{\hspace{0.1em}\ }
\def\trs{\hspace{0.15em}\ }
\def\qss{\hspace{0.2em}\ }
\def\pss{\hspace{0.3em}\ }
\def\oss{\hspace{0.4em}\ }

\def\halfff{\hspace*{0.025em}}
\def\fff{\hspace*{0.05em}}
\def\dff{\hspace*{0.1em}}
\def\trf{\hspace*{0.15em}}
\def\qff{\hspace*{0.2em}}
\def\pff{\hspace*{0.3em}}
\def\off{\hspace*{0.4em}}

\newcommand{\hnsp}{\hspace*{-0.05em}}
\newcommand{\nsp}{\hspace*{-0.1em}}
\newcommand{\nnsp}{\hspace*{-0.15em}}
\newcommand{\snsp}{\hspace*{-0.175em}}
\newcommand{\dnsp}{\hspace*{-0.2em}}

\renewcommand{\leq}{\leqslant}
\renewcommand{\geq}{\geqslant}
\newcommand{\id}{\mathop{\mbox{id}}\nolimits}

\newcommand{\zzz}{\mathbf{Z}}
\newcommand{\rrr}{\mathbb{R}}

\newcommand{\ftwo}{\mathbb{F}_{\dff 2}}

\newcommand{\hclass}[1]{[\dff #1 \dff]}
\newcommand{\fclass}[1]{[\snsp [\dff #1\dff]\snsp ]}
\newcommand{\num}[1]{\left|\qff #1 \qff\right|}
\newcommand{\norm}[1]{\|\qff #1 \qff\|}
\newcommand{\sco}[1]{\langle\dff #1 \dff\rangle}

\newcommand{\bd}{\operatorname{bd}}

\renewcommand{\int}{\operatorname{int}}
\newcommand{\st}{{\operatorname{st}\trf}}
\newcommand{\bs}{{\operatorname{b}\fff}}
\newcommand{\ibs}[1]{{\operatorname{b}^{ #1}\dff}}
\newcommand{\cs}{{\operatorname{c}\trf}}
\newcommand{\bigfclass}[1]{\biggl[\dnsp\nsp \biggl[\dff #1\dff\biggr]\dnsp\nsp \biggr]}

\newcommand{\toto}{\longrightarrow}
\newcommand{\ttoo}{\hspace*{0.2em}\longrightarrow\hspace*{0.2em}}

\begin{document}

\setlength{\baselineskip}{12pt plus 0pt minus 0pt}
\setlength{\parskip}{12pt plus 0pt minus 0pt}
\setlength{\abovedisplayskip}{12pt plus 0pt minus 0pt}
\setlength{\belowdisplayskip}{12pt plus 0pt minus 0pt}

\newskip\smallskipamount \smallskipamount=3pt plus 0pt minus 0pt
\newskip\medskipamount   \medskipamount  =6pt plus 0pt minus 0pt
\newskip\bigskipamount   \bigskipamount =12pt plus 0pt minus 0pt

\title{{The\pss lemmas\pss of\pss Alexander\qss and\pss Sperner}}
\date{} 
\author{\textnormal{Nikolai\qss V.\qss Ivanov{}}}

\footnotetext{\hspace*{-0.65em}\copyright\oss 
Nikolai\qss V.\qss Ivanov,\oss 2019.\trs 
Neither the work reported in the present paper\halfff,\qss
nor its preparation were supported by any corporate entity.}

\footnotetext{\hspace*{-0.65em}The author\dss is\dss grateful to\dss M.\qss Prokhorova\dss 
for correction of\dss several\dss infelicities and\dss the suggestion\dss to expand\dss the
discussion of\trs graphs\sss and\dss path-following\sss algorithms.}

\maketitle

\vspace*{8ex}

{\renewcommand{\baselinestretch}{1}
\selectfont

\myit{\hspace*{0em}\large Contents}\vspace*{1.5ex} \\ 
\hbox to 0.8\textwidth{\myit{Preface}\hfil 2} \hspace*{0.5em} \vspace*{1.5ex}\\
\hbox to 0.8\textwidth{\myit{\phantom{1}1.}\hspace*{0.5em} Simplicial\sss complexes and\sss chains\hfil 5} \hspace*{0.5em} \vspace*{0.25ex}\\
\hbox to 0.8\textwidth{\myit{\phantom{1}2.}\hspace*{0.5em} A\dss lemma\sss of\qss Alexander\hfil 9} \hspace*{0.5em} \vspace*{0.25ex}\\
\hbox to 0.8\textwidth{\myit{\phantom{1}3.}\hspace*{0.5em} Brouwer's\dss invariance of\dss 
dimension\dss theorem\hfil 12}  \hspace*{0.5em}\vspace*{0.25ex}\\
\hbox to 0.8\textwidth{\myit{\phantom{1}4.}\hspace*{0.5em} Brouwer's\dss invariance of\dss 
domain\dss theorem\hfil 16} \hspace*{0.5em} \vspace*{0.25ex}\\
\hbox to 0.8\textwidth{\myit{\phantom{1}5.}\hspace*{0.5em} Brouwer's\qss fixed-point\dss theorem\hfil 19}\hspace*{0.5em} \vspace*{0.25ex}\\
\hbox to 0.8\textwidth{\myit{\phantom{1}6.}\hspace*{0.5em} Homology\dss groups 
and\dss their\dss topological\dss invariance\hfil 25} \hspace*{0.5em} \vspace*{0.25ex}\\
\hbox to 0.8\textwidth{\myit{\phantom{1}7.}\hspace*{0.5em} Sperner's\qss lemma and\dss its\dss
combinatorial\dss proof\hfil 30}  \hspace*{0.5em} \vspace*{0.25ex}\\
\hbox to 0.8\textwidth{\myit{\phantom{1}8.}\hspace*{0.5em} Cochains\sss and\trs 
Sperner's\qss lemma\hfil 32}  \hspace*{0.5em} \vspace*{0.25ex}\\
\hbox to 0.8\textwidth{\myit{\phantom{1}9.}\hspace*{0.5em} Graphs\sss and\dss path-following\dss
algorithms\hfil 36} \hspace*{0.5em} \vspace*{0.25ex}\\
\hbox to 0.8\textwidth{\myit{10.}\hspace*{0.5em} Cohomology\dss groups\hfil 40}  \hspace*{0.5em}  \vspace*{1.5ex}\\
\myit{Appendices}\hspace*{0.5em}  \hspace*{0.5em} \vspace*{1.5ex}\\
\hbox to 0.8\textwidth{\myit{\phantom{1}A.1.}\hspace*{0.5em} Barycentric subdivisions\hfil 47}\hspace*{0.5em} \vspace*{0.25ex}\\
\hbox to 0.8\textwidth{\myit{\phantom{1}A.2.}\hspace*{0.5em} Stellar subdivisions\hfil 50}\hspace*{0.5em}  \vspace*{1.5ex}\\
\hbox to 0.8\textwidth{\myit{References}\hspace*{0.5em}\hfil 54}\hspace*{0.5em}  \vspace*{0.25ex}  

}

\vspace*{2ex}

\renewcommand{\baselinestretch}{1}
\selectfont

{\small
Most\dss of\trs theorems and\dss lemmas in\dss this paper are\sss
unnumbered,\pss
but\sss often\dss have names such as\qss ``Alexander's\dss lemma''\pss
or\qss ``\dnsp$\partial\partial$\dnsp-theorem''.\oss
They\sss are\sss
referred\dss to by\dss their names or\dss their\dss places in sections.\oss
For\dss the numbered\dss theorems and\dss formulas,\pss continuous numbering\dss is\dss used.\oss

}

\newpage
\mynonumbersection{Preface}

\vspace*{6pt}
The\qss Brouwer\qss fixed-point\qss theorem,\oss 
due to\qss L.E.J.\dss Brouwer\qss \cite{b},\oss 
is one of\dss the best\dss known and\dss useful\sss theorems\sss in\sss topology\halfff.\oss
According\dss to\oss J.\dss Dieudonn\'{e}\qss \cite{d},\pss 
this\sss theorem is one of\dss the\vspace{-10.5pt}

\begin{quoting}
\ldots\oss epoch-making results of\dss Brouwer in\dss 1910\dff--1912,\oss
which may rightly be called the\qss \emph{first\dss proofs}\qss in algebraic topology,\oss
since Poincar\'{e}\halfff's papers can only be considered as blueprints for theorems to come.\oss\vspace{-1.5pt}

\ldots\oss In a rapid succession of\dss papers published in less\sss than\dss two years,\oss
the\qss ``Brouwer theorems''\qss
(as\sss they\sss are still called)\qss made him famous overnight\halfff.\oss
\end{quoting}

\vspace{-10.5pt}
See\qss \cite{d},\oss pp.\qss 161\dss and\dss 167\dss respectively.\oss
Very soon\dss Brouwer's\dss methods turned\sss into\sss the main\dss tool of\dss the\trs
\emph{combinatorial\dss topology}\qss of\dss the time,\oss
which in\sss the late\dss 1920is\dss morphed\sss into the\trs \emph{algebraic\dss topology}\halfff.\oss
It\dss is hardly\dss surprising\dss that\dss Brouwer's\dss methods 
were considered as sophisticated\dss
by\dss many of\dss his contemporaries.\oss
The algebraic topology\halfff,\pss including\sss the modern versions of\dss
Brouwer's\dss proof\dss of\dss his fixed-point\sss theorem,\pss 
also has reputation of\trs being sophisticated.\oss

The unwillingness of\dss non-topologists to read at\sss least\sss 100\dss pages,\qss
usually preceding\dss the proof\dss of\dss Brouwer theorem
in algebraic topology textbooks,\qss is understandable.\oss
The most popular way\dss to avoid\dss this onerous task is\dss to prove a\qss ``combinatorial''\qss
lemma of\qss Sperner\qss \cite{s}\qss and\dss then deduce Brouwer\dss theorem
from it\dss by\dss an argument of\qss
Knaster-Kuratowski-Mazurkiewich\qss \cite{kkm}.\oss 
The latter seems\sss to be so simple\sss that\dss its authors are often\dss not\sss even\dss mentioned,\pss
creating\dss the impression\dss that\dss the whole proof\dss   
is due\sss to\dss Sperner\halfff.\oss
In any case,\oss Sperner's\dss lemma is invariably praised as an ingenious
and\dss surprising\dss replacement\sss of\dss the 
machinery of\dss algebraic topology at\dss least\dss in\dss the 
proof\dss of\dss the\qss Brouwer\dss fixed-point\dss theorem.\oss

Another popular way\dss to avoid\dss learning algebraic topology
is\dss based\sss on\sss an\sss analytical\dss proof\dss of\dss Brouwer's\dss theorem
due to\dss Dunford-Schwartz\qss \cite{ds}.\oss
The\dss Dunford-Schwartz\dss proof\dss turned out\dss to be\sss
the usual\dss topological\dss proof\trs in a disguise\qss \cite{i1}.\pss
It\dss is a cochain-level\sss 
version of\dss the standard\dss proof\dss based\sss on\dss de Rham\dss cohomology\qss 
(in\dss de Rham\dss theory cochains are nothing else but\sss differential forms),\oss 
written in the language of\dss elementary\dss multivariable calculus.\oss
Remarkably\halfff,\pss this is\sss true not\sss only in a vague\qss ``moral''\qss sense,\pss
but\dss also on\dss the level of\dss minutiae details.\oss 

After discovering\dss this\dss  
I\dss started\dss to suspect\sss that\dss 
the proof\dss based on\dss Sperner's\dss lemma
is another disguise of\dss the usual\dss topological\dss proof\halfff.\oss
Online discussions with\dss the late\dss A.\dss Zelevinsky\sss and\sss with\dss F.\dss Petrov\dss 
kindled\dss
my\dss interest\dss further\halfff.\oss
The suspicion\dss turned out\dss to be correct\halfff,\pss
and\dss the\sss two situations\sss to be very\sss similar\halfff.\oss
Moreover\halfff,\oss this proof\dss 
turned out\dss to be\sss
much closer to\sss the ideas of\dss the classical\sss
combinatorial\qss (nowadays algebraic)\qss topology\dss than\dss
Dunford-Schwartz\dss proof\halfff.\dff\oss
Sperner's\trs lemma\dss turned out\dss to be a cochain-level\sss
version of\trs the standard\qss (simplicial\fff)\qss cohomological\sss arguments.\oss
The\dss Knaster-Kuratowski-Mazurkiewich\dss argument\dss turned out\dss to be
closely\dss related\dss to a long\dss proven\dss tool of\dss topologists,\oss
the simplicial\sss approximation of\dss continuous maps.\oss
A condensed account\sss of\trs these ideas appeared as\sss the note\qss \cite{i2}.\oss

The real surprise was awaiting\dss in\dss the footnote on\dss p.\qss 48\qss of\dss 
the cited above book\qss \cite{d}\qss
by\pss
J.\qss Dieudonn\'{e},\oss
and\dss in\dss the papers and\dss books\sss to which\dss this\sss footnote led.\qff\oss
In\dss 1928,\pss when\dss it\dss was\sss published,\oss
Sperner's\qss lemma\dss was\sss hardly\dss new\dss or\dss surprising\halfff.\oss

About\dss two years before\trs Sperner's\trs paper\halfff,\oss
Alexander\qss published a fundamental\dss paper\qss \cite{a2}\qss
devoted\dss to a proof\dss of\dss the topological\dss invariance of\trs combinatorially defined\qss
\emph{Betti\dss numbers}\qss and\qss \emph{torsion coefficients}\qss of\dss polyhedra.\oss
This was,\pss in\dss fact\halfff,\pss a proof\dss of\trs the\sss topological\dss
invariance of\qss \emph{homology\dss groups}\qss of\dss polyhedra,\pss
but\dss this\dss language was still\dss in\dss the future.\oss
A key element of\dss Alexander's\qss proof\qss is\dss a\sss lemma at\dss the end of\trs his paper\qss
(see\qss \cite{a2},\pss p.\qss 328\halfff).\oss
In\dss the above mentioned\dss footnote\pss J.\qss Dieudonn\'{e}\qss
writes\sss that\qss\vspace{-8pt}

\begin{quoting}
This lemma is a special case of\dss the Sperner\dss lemma,\pss 
proved\dss two years later\dss by\dss the same\sss method\qss
(\dff\cite{ah},\oss p.\qss 376\dff).
\end{quoting}

\vspace{-8pt}
This seems\sss to be an oversimplification,\pss certainly\dss suitable for a footnote.\oss
The methods are\qss ``morally''\qss the same,\oss
but\dss this\dss is\dss hardly\sss obvious.\oss
Praising\trs Sperner's\dss lemma as a double-count\-ing\sss
alternative\sss to\sss the methods of\dss 
algebraic\sss topology\dss
is\dss hardly\sss compatible even\dss with\dss the claim\dss that\dss the methods
are\qss ``morally''\qss the same.\oss

Alexander\sss stated\sss and\sss proved\sss 
his lemma in terms of\trs  
the combinatorial\dss topology\halfff.\oss
In\dss this context\qss 
Alexander's\qss lemma\dss appears\dss inevitably\halfff,\pss 
being\sss exactly\dss what\dss is\sss needed\dss for\dss 
his\dss proof\dss of\dss the topological invariance of\qss homological\dss invariants.\oss
While\sss the assumptions of\trs Sper\-ner's\dss lemma
look\sss somewhat strange,\pss
and are,\pss apparently\halfff,\pss interpreted\dss as a strike of\dss a genius by some,\pss
the corresponding\sss assumptions of\qss Alexander's\dss lemma\sss 
are only\dss natural.\oss

Another surprise is\dss the simplicity\sss of\qss Alexander's\dss proof\halfff.\oss
His proof\dss is simpler\dss than\dss Sperner's\dss one,\pss
at\dss least\dss for\dss mathematicians\sss comfortable with 
using\dss linear algebra over\dss the field
of\dss two elements\dss $\ftwo$\dss in\dss the spirit\sss 
of\trs linear algebra methods in combinatorics.

Some\sss things\dss
are needed\dss to be pointed\sss out\dss
on\dss a\sss technical\dss level\qss (as opposed\dss to\sss the\qss ``moral''\qss one).\oss
First\halfff,\pss Alexander's\dss proof\dss works without\sss any\sss changes\sss in\dss
the more general\sss case considered\dss by\dss Sperner\halfff.\oss
Second,\oss the special case considered\dss by\dss Alexander\dss is sufficient\dss
for all\dss topological\sss applications.\oss
Moreover\halfff,\oss Sperner\qss himself\dss silently\dss used only\dss this special\sss case\qss
(as also\qss Knaster-Kuratowski-Mazurkiewich).\oss
In\dss more details,\oss both\qss Alexander's\dss and\dss Sperner's\dss lemmas\sss
are concerned\dss 
with\dss subdivisions of\dss a simplex\dss into smaller\sss simplices.\oss
For applications one needs subdivisions into
arbitrarily\sss small\sss simplices.\oss
At\dss the\sss time\sss
the only\sss way\dss to get\dss such subdivisions\dss
was\sss to construct\dss them as\sss the so-called\qss \emph{iterated\dss barycentric subdivisions}.\oss
While\dss Sperner\dss simply\dss ignores\sss the question of\dss existence,\oss
Alexander\dss works with\dss the iterated\dss barycentric subdivisions.

At\dss the same\sss time\dss Alexander's\dss result\dss is\sss stronger\dss than\dss
Sperner's\dss one.\oss
Sperner\dss proves\sss that\dss the number of\dss simplices with a desirable property\dss
is odd and\dss hence non-zero.\oss
There is a natural\dss way\dss to assign\dss either\dss $1$\dss or\dss $-\dff 1$\dss
to\sss each of\dss these simplices,\oss
and\dss the sum of\dss these numbers\sss turns out\dss to be\dss $1$\nnsp.\oss
While\sss this result\dss immediately\dss follows\dss from\dss Alexander's\dss lemma,\oss
it\dss was published\dss in\dss 1961\dss as a\qss ``strengthening of\dss Sperner's\dss lemma''\pss
\cite{bc}.\oss 

Apparently,\oss initially\dss it\dss was well\dss understood\dss how\dss the lemmas of\qss
Alexander\dss and\dss Sperner\dss are related.\oss
The\dss book\qss \cite{ah}\qss by\qss P.\dss Alexandroff\qss and\qss H.\dss Hopf\halfff,\oss
referred\dss to by\qss J.\qss Dieu\-donn\'{e},\pss 
was comissioned\dss in\dss 1928\trs by\trs R.\dss Courant\dss for his book series\qss
\emph{Die\dss Grundlehren\dss der\dss Mathematischen\dss Wissen\-schaften},\pss
published\dss by\trs Springer-Verlag.\oss
It\dss was published\sss in\dss 1935.\oss
For quite a while it\dss was a definitive monograph\sss in\dss topology\halfff.\oss
Sperner's\trs lemma appears in\qss \cite{ah}\qss in an\dss Appendix\dss to Chapter\dss IX\dss
devoted\dss to\qss ``elementary''\qss proofs of\qss Brouwer's\dss fixed\dss point\dss theorem
and\dss related\dss results.\oss
The\sss proofs\sss given in\qss \cite{ah}\qss are not\sss quite elementary\halfff:\oss
they\dss are\dss based\sss on a version of\qss Alexander's\dss lemma,\pss
proved\dss earlier\dss in\dss the book\dss by\dss using\sss chains\sss and\dss Alexander's\dss methods.\oss
Sperner's\dss paper\qss \cite{s}\qss is\dss referred\dss to only\dss in a footnote,\oss
while\dss Alexander's\dss papers\qss \cite{a1}\halfff,\pss
\cite{a2}\qss  are listed at\dss the end\dss of\dss the book among\dss the main
references for\dss this\dss Chapter\halfff.\oss\vspace{2pt}

In\qss 1932\qss P.\dss Alexandroff\qss published a short\dss book\qss \cite{pa1},\pss 
a sort\dss of\dss popular\dss introduction\dss to\sss the\sss basic notions 
of\dss topology\halfff,\oss
which at\dss the same\sss time looks\sss like 
a\sss blueprint\dss for\dss parts of\qss \cite{ah}.\oss 
It\dss was originally\dss intended\dss to be an appendix\dss to\dss 
Hilbert's\qss \emph{Anschauliche\dss Geometrie}\oss
\cite{hc}\halfff,\pss
and\dss the preface\sss to\qss \cite{pa1}\qss 
was written\sss by\qss Hilbert\qss himself\halfff.\oss
The book culminates in an  
outline of\dss proofs of\dss 
the\sss topological\dss invariance of\dss dimension and\sss of\dss homology\dss groups
based\dss on\dss methods of\qss Brouwer\halfff,\pss
Lebesgue,\pss and\qss Alexander\halfff.\oss
As in\qss \cite{ah},\pss
Alexander's\dss lemma appears in the form of\dss the last\dss two out\sss of\trs three\qss
\emph{conservation\dss theorems}.\oss 
Sperner's\qss lemma is not\sss even\dss mentioned\halfff.\oss\vspace{2pt}

Later on\dss the fates of\qss Alexander's\dss and\dss Sperner's\dss lemmas diverged.\oss
Alexander's\dss lemma\dss disappeared\dss in\dss the 
sky of\dss more and\dss more abstract\sss and\dss powerful\dss
machinery of\dss algebraic topology\halfff.\oss
It\trs became customary\dss to prove\dss Brouwer's\dss fixed-point\dss theorem as
an illustration of\trs the power of\trs this machinery\halfff.\oss
For example,\pss Spanier\qss \cite{sp}\qss 
proves it\sss only\dss on\dss p.\dss 194\halfff,\oss
and\dss Hatcher\qss \cite{h}\qss 
states\qss Brouwer's\dss fixed-point\dss theorem\sss on\dss p.\qss 114\qss
and\dss completes\sss the proof\trs on\dss p.\dss 124.\oss
Sperner's\qss lemma\dss became a tool of\dss choice
in\dss more set-theoretic branches of\dss topology\sss
such as\sss the dimension theory,\oss
although\dss algebraic\dss topology\dss triumphantly\dss returned\dss 
to\sss the dimension\dss theory\dss in\dss the works\sss of\qss
A.N.\dss Dranishnikov\qss 
\cite{dr1},\pss \cite{dr2}.\oss
Sperner's\qss lemma\dss became a tool\sss of\dss choice also
in\sss combinatorics,\oss game\sss theory\halfff,\oss 
and\dss mathematical\sss economics.\oss\vspace{2pt}

The rest\sss of\trs the paper is devoted\dss to the mathematical details of\dss this story\halfff.\oss
No familiarity with algebraic topology is assumed,\pss
and,\qss perhaps,\qss the paper can serve as an invitation to it\halfff.\oss
We start\dss with\dss the basic notions of\dss simplices,\oss
simplicial\sss complexes,\oss and chains,\pss 
and\dss then\dss prove\dss Alexander's\dss lemma.\oss
As\dss the first\sss applications\sss we prove\dss
Brouwer's\trs invariance of\dss dimension and\sss invariance of\dss domains\sss theorems\sss
following\dss Lebesgue\dss ideas\sss in\dss the form given\dss to\sss them\dss by\dss Sperner\halfff.\oss
But\dss we refer\dss to\qss Alexander's\dss lemma\dss instead\sss of\qss Sperner's\dss one.\oss
Next\halfff,\pss we\dss introduce simplicial\sss approximations 
and use\sss them\dss to\dss prove\dss
Brouwer's\dss fixed-point\dss theorem.\oss
Our main application of\dss simplicial\sss approximations\dss is\dss to\sss the beautiful\qss
Alexander's\dss proof\dss
of\dss the\sss topological\dss invariance of\trs homology\dss groups.\oss
A\sss technical\dss part\sss of\trs this proof\trs is\dss relegated\dss to\dss
Appendix\qss \ref{stellar}.\oss
Brouwer's\dss fixed-point\dss theorem\dss is\dss also proved\dss by\qss
Knaster-Kuratowski-Mazurkiewich\qss argument\halfff,\pss
which after a closer\sss examination\dss turns out\dss to be a version of\dss
the proof\dss based on simplicial\sss approximations.\pss
Section\qss \ref{sperner-lemma}\qss is\dss devoted\dss to\dss Sperner's\dss lemma and\dss
its\dss combinatorial\dss proof\halfff,\pss
and\qss Sections\qss \ref{dualizing}\qss and\qss \ref{cohomology}\qss
to\sss their cohomological\dss interpretation.\oss
In\qss Section\qss \ref{graphs}\qss we explain\dss how\sss classical\dss proofs\sss
of\pss Sperner's\trs lemma\dss
lead\sss to\sss so-called\qss \emph{path-following\dss algorithms}.\oss

\mysection{Simplicial\qss complexes\qss and\qss chains}{simplicial-complexes}

\myuppar{Geometric simplices.}
Recall\dss that\sss the\trs \emph{convex\dss hull}\pss of\dss points\trs 
$x_{\dff 0}\fff,\pff x_{\fff 1}\fff,\pff \ldots\fff,\pff x_{\fff n}\off \in\off \rrr^{\fff d}
$\qss 
in a Euclidean space\dss $\rrr^{\fff d}$\dss is\dss 
the set of\dss all\dss linear combinations\vspace*{3.375pt}
\begin{equation}
\label{convex-hull}
\quad
\sum\nolimits_{i\qff =\qff 0}^{n}\qff a_{\fff i}\dff x_{\fff i}
\end{equation}

\vspace*{-8.625pt}
such\dss that\dss the coefficients\dss $a_{\fff i}$\dss 
are real\dss and\dss non-negative\dss and\vspace*{3pt}
\[
\quad
\sum\nolimits_{i\qff =\qff 0}^{n}\qff a_{\fff i}
\off\off =\off\off
1\dff.
\]

\vspace*{-9pt}
The points\trs
$x_{\dff 0}\fff,\pff x_{\fff 1}\fff,\pff \ldots\fff,\pff x_{\fff n}$\qss
are said\dss to be\qss \emph{affinely\dss independent}\pss
if\dss the two relations\vspace*{4.5pt}
\[
\quad
\sum\nolimits_{i\qff =\qff 1}^{n\qff +\qff 1}\qff a_{\fff i}\dff x_{\fff i}
\off\off =\off\off
0
\hspace*{1.5em}\mbox{and}\hspace*{1.5em}
\sum\nolimits_{i\qff =\qff 1}^{n\qff +\qff 1}\qff a_{\fff i}
\off\off =\off\off
0
\]

\vspace*{-7.5pt}
together imply that all coefficients\qss $a_{\fff i}\qff =\qff 0$\nnsp.\oss
If\dss this is the case,\pss then the presentation of\dss a point of\dss the convex hull of\qss
$x_{\dff 0}\fff,\pff x_{\fff 1}\fff,\pff \ldots\fff,\pff x_{\fff n}$\qss
in the form\qss (\ref{convex-hull})\qss is unique.\oss
In\dss this case\dss the numbers\qss
$a_{\dff 0}\fff,\pff a_{\fff 1}\fff,\pff \ldots\fff,\pff a_{\fff n}$\qss
are called\dss the\qss
\emph{barycentric coordinates}\qss
of\dss the point\qss (\ref{convex-hull}).\oss
A trivial verification shows that the points\qss
$x_{\dff 0}\fff,\pff x_{\fff 1}\fff,\pff \ldots\fff,\pff x_{\fff n}$\qss
are affinely independent\sss if\dss and only\sss if\dss the vectors\vspace*{3pt}
\[
\quad
x_{\fff 1}\qff -\qff x_{\dff 0}\fff,\off
x_{\fff 2}\qff -\qff x_{\dff 0}\fff,\off
\ldots\fff,\off
x_{\fff n}\qff -\qff x_{\dff 0}
\]

\vspace*{-9pt}
are linearly independent\sss in\dss $\rrr^{\fff d}$\dnsp.\oss
A\qss \emph{geometric $n$\dnsp-simplex}\halfff,\pss
or\sss a\qss \emph{geometric\dss simplex\dss of\qss dimension} $n$ in\dss $\rrr^{\fff d}$\dss
is defined as the convex hull of\dss $n\qff +\qff 1$ affinely independent points in $\rrr^{\fff d}$\nsp\dnsp,\oss
called\dss its\qss \emph{vertices}.\oss
A\qss \emph{geometric\dss simplex}\qss is defined as a geometric $n$\dnsp-simplex for some $n$\nnsp.\oss
The convex hulls of\dss subsets of\dss the set of\dss vertices of\dss
a geometric simplex\dss $\sigma$\dss are called\dss its\qss \emph{faces}.\oss
A face of\dss $\sigma$\dss is said\dss to be\qss \emph{proper}\pss 
if\dss it is not equal\dss to\dss $\sigma$\dnsp.\oss
Each face is also a geometric simplex\halfff.\oss
A face\sss is said\dss to be an\dss \emph{$n$\dnsp-face}\pss if\trs it\sss is\sss a
geometric $n$\dnsp-simplex\halfff.\oss
The union\dss $\bd \sigma$\dss of\dss all\dss proper\dss faces of\dss $\sigma$\dss
is called\dss the\qss \emph{geometric boundary}\qss of\dss $\sigma$\nnsp.\oss
It\dss is\dss the boundary of\dss $\sigma$\dss in\dss the most\dss naive sense.\oss

\myuppar{Geometric\dss simplicial\dss complexes.}
A\qss \emph{geometric simplicial complex}\qss in $\rrr^{\fff d}$
is a finite collection $S$ of\dss geometric simplices in\dss $\rrr^{\fff d}$\dss
such that\qss if\pss $\sigma\qff \in\qff S$\qss
and $\tau$ is a face of $\sigma$\nnsp,\oss 
then\qss $\tau\qff \in\qff S$\dnsp,\oss
and\dss if\qss $\sigma\fff,\pff \sigma'\qff \in\qff S$\qss
then\qss $\sigma\qff \cap\qff \sigma'$\qss is a face of\dss both
$\sigma$ and $\sigma'$\dnsp.\oss
The vertices of\dss geometric simplices\qss $\sigma\qff \in\qff S$\qss
are called\qss \emph{vertices}\qss of\dss $S$\nnsp,\oss
and\dss the set\dss of\dss vertices of\dss $S$\dss 
is denoted\dss by\dss $v\dff(\dff S\dff)$\nnsp.\oss
The\qss \emph{dimension}\qss of\dss $S$\dss is\sss the maximal $n$ such that\dss $S$\dss
contains an $n$\dnsp-simplex\halfff.\oss
A geometric simplicial complex\dss $Q$\dss is said\dss to be a\qss \emph{subcomplex}\qss of\dss $S$\dss
if\dss every simplex of\dss $Q$\dss is also a simplex of\dss $S$\nnsp.\oss 

The union of\dss all geometric simplices of\dss $S$\dss is denoted\dss by\sss
$\norm{S}$\sss and\sss called\dss the\qss \emph{polyhedron}\pss of\qss $S$\nnsp,\oss
and\dss
$S$\dss is\dss said\dss to be
a\qss \emph{triangulation}\qss of\dss $\norm{S}$\dnsp.\oss
If\dss $Q$\dss is a subcomplex of\dss $S$\nnsp,\oss
then\qss $\norm{Q}\qff \subset\qff \norm{S}$\nnsp.\oss
Most\sss of\dss interesting\qss
(finitely dimensional)\qss topological spaces are homeomorphic\sss to polyhedra,\oss and
such spaces 
are\sss 
most\sss accessible to combinatorial\sss and\sss algebraic methods.\oss

\myuppar{Abstract\dss simplicial\dss complexes.}
A substantial\dss part\sss of\dss the theory of\dss geometric simplicial\sss
complexes is\dss
purely combinatorial.\oss
There is a combinatorial counterpart\sss of\dss the notion
of\dss a geometric simplicial complex,\oss
namely,\oss the notion of\dss
an\qss \emph{abstract\dss simplicial complex}.\oss
It\dss is defined as
a finite collection $K$ of\dss subsets of\dss
a finite set\dss $v\dff(\dff K\dff)$\dss such that
if\qss
$\sigma\qff \in\qff K$\qss and\qss
$\sigma'\qff \subset\qss \sigma$\dnsp,\oss
then\qss
$\sigma'\qff \in\qff K$\nnsp,\oss
and\dss $v\dff(\dff K\dff)$\dss  is equal\dss to the union of\dss all subsets in $K$\dnsp.\oss
The elements of\dss $v\dff(\dff K\dff)$\dss  are called\dss the\qss \emph{vertices},\oss
and\dss the elements of\dss $K$ the\qss \emph{simplices}\qss of\dss
$K$\nnsp.\oss

A simplex $\sigma'$ is said\dss to be a\qss \emph{face}\qss of\dss
a simplex $\sigma$ if\qss $\sigma'\qff \subset\qff \sigma$\nnsp.\oss
A face of\dss $\sigma$\dss is said\dss to be\qss \emph{proper}\pss 
if\qss it\dss is\sss not\sss equal\dss to\dss $\sigma$\nnsp.\oss
A simplex\qss $\sigma\qff \in\qff K$\qss is said\dss to be an\dss
\emph{$n$\dnsp-simplex}\oss
or\dss a\qss \emph{simplex\dss of\dss dimension}\qss $n$\qss 
if\pss $\num{\sigma}\off =\off n\qff +\qff 1$\nnsp,\oss
where,\pss as usual\halfff,\oss we denote by\dss $\num{\sigma}$\dss
the number of\dss elements of\dss $\sigma$\nnsp.\oss
The\qss \emph{dimension}\qss of\dss $K$\dss is\dss the maximal $n$ such that\dss $K$\dss
contains an $n$\dnsp-simplex\halfff.\oss\vspace{0pt}

\myuppar{Geometric\sss simplicial\sss complexes and\sss abstract\sss ones.}
A geometric simplicial complex $S$ leads to an abstract 
simplicial complex $a\dff (\dff S\dff)$
having as its set of\dss vertices the set\dss $v\dff(\dff S\dff)$\dss of\dss vertices of\dss $S$\nnsp,\oss 
and as its $n$\dnsp-simplices the sets of\dss vertices of\dss geometric 
$n$\dnsp-simplices of\dss $S$\nnsp.\oss
The geometric simplicial complex $S$ can be recovered from $a\dff (\dff S\dff)$
as the set of\dss convex hulls of\dss simplices of\dss $a\dff (\dff S\dff)$\nnsp.\oss
But\dss the fact\sss that\sss the vertices of\dss $a\dff (\dff S\dff)$\dss are
the points of\dss $\rrr^{\fff d}$\dss is\sss better to be ignored\dss
to the extent\dss possible.\oss 
From such a point of\dss view
the complex $a\dff (\dff S\dff)$ 
encodes the combinatorics of\dss
$S$\dnsp,\oss 
understood as the pattern of\dss intersections 
of\dss geometric simplices of\dss $S$\dnsp.\oss

The combinatorial\dss part of\dss the theory deals\dss
not\sss with\dss geometric simplicial complexes\dss $S$\nnsp,\oss
but\dss with corresponding abstract\sss simplicial complexes\dss
$a\dff (\dff S\dff)$\nnsp.\oss
Since\dss there is a\sss tautological one-to-one correspondence
between\dss the simplices of\dss $S$\dss and\dss the simplices of\dss
$a\dff (\dff S\dff)$\nnsp,\oss
which respects\sss the property of\dss being a face,\oss
usually\dss there is no need\dss to distinguish\dss between\dss $S$\dss
and\dss $a\dff (\dff S\dff)$\nnsp.\oss
Also,\oss some definitions and arguments apply equally well\dss to
both geometric and abstract\sss simplicial complexes.\oss 
In such situations we will speak simply about\qss \emph{simplicial\sss complexes}.\oss\vspace{0pt}

\myuppar{Simplicial\dss maps.}
This notion is easier\dss to introduce in\dss the context\sss of\dss
abstract\dss complexes.\oss
Let\pss 
$K\dff,\off L$\qss be abstract\sss simplicial complexes.\oss
\emph{Simplicial\dss maps}\oss 
$\varphi\dff \colon\dff K\ttoo L$\oss
are defined as maps\oss\vspace{3pt}
\[
\quad
\varphi\dff \colon\dff v\dff(\dff K\dff)\qff \ttoo\qff v\dff(\dff L\dff)
\]

\vspace{-9pt}
taking simplices of\dss $K$\dss to simplices of\dss $L$\nnsp.\oss
When\dss $L$\dss consists of\dss a single simplex and\dss its faces,\oss
every subset of\dss $v\dff(\dff L\dff)$\dss is a simplex and every map\qss
$v\dff(\dff K\dff)\toto v\dff(\dff L\dff)$\qss is a simplicial\dss map\qss
$K\toto L$\nnsp.

For\dss 
geometric simplicial complexes\pss
$S\fff,\off S'$\pss
\emph{simplicial\dss maps}\oss 
$S\ttoo S'$\oss
are defined\sss simply\sss as simplicial\sss maps\qss
$a\dff(\dff S\dff)\ttoo a\dff(\dff S'\dff)$\nnsp.\oss
In other words,\oss a simplicial\dss map\qss
$\varphi\dff \colon\dff
S\ttoo S'$\oss
is a map\qss 
$\varphi\dff \colon\dff
v\dff(\dff S\dff)\ttoo v\dff(\dff S'\dff)$\dss 
such\dss that\dss $\varphi$\dss takes\sss the set\sss of\dss vertices of\dss
each simplex of\dss $S$\dss into\sss the set\sss of\dss vertices of\dss
some simplex of\dss $S'$\nnsp.\oss
Obviously,\oss a simplicial\dss map\qss 
$\varphi
\dff \colon\dff
S\ttoo S'$\qss 
defines a map from\dss
the set\sss of\dss simplices of\dss $S$\dss to\sss the set\sss of\dss
simplices of\dss $S'$\nnsp.\oss
The latter map\sss is\sss also denoted\dss by\dss $\varphi$\nnsp.\oss

While we treat\dss simplicial\dss maps as combinatorial\dss
objects,\oss their\qss \emph{raison d'\^{e}tre}\qss
is\sss the fact\dss that\dss they\sss are combinatorial\sss analogues
of\dss continuous maps.\oss
It\dss is comforting\dss to know\dss that\sss
a simplicial\dss map\qss
$\varphi
\dff \colon\dff
S\ttoo S'$\qss
canonically\sss extends\sss to a
continuous map\qss
$\norm{\nnsp\varphi\nnsp}\dff \colon\dff
\norm{S}\ttoo \norm{S'\nsp}$\nnsp.\oss
It\dss is defined as\sss follows.\oss
Let\qss
$\{\qff x_{\fff 0}\fff,\pff x_{\fff 1}\fff,\pff \ldots\fff,\pff x_{\fff n} \qff\}$\pss
be\sss the set\sss of\dss vertices of\dss a simplex of\dss $S$\nnsp.\oss 
If\qss
$a_{\fff 0}\fff,\pff a_{\fff 1}\fff,\pff \ldots\fff,\pff a_{\fff n}$\qss
are\dss non-negative numbers with\dss the sum\dss $1$\nnsp,\oss
then\vspace{4.5pt}
\[
\quad
\norm{\nnsp\varphi\nnsp}
\qff \colon\off
\sum\nolimits_{i\qff =\qff 0}^{n}\qff a_{\fff i}\dff x_{\fff i}
\off\off \longmapsto\off\off 
\sum\nolimits_{i\qff =\qff 0}^{n}\qff a_{\fff i}\trf \varphi\dff(\dff x_{\fff i} \dff)
\qff.
\]

\vspace{-7.5pt}
It\dss is easy\dss to see\sss that\dss $\norm{\nnsp\varphi\nnsp}$\dss
is correctly\sss defined\sss and\dss continuous.\oss
Somewhat\sss surprisingly,\oss we will\sss use\sss not\dss the maps\dss
$\norm{\nnsp\varphi\nnsp}$\nnsp,\oss
but\dss simplicial\sss maps\dss $\varphi$\dss themselves as\qss \emph{models}\qss of\dss continuous maps.

\myuppar{Subdivisions of\dss geometric simplicial complexes.}
In order\dss to find a good enough simplicial\dss model of\dss 
$\norm{S}$\dss considered as a topological\sss space or\sss of\dss a continuous map\qss 
$\norm{S}\ttoo \norm{T}$\nnsp,\oss
one usually needs to replace $S$ by\sss a 
simplicial complex having\dss the same polyhedron,\oss but\sss smaller simplices.\oss
A geometric simplicial complex\dss $S'$\dss is said\dss to be a\qss \emph{subdivision}\qss
of\dss a geometric simplicial complex\dss $S$\dss if\dss
every simplex of\dss $S'$\dss is contained\sss in a simplex of\dss $S$\dss
and every\qss $\sigma\qff \in \qff S$\qss is equal\dss to the union of\dss 
simplices of\trs $S'$\dss contained\sss in\dss $\sigma$\dnsp.\oss
If\qss $S'$\dss is a subdivision of\dss $S$\dnsp,\oss 
then,\oss obviously\halfff,\pss
$\norm{S'}\off =\off \norm{S}$\nnsp.\oss
Given\dss $S'$\dss and\qss $\sigma\qff \in\qff S$\dnsp,\oss
let\dss $S'\fff(\dff \sigma\dff)$\dss be the set of\dss simplices of\dss $S'$\dss
contained\sss in\dss $\sigma$\dnsp.\oss
Then\dss $S'\fff(\dff \sigma\dff)$\dss is a geometric simplicial complex and\pss
$\norm{S'\fff(\dff \sigma\dff)}\off =\off \sigma$\dnsp,\qff\oss
i.e.\qss $S'\fff(\dff \sigma\dff)$\sss is a triangulation of\dss $\sigma$\dnsp.\oss
Clearly,\oss the polyhedron\dss $\norm{S}$\dss is
equal\dss to\sss the union of\dss polyhedra\dss $\norm{S'\fff(\dff \sigma\dff)}$\nnsp.

\myuppar{Chains and\dss boundaries.}
The key element of\dss the combinatorial\dss
structure of\dss simplicial complexes 
is the relation\oss 
\emph{``\nsp$\tau$\qss is\dss a\dss face\dss of\qss $\sigma$\dnsp''}\oss
between\dss two simplices\dss $\tau\fff,\pff \sigma$\nnsp.\oss
The boundary of\dss a geometric figure 
has dimension\dss less by $1$ than\dss 
the dimension of\dss the geometric figure\dss itself\halfff,\oss
and\dss 
the geometric intuition\sss
suggests\sss to concentrate on\dss the case
when\dss the dimension of\dss $\tau$\dss is\dss $1$\dss less\sss
than\dss the dimension of\dss $\sigma$\nnsp.\oss
In\dss this case\sss algebraic\dss topology\sss 
suggests\sss to encode\sss this relation\dss
by a map assigning to an $n$\dnsp-simplex\sss 
the formal sum of\dss all\dss its $(\fff n\dff -\dff 1 \fff)$\dnsp-faces.\oss
Such formal sums with coefficients in a fixed abelian group are known as\qss \emph{chains}.\oss
For\dss the purposes of\dss this paper\sss it\dss is sufficient\dss
to\sss consider only\dss the chains with coefficients in\qss
$\ftwo\off =\off \zzz\dff/2\dff \zzz$\nnsp.\oss
Such chains can be identified with sets of\dss simplices of\dss
the same dimensions,\pss namely,\oss
with sets of\dss simplices appearing with\dss non-zero coefficients in\dss the formal sum,\oss 
and\sss can\dss be\sss thought\sss of\dss 
as geometric figures,\oss say,\oss 
as complexes consisting of\dss all\dss faces of\dss simplices with\dss non-zero coefficients.\oss 

Let\dss us\dss turn\dss to formal definitions and consider
a simplicial complex\dss $X$\nnsp,\oss
either geometric or abstract\halfff.\oss
For an integer\qss $m\qff \geq\qff 0$\qss
let\dss $C_{\fff m}\fff(\dff X \dff)$\dss 
be the vector space over\dss $\ftwo$\dss
having\dss the set of\dss all\dss $m$\dnsp-simplices of\dss $X$\dss
as its basis.\oss
The elements of\dss $C_{\fff m}\fff(\dff X \dff)$\dss
are called\dss the\dss \emph{$m$\dnsp-chains}\qss of\dss $X$\nnsp.\oss
The\qss \emph{boundary}\qss $\partial\fff\sigma$\dss 
of\dss an $m$\dnsp-simplex $\sigma$ of\dss $K$\dss 
is defined as the sum of\dss all $(\fff m\dff -\dff 1 \fff)$\dnsp-faces\dss
of\dss $\sigma$\dnsp.\fff\oss
Extending\dss the map\qss
$\sigma\off \longmapsto\off \partial\fff\sigma$\qss
by\dss linearity we get\dss the\qss \emph{boundary\dss operator}\vspace{3pt}
\[
\quad
\partial
\qff \colon\qff
C_{\fff m}\fff(\dff X \dff)
\off \ttoo\off
C_{\fff m\dff -\dff 1}\fff(\dff X \dff)\dff.
\]

\vspace{-9pt}
If\dss $X$\dss is a geometric simplicial complex,\oss
then\dss for every\dss $m$\dss there is a canonical\dss isomorphism\dss
between\dss the spaces of\dss $m$\dnsp-chains of\dss $X$\dss
and\sss of\dss $a\dff(\dff X\dff)$\nnsp.\oss
These isomorphisms respect\dss the\sss boundary operators,\oss
and we will use\sss them\dss to identify chains of\dss $X$\dss
with chains of\dss $a\dff(\dff X\dff)$\nnsp.\oss

\myuppar{Induced\sss maps.}
Let\pss $X\fff,\pff Y$\pss 
be two simplicial complexes\qss
(either geometric or abstract)\qss 
and\dss let\qss 
$\varphi\dff \colon\dff X\ttoo Y$\oss
be a simplicial\dss map.\oss
For\sss an integer\qss $m\qff \geq\qff 0$\qff and\sss
an $m$\dnsp-simplex $\sigma$ of\qss $X$\qss let\vspace{3pt}
\[
\quad
\varphi_*\dff(\dff\sigma\dff)
\off =\off
\varphi\dff(\dff\sigma\dff)
\]

\vspace*{-12pt}
if\qss $\varphi\dff(\dff\sigma\dff)$\dss
is\trs an\qss $m$-simplex\fff.\oss 
Otherwise\dss $\varphi\dff(\dff\sigma\dff)$\dss is a simplex of\dss
dimension\qss $<\qff m$\qss and we set \vspace*{0pt}
\[
\quad
\varphi_*\dff(\dff\sigma\dff)
\off =\off
0
\qff.
\]

\vspace{-9pt}
Informally,\oss if\qss the dimension of\qss $\varphi\dff(\dff\sigma\dff)$\dss
is\qss $<\qff m$\nnsp,\oss then\dss $\varphi\dff(\dff\sigma\dff)$\dss
is equal\dss to zero\qss \emph{as an $m$\dnsp-simplex}.\oss
Let\dss us define\sss the\qss \emph{induced\dss map}\qss\vspace*{1.5pt}
\[
\quad
\varphi_*
\qff \colon\qff
C_{\fff m}\fff(\dff X \dff)
\off \ttoo\off
C_{\fff m}\fff(\dff Y \dff)\dff.
\]

\vspace*{-10.5pt}
as\dss 
the extension of\dss the map\qss
$\sigma\qff \longmapsto\qff \varphi_*(\fff\sigma\fff)$\oss
by\dss linearity.\oss
The basic property of\dss induced\dss maps is\dss the fact\dss
that\dss they commute with\dss the boundary operators in\dss
the sense of\dss the following\dss theorem.

\myuppar{Theorem\qss 1.}
$\partial\qff \circ\qff \varphi_*
\off =\off
\varphi_*\qff \circ\qff \partial$\nsp,\oss
\emph{i.e.}
\vspace{3pt} 
\begin{equation}
\label{induced-maps}
\quad
\partial\trf\bigl(\dff \varphi_*\dff(\dff \sigma \dff) \dff\bigr)
\off =\off
\varphi_*\dff\bigl(\dff \partial\trf(\dff \sigma \dff) \dff\bigr)
\end{equation}

\vspace{-9pt}
\emph{for\dss all\sss 
simplices $\sigma$ of\pss $X$\nnsp.}\oss

\proof
Let\dss $\sigma$\dss be an $m$\dnsp-simplex\dss of\qss $X$\nnsp.\oss
Then\dss $\varphi\dff(\dff \sigma \dff)$\dss is a simplex of\dss $Y$\dss
of\dss dimension\qss $\leq\qff m$\nnsp.

\myepar{Case\qss 1.\oss The\dss dimension\dss of\qss $\varphi\dff(\dff \sigma \dff)$\dss 
is equal\dss to\dss $m$\nnsp.}\oss
In\dss this case\dss
$\varphi$\dss is injective on $\sigma$ and\dss hence\dss 
$\varphi\dff(\dff\tau\dff)$\dss 
is\dss an $(\fff m\dff -\dff 1 \fff)$\dnsp-simplex\dss
for\dss all\dss $(\fff m\dff -\dff 1 \fff)$\dnsp-faces $\tau$ of\dss $\sigma$\dnsp.\oss
This implies\qss (\ref{induced-maps})\qss for such $\sigma$\dnsp.\oss

\myepar{Case\qss 2.\oss The\dss dimension\dss of\qss $\varphi\dff(\dff \sigma \dff)$\dss 
is\qss $\leq\qff m\qff -\qff 2$\nnsp.}\oss
In\dss this case\dss
$\varphi_*\dff(\dff\sigma\dff)
\off =\off
0$\oss 
by\dss the definition.\oss 
If\dss $\tau$\dss is an $(\fff m\dff -\dff 1\fff)$\dnsp-face of\dss $\sigma$\nnsp,\oss
then\dss 
$\varphi_*\dff(\dff\tau\dff)$\dss
is contained\dss in\dss $\varphi_*\dff(\dff\sigma\dff)$\dss
and\dss hence\sss the dimension of\trs
$\varphi_*\dff(\dff\tau\dff)$\dss is\qss $\leq\qff m\qff -\qff 2$\nnsp.\oss
In\dss particular\halfff,\oss
$\varphi\fff(\fff\tau\fff)$\dss is not\dss an 
$(\fff m\dff -\dff 1 \fff)$\dnsp-simplex and\dss hence\qss
$\varphi_*\dff(\dff\tau\dff)\off =\off 0$\nnsp.\oss
Therefore,\oss
in\dss this case the both sides of\qss (\ref{induced-maps})\qss 
are equal\dss to $0$\nnsp.\qff\oss

\myepar{Case\qss 3.\oss The\dss dimension\dss of\qss $\varphi\dff(\dff \sigma \dff)$\dss 
is equal\dss to\qss $m\qff -\qff 1$\nnsp.}\oss
In\dss this case\qss $\num{\varphi\dff(\dff \sigma\dff)}\off =\off m$\nnsp,\oss
and\dss since\qss $\num{\sigma}\off =\off m\qff +\qff 1$\nnsp,\oss
there is a unique pair\qss $a\fff,\pff b\qff \in\qff \sigma$\qss
such that\pss
$\varphi\dff(\dff a\dff)
\off =\off
\varphi\dff(\dff b\dff)$\pss
and\oss
$a\off \neq\off b$\nnsp.\qff\oss
This implies\sss that\qss
$\varphi\dff(\dff\sigma\qff \smallsetminus\qff \{\dff a \qff\}\dff)
\off =\off
\varphi\dff(\dff\sigma\qff \smallsetminus\qff \{\dff b \qff\}\dff)
\off =\off
\varphi\dff(\dff\sigma\dff)$\qss
and\qss 
$\num{\varphi\dff(\dff \tau \dff)}\qff \leq\qff m\qff -\qff 1$\pss
if\qss $\tau$\trs is\dss an $(\fff m\dff -\dff 1 \fff)$\dnsp-face\dss
of\dss $\sigma$ different\trs from\qss
$\sigma\qff \smallsetminus\qff \{\dff a \qff\}$\nnsp,\pss
$\sigma\qff \smallsetminus\qff \{\dff b \qff\}$\nnsp.\qff\oss
Therefore\vspace{3pt}
\[
\quad
\varphi_*\dff\bigl(\trf \partial\dff(\dff\sigma\dff) \dff\bigr)
\off =\off
\varphi\dff(\dff\sigma\qff \smallsetminus\qff \{\dff a \qff\}\dff)
\off +\off
\varphi\dff(\dff\sigma\qff \smallsetminus\qff \{\dff b \qff\}\dff)
\off =\off
2\dff \varphi\dff(\dff\sigma\dff)
\qff.
\]

\vspace{-9pt}
Since we work over\dss $\ftwo$\nnsp,\oss
it follows that\pss
$\varphi_*\dff(\trf \partial\dff(\dff\sigma\dff) \dff)
\off =\off
0$\dnsp.\oss
Also,\pss
$\varphi_*\dff(\dff \sigma\dff)\off =\off 0$\qss
because\dss the dimension of\qss $\varphi_*\dff(\dff \sigma\dff)$\dss
is\qss $m\qff -\qff 1$\nnsp.\oss 
Therefore,\oss both sides of\qss (\ref{induced-maps})\qss 
are equal\dss to\dss $0$\nnsp.\oss  \eproof

\mysection{A\qss lemma\qss of\pss Alexander}{alexander-lemma}

\myuppar{Subdivisions of\dss a geometric simplex\halfff.}
Let\qss
$I
\off =\off
\{\qff 0\fff,\pff 1\fff,\pff 2\fff,\pff \ldots\fff,\pff n\qff\}$\nnsp.\oss
Let $\delta$ be the geometric $n$\dnsp-simplex\dss in $\rrr^{\fff d}$ with the vertices\qss
$v_{\fff 0}\fff,\pff v_{\fff 1}\fff,\pff \ldots\fff,\pff v_{\fff n}$\nnsp.\oss
For each\qss $i\qff \in\qff I$\qss
let\dss $\delta_{\fff i}$\dss be the 
$(\fff n\dff -\dff 1 \fff)$\dnsp-face of\dss $\Delta$\dss
having as its vertices all\dss points 
$v_{\fff 0}\fff,\pff \ldots\fff,\pff v_{\fff n}$ 
except $v_{\fff i}$\nnsp.\oss 
Every proper\dss face of\dss $\delta$\dss is contained\dss in\dss
$\delta_{\fff i}$\dss for some\dss $i\qff \in\qff I$\nnsp.\oss
The geometric boundary\dss $\bd\fff \delta$\dss is equal\dss to\sss
the union of\dss all\dss faces\dss $\delta_{\fff i}$\nnsp.\oss
Let\dss $\Delta$\dss be\sss the geometric simplicial complex consisting of\dss
the simplex\dss $\delta$\dss and\sss all\dss its faces,\oss
and\dss let\dss $\bd \Delta$\dss be\sss the complex
consisting of\dss all\dss proper faces of\dss $\delta$\nnsp.\oss
Then\dss $\delta$\dss is an $n$\dnsp-simplex of\dss $\Delta$\dss and\vspace{3pt}
\begin{equation}
\label{boundary-delta}
\quad
\partial\dff \delta
\off =\qff\off
\sum\nolimits_{i\qff =\qff 0}^n\qff \delta_{\fff i}
\off.
\end{equation}

\vspace{-9pt}
Let\dss $T$\dss be a subdivision of\dss $\Delta$\nnsp,\oss
or\halfff,\oss
what\dss is\dss the same,\oss
a\sss triangulation of\dss the simplex\dss $\delta$\nnsp.\oss

Clearly\halfff,\oss the dimension of\trs $T$\sss is $n$\nnsp.\oss
Since simplices are convex,\oss
every simplex of\trs $T$\sss contained\dss in $\bd\fff \delta$ 
is contained\dss in $\delta_{\fff i}$ for some\qss $i\qff \in\qff I$\nnsp.\oss
Let\dss $\bd\fff T$\sss be the subcomplex of\dss $T$\dss
consisting of\dss 
simplices  
contained\dss in $\bd\fff \delta$\nnsp,\oss
and\dss for each\qss
$i\qff \in\qff I$\pss
let\dss $T_i$\sss be the subcomplex of\dss $T$\dss consisting of\dss 
simplices  
contained\dss in $\delta_{\fff i}$\nnsp.\oss
Then\dss $\bd\fff T$\dss is a triangulation of\dss  $\bd\fff \delta$\dss
and\dss
$T_i$\dss is a triangulation of\dss $\delta_{\fff i}$\dss
for each\qss $i\qff \in\qff I$\nnsp.\oss
Obviously,\oss $\bd\fff T$\sss
is equal\dss to the union of\dss the complexes\dss $T_i$\nnsp.\vspace{6pt}

\emph{For\dss the rest\sss of\dss the paper we will\dss keep\sss the above notations\qss 
$I\fff,\pff \delta\fff,\pff \delta_{\fff i}\fff,\pff \Delta\fff,\pff T$\nnsp,\oss etc.}\oss\vspace{0pt}

\myuppar{Theorem\qss 2.}
\emph{Let\qss $Q$\dss be a geometric simplicial complex.\oss 
Suppose\sss that\oss
$\varphi\dff \colon\dff Q\qff \ttoo \Delta$\oss
is\dss a simplicial\dss map and\pss $\alpha\qff \in\qff C_{\dff n}\dff(\dff Q\dff)$\nnsp.\qff\oss
If\qff\oss
$\varphi_*\dff(\dff \partial\dff \alpha \dff)
\off =\off
\partial\fff \delta$\nnsp,\qff\oss
then\oss
$\varphi_*\dff(\dff \alpha \dff)
\off =\off
\delta$\nnsp.}

\proof
By\dss the definition,\pss $C_{\dff n}\dff(\dff \Delta\dff)$\dss
is a one-dimensional\sss vector space over\dss $\ftwo$\dss
with\dss $\delta$\dss forming a basis.\oss
Therefore\qss 
$\varphi_*\dff(\dff \alpha \dff)
\off =\off
c\dff \delta$\qss
for\dss some coefficient\dss $c$\qss
(of\dss course,\pss $c\off =\off 0$\dss or\dss $1$\nnsp,\oss 
but\dss the proof\dss does not\sss depend on\dss this).\oss
It\dss follows\sss that\qss
$\partial\dff (\dff \varphi_*\dff(\dff \alpha \dff) \dff)
\off =\off
c\trf (\dff \partial\fff \delta \dff)$\nnsp.\oss
But\dss
Theorem\qss 1\qss implies\sss that\qss
$\partial\dff (\dff \varphi_*\dff(\dff \alpha \dff) \dff)
\off =\off\qff
\varphi_*\dff(\dff \partial\dff \alpha \dff)
\off =\off\dff
\partial\fff \delta$\nnsp.\oss
Hence\qss $c\off =\off 1$\qss
and\qss
$\varphi_*\dff(\dff \alpha \dff)
\off =\off
\delta$\nnsp.\oss  \eproof

\myuppar{Subdivision of\dss chains.}
Let\dss $S$\dss be a  geometric simplicial complex and\dss let\dss $S'$\dss
be a subdivision of\dss a $S$\nnsp.\oss
Given a geometric $m$\dnsp-simplex $\sigma$ of\dss $S$\nnsp,\qff\oss
let\dss\vspace{3.5pt} 
\[
\quad
\fclass{\sigma}
\off =\off
\sum\nolimits_{\qff \sigma'\qff \in\pff S'\fff(\dff \sigma\dff)}\qff \sigma'
\off.
\]

\vspace{-8.5pt}
i.e.\qss let\dss $\fclass{\sigma}$\dss
be\dss the sum of\dss all\dss
geometric $m$\dnsp-simplices\dss $\sigma'$\dss of\qss $S'$\dss
contained\dss in\dss $\sigma$\nnsp.\oss
Then\qss 
$\fclass{\sigma}$\qss is an $m$\dnsp-chain of\dss $S'$\nnsp,\oss
called\dss the\qss \emph{subdivision of\pss $\sigma$\sss with\sss respect\dss to\dss $S'$\nnsp.}\oss
Extending\dss the map\qss
$\sigma\off \longmapsto\off \fclass{\sigma}$\oss 
by\dss linearity\dss leads\dss to a map\vspace{0pt}
\[
\quad
C_{\dff m}\dff (\dff S\dff)
\qff \ttoo\qff
C_{\dff m}\dff (\dff  S'\trf)
\qff
\]

\vspace{-12pt}
denoted\dss by\qss
$\alpha\qff \longmapsto\qff \fclass{\alpha}$\nnsp.\oss
The chain\dss $\fclass{\alpha}$\dss is called\dss the\qss
\emph{subdivision of\qss $\alpha$\dss with respect\dss to\dss $S'$\nnsp.}\oss

\myuppar{The non-branching\dss property\halfff.}
If\qss $\tau$\dss is an $(\fff n\dff -\dff 1 \fff)$\dnsp-simplex of\dss
$T$\dnsp,\oss
then either\pss
$\tau\pff \subset\pff \bd\fff \delta$\pss
and\dss then\pss
$\tau\pff \subset\pff \delta_{\fff i}$\pss for some $i$\nnsp,\pss
or\qss $\tau\pff \not\subset\pff \bd\fff \delta$\nnsp.\oss
In the first case $\tau$ is a face of\dss
exactly one geometric $n$\dnsp-simplex of\dss $T$\dnsp,\pss
and\dss in the second case $\tau$ is a face of\dss
exactly two geometric $n$\dnsp-simplices of\dss $T$\dnsp.\oss
We will accept\dss this property as geometrically obvious and call\dss it\dss
the\qss \emph{non-branching\dss property}.\oss    
Let\sss us consider\dss the subdivisions of\dss chains of\dss $\Delta$\dss
with\sss respect\dss to\dss $T$\nnsp.\oss
The non-branching\dss property
means\sss that\sss an
$(\fff n\dff -\dff 1\fff)$\dnsp-simplex of\dss $T$\dss
enters\sss the boundary\dss $\partial\dff \fclass{\delta}$\dss
with\dss the coefficient $1$ if\dss it\dss is contained\dss in\dss
$\bd\fff \delta$\nnsp,\oss
and\sss with\dss the coefficient $2$ otherwise.\oss
It\dss follows\sss that\vspace{4.5pt}
\begin{equation}
\label{boundary-fclass}
\quad
\partial\qff  
\fclass{\delta} 
\off =\qff\off
\sum\nolimits_{i\qff =\qff 0}^n\qff \fclass{\delta_{\fff i}}
\off.
\end{equation}

\vspace{-7.5pt}
\myuppar{Lemma.}\oss
$\partial\dff \fclass{\delta}
\off =\off
\fclass{\dff\partial\dff \delta\dff}$\nnsp.\oss

\proof
Since\dss the subdivision map\dss is linear\dss by\dss the definition,\oss\vspace{4.5pt}
\[
\quad
\sum\nolimits_{i\qff =\qff 0}^n\qff \fclass{\delta_{\fff i}}
\off =\off\qff
\bigfclass{\sum\nolimits_{i\qff =\qff 0}^n\qff \delta_{\fff i}}
\off.
\]

\vspace{-7.5pt}
It\dss remains\sss to combine\sss this equality\dss 
with\qss (\ref{boundary-delta}).\oss  \eproof

\myuppar{Corollary\halfff.}
\emph{If\oss  
$S'$\dss is a subdivision of\qss $S$\dss as above,\oss then}\oss
$\partial\dff \fclass{\alpha}
\off =\off
\fclass{\dff\partial\dff \alpha\dff}$\nnsp.\oss

\proof
Lemma\sss implies\sss that\dss this\sss is\sss true\sss when\dss  $\alpha$\dss
consists of\dss only\sss one simplex of\trs $S$\dnsp.\oss
The general case\sss follows\sss by\dss linearity.\oss  \eproof

\myuppar{Pseudo-identical\dss maps.}
If\qss $S'$\dss is\dss a\sss subdivision of\dss 
a geometric simplicial complex\dss
$S$\dss as above,\oss
then\qss
$\norm{S}\off =\off \norm{S'}$\nnsp,\oss
but\dss there is
no natural simplicial\dss map\qss
$S\ttoo S'$\nnsp.\oss
The subdivision of\dss chains\qss
$\alpha\qff \longmapsto\qff \fclass{\alpha}$\qss
serves as a substitute of\dss
such\dss map.\oss
There is
no natural\sss simplicial\dss map\qss $S'\ttoo S$\qss
either\halfff,\oss
but\dss there are\sss some distinguished\dss maps\qss $S'\ttoo S$\nnsp,\oss
namely,\oss
the\qss \emph{simplicial approximations}\qss of\dss
the identity\dss map\qss
$\id\dff \colon\dff
\norm{S'}\ttoo \norm{S}$\nnsp.\oss
In\dss the situation at\dss hand\dss the notion of\dss
simplicial approximation\dss reduces to\sss the following\halfff.\oss

Let\dss $S$\dss be a geometric simplicial complex
and\dss let\qss $x\qff \in\qff \norm{S}$\nnsp.\oss
Then\dss $x$\dss belongs\sss to at\dss least\sss one simplex of\dss $S$\nnsp,\oss
and since\sss the intersection of\dss simplices is also a simplex,\oss
there is unique minimal\qss (with respect\dss to\sss the inclusion)\qss
simplex of\dss $S$\dss containing\dss $x$\nnsp,\oss
called\dss the\qss \emph{carrier}\qss of\dss $x$\dss in\dss $S$\nnsp.\oss

Suppose now\dss that\dss $S'$\dss is a subdivision of\dss $S$\nnsp.\oss
Following\trs Alexander\qss \cite{a2},\oss
we will say\dss that\sss a simplicial\dss map\trs
$\varphi\dff \colon\dff
S'\ttoo S$\qss
is\qss \emph{pseudo-identical}\oss
if\qss for every vertex\sss $v$\sss of\qss $S'$\dss
its image\dss $\varphi\dff(\dff v\dff)$\dss
is a vertex of\dss the carrier of\dss $v$\dss in\dss $S$\nnsp.\oss
Equivalently,\oss
$\varphi$\dss is\qss \emph{pseudo-identical}\dff\oss
if\qss $\varphi\dff(\dff v\dff)$\qss is a vertex of\dss $\sigma$\dss 
for every vertex\sss $v$\sss of\qss $S'$\dss
and\sss every simplex\dss $\sigma$\dss of\dss $S$\dss containing\dss $v$\nnsp.\oss

Pseudo-identical\dss maps always exist\halfff.\oss 
Indeed,\oss let\qss
$\varphi\dff \colon\dff
v\dff(\dff S'\dff)\qff \ttoo\qff v\dff(\dff S\dff)$\qss
be\dss an\sss arbitrary\dss map\dss such\dss that\dss $\varphi\dff(\dff v\dff)$\dss
is a vertex of\dss the carrier of\dss $v$\qss for every vertex\sss $v$\sss of\qss $S'$\nnsp.\oss
Then\dss $\varphi$\dss is a pseudo-identical\sss simplicial\dss map.\oss
It\dss is sufficient\dss to check\sss that\dss 
$\varphi$\dss is a simplicial\dss map.\oss 
Let\sss $\tau$\sss be a simplex of\dss $S'$\nnsp.\oss
Since\dss $S'$\dss is a subdivision of\dss $S$\nnsp,\oss
the simplex\sss $\tau$\sss is contained in a simplex\dss $\sigma$\dss
of\dss $S$\nnsp.\oss
Every vertex of\dss $\tau$\dss belongs\sss to\sss $\sigma$\sss
and\dss hence its carrier is either\sss $\sigma$\sss or a proper 
face of\dss $\sigma$\nnsp.\oss
It\dss follows\sss that\dss $\varphi$\dss maps all vertices of\dss $\tau$\sss
to vertices of\dss $\sigma$\nnsp.\oss
It\dss follows\sss that\dss $\varphi$\dss
is a simplicial\dss map.\oss

\myuppar{Alexander\fff's\qss lemma.}
\emph{Let\qss $S'$\dss be\dss a\dss subdivision of\oss 
$S$\dss as above.\oss
If\oss
$\varphi\dff \colon\dff
S'\qff \ttoo\qff S$\qss
is\dss a pseudo-identical\dss simplicial\qss map,\oss
then\qss
$\varphi_*\dff(\dff \fclass{\alpha}\dff)
\off =\off
\alpha$\qss
for every chain\dss $\alpha$\dss of\oss $S$\nnsp.\oss}

\proof
The theorem is\dss trivially\dss true for $0$\dnsp-chains.\oss
Arguing\dss by\sss induction,\oss 
we may assume\dss that\dss it\dss is\dss true
for $m$\dnsp-chains with\qss $m\qff \leq\qff n\qff -\qff 1$\nnsp.\oss
Let\dss $\sigma$\dss be an $n$\dnsp-simplex of\dss $S$\dss
considered as an $n$\dnsp-chain.\oss
Then\dss $\partial\fff \sigma$\dss is an
$(\fff n\dff -\dff 1\fff)$\dnsp-chain
and\dss by\dss the inductive assumption\qss\vspace{0pt}
\[
\quad
\varphi_*\dff(\dff \fclass{\partial\fff \sigma}\dff)
\off =\off
\partial\fff \sigma
\qff.
\]

\vspace{-12pt}
On\dss the other hand,\oss
by applying\dss the above\dss lemma\dss  to\qss 
$\delta\off =\off \sigma$\nnsp,\oss 
we see\sss that\qss\vspace{0pt} 
\[
\quad
\partial\qff \fclass{\sigma}
\off =\off
\fclass{\partial\dff \sigma} 
\]

\vspace{-13.75pt}
and\dss hence\vspace{-1.75pt}
\[
\quad
\varphi_*\dff(\dff \partial\qff \fclass{\sigma} \dff)
\off =\off
\varphi_*\dff(\dff \fclass{\partial\dff \sigma}\dff)
\off =\off
\partial\dff \sigma
\qff.
\]

\vspace{-12pt}
We are almost\dss ready\dss to apply\dss Theorem\qss 2.\oss
Recall\dss that\dss
the collection\dss $S'\fff(\dff \sigma\dff)$\dss
of\dss all\sss simplices of\dss $S'$\dss contained\dss in\sss $\sigma$\sss
is a triangulation of\dss $\sigma$\nnsp.\oss
The subdivisions\dss
$\fclass{\partial\dff \sigma}$\dss of\qss
$\partial\dff \sigma$\dss with respect\dss to\dss $S'$\sss
and\sss with respect\dss
to\dss $S'\fff(\dff \sigma\dff)$\dss are\sss obviously\dss the same.\oss
Since\sss the simplicial\dss map\dss $\varphi$\dss is pseudo-identical,\oss
it\qss induces a simplicial\dss map\qss
$S'\fff(\dff \sigma\dff) \ttoo \sigma$\dnsp,\oss
which\dss is\dss also\dss pseudo-identical.\oss
The corresponding\dss induced\dss map\dss is simply\dss the
restriction of\dss $\varphi_*$\nnsp.\oss
By applying\dss Theorem\qss 2\qss to\sss this simplicial\dss map
and\dss the chain\dss $\fclass{\sigma}$\dss
in\dss the roles of\dss $\varphi$\dss and\dss $\alpha$\dss
respectively,\oss
we conclude\sss that\vspace{0pt}
\[
\quad
\varphi_*\dff(\dff \fclass{\sigma} \dff)
\off =\off
\sigma
\qff.
\]

\vspace{-12pt}
This proves\sss the\sss theorem\dss for $n$\dnsp-chains consisting of\dss one simplex.\oss
By\dss linearity\dss this implies\dss that\dss the\sss theorem\dss is\dss true for
all\sss $n$\dnsp-chains.\oss
An application of\dss induction completes\sss the proof\halfff.\oss  \eproof

\myuppar{Remarks.}
Alexander\fff's\qss lemma\qss is\dss the last\sss and\sss crucial\dss lemma in\dss the paper\qss \cite{a2}\qss
by\qss Alexander\halfff,\oss
devoted\dss to his second\qss
(and\dss the first\sss completely satisfactory)\qss 
proof\dss of\dss the\sss topological\dss
invariance of\trs Betti numbers\qss
(essentially,\qss of\dss the homology\dss groups)\qss 
of\dss polyhedra.\oss
Theorem\qss 2\qss together with\dss its proof\trs is\dss a\sss part\sss
of\qss Alexander's\dss proof\dss of\dss this lemma.\oss
Apparently,\oss it\dss was\qss Alexandroff\pss \cite{pa1}\qss
who elevated\dss this part\sss of\qss Alexander's\trs proof\trs to\sss a\dss theorem.\oss

Alexandroff\pss \cite{pa1}\qss
and\qss Alexandroff--Hopf\pss \cite{ah}\qss
viewed\qss Theorem\qss 1\qss
as\qss ``Theorem of\dss
conservation of\pss boundaries\sss by\sss simplicial\dss maps''\qss
(see\qss \cite{ah},\oss Chapter\qss IV,\oss Section\qss 3.7),\pss
Theorem\qss 2\qss as its counterpart\halfff,\pss
and\qss Alexander\fff's\qss lemma\qss as a natural extension of\qss Theorem\qss 2.\oss
Alexandroff\dss and\qss Hopf\qss gave all 
of\trs them\dss the name of\qss \emph{conservation\dss theorems}\pss
(\emph{``Erhaltungsatzes''}\pss in\dss German).

\mysection{Brouwer's\qss invariance\qss of\pss dimension\qss theorem}{dimension}

\myuppar{Systems of\dss sets and coverings.}
Let\dss $X\qff \subset\qff\rrr^{\dff d}$\nnsp.\oss
A\qss \emph{covering}\qss of\dss $X$\dss is\sss 
a finite system\qss\vspace{2pt}
\begin{equation}
\label{covering}
\quad
F_{\dff 0}\fff,\off F_{\dff 1}\fff,\off F_{\dff 2}\fff,\off \ldots\fff,\off F_{\fff s}\qss
\end{equation}

\vspace{-10pt}
of\dss subsets\dss $\rrr^{\dff d}$\dss such\dss that\dss $X$\dss is contained\dss
in\dss their union.\oss
The covering\qss (\ref{covering})\qss is said\dss to be\qss \emph{closed}\qss if\dss all\sss
$F_{\dff i}$\dss are closed.\oss
Let\dss $\varepsilon\qff >\qff 0$\nnsp.\oss
The covering\qss (\ref{covering})\qss is said\dss to be an $\varepsilon$\dnsp-covering\dss
if\trs the diameter of\dss every\dss $F_{\dff i}$\dss is\qss $<\qff \varepsilon$\nnsp.\oss
The\qss \emph{order}\pss of\dss a system of\dss sets\qss
(say,\pss a covering)\qss 
is\sss the maximal\dss number\dss $m$\dss such\dss that\sss there are\dss $m$\dss
different sets\sss in\dss the system\dss having\dss non-empty intersection.\oss

\myuppar{Lemma.}
\emph{Suppose\sss that\dss $m$\dss is\sss the order of\dss a closed\sss 
$\varepsilon$\dnsp-covering of\pss $\delta$\nnsp.\oss
If\pss $\varepsilon$\dss is sufficiently\dss small,\oss
then\dss there is a closed covering of\pss $\delta$\dss having\dss the order\qss
$\leq\qff m$\qss and\dss
consisting\dss of\pss $n\qff +\qff 1$\qss sets}\vspace{2pt}
\begin{equation}
\label{other-covering}
\quad
F_{\dff 0}\fff,\off F_{\dff 1}\fff,\off F_{\dff 2}\fff,\off \ldots\fff,\off F_{\fff n}
\qff,
\end{equation}

\vspace{-10pt}
\emph{such\dss that\qss $v_{\fff i}\qff \in\qff F_{\dff i}$\qss
and\pss $F_{\dff i}$\dss is disjoint\dss from\qss $\delta_{\dff i}$\qss
for every\qss $i\qff \in\qff I$\nnsp.\oss}

\proof
Let\qss $\varepsilon\qff >\qff 0$\qss be so small\dss that\dss no set\sss
of\dss diameter\qss $<\qff \varepsilon$\qss can\dss intersect\sss all\dss
$(\fff n\dff -\dff 1\fff)$\dnsp-faces\dss $\delta_{\dff i}$\dss of\dss $\delta$\nnsp.\oss
Then,\oss in\dss particular\halfff,\oss
no set\sss of\dss diameter\qss $<\qff \varepsilon$\qss can simultaneously
contain\sss a\sss vertex\dss $v_{\fff i}$\dss and a point\dss of\dss 
the\sss $(\fff n\dff -\dff 1\fff)$\dnsp-face\dss $\delta_{\fff i}$\dss
opposite\sss to it\halfff.\oss
Suppose\sss that\qss (\ref{covering})\qss is an
$\varepsilon$\dnsp-covering of\dss $\delta$\nnsp.\oss
By\dss the choice of\dss $\varepsilon$\dss
the\qss $n\qff +\qff 1$\qss vertices of\dss $\delta$\dss
belong\dss to\qss $n\qff +\qff 1$\qss different\sss sets\dss $F_{\dff i}$\nnsp.\oss
After\dss renumbering\dss the sets\dss $F_{\dff i}$\nnsp,\oss
if\dss necessary,\oss
we may assume\sss that\qss $v_{\fff i}\qff \in\qff F_{\dff i}$\qss
for all\qss $i\qff \in\qff I$\nnsp.\oss
Then\dss $F_{\dff i}$\dss is disjoint\dss from\dss $\delta_{\dff i}$\nnsp,\oss
again\dss by\dss the choice of\dss $\varepsilon$\nnsp.\oss

Suppose\sss that\dss the number\dss $s$\dss of\dss sets\qss
(\ref{covering})\qss is\qss $>\qff n\qff +\qff 1$\nnsp,\oss
and consider some set\dss $F_{\dff k}$\dss with\qss $k\qff >\qff n$\nnsp.\oss
By\dss the choice of\dss $\varepsilon$\dss the set\dss $F_{\dff k}$\dss
is disjoint\dss from some\sss $(\fff n\dff -\dff 1\fff)$\dnsp-face\dss $\delta_{\fff i}$\nnsp.\oss
Let\dss us consider\dss the union\qss $F_{\dff k}\qff \cup\qff F_{\fff i}$\nnsp.\oss
Clearly,\oss
$v_{\fff i}\qff \in\qff
F_{\dff k}\qff \cup\qff F_{\fff i}$\nnsp.\oss
On\dss the other\dss hand,\oss
both\dss $F_{\dff k}$\dss and\dss $F_{\dff i}$\dss
are disjoint\dss from\dss $\delta_{\dff i}$\dss
and\dss hence\sss the union\qss $F_{\dff k}\qff \cup\qff F_{\fff i}$\qss
is also disjoint\dss from\dss $\delta_{\dff i}$\nnsp.\oss
Let\dss us replace\sss the sets\qss $F_{\dff k}\fff,\off F_{\fff i}$\qss
by\dss their union\qss $F_{\dff k}\qff \cup\qff F_{\fff i}$\qss
and\dss rename\sss this union as\dss $F_{\fff i}$\nnsp.\oss 
This results in a new covering of\dss $\delta$\nnsp,\oss
which,\oss as we\dss just\sss saw,\oss satisfies\sss the last\sss
condition of\dss the lemma.\oss
Clearly,\oss this operation cannot\dss increase\sss the order of\dss the covering\halfff.\oss
By\dss repeating\dss this process we will eventually arrive at\sss
a covering consisting of\qss $n\qff +\qff 1$\qss sets
and satisfying\dss the\sss two other conditions of\dss the lemma also.\oss  \eproof\vspace{-0.125pt}

\myuppar{Lebesgue\dss lemma\sss for closed\sss sets.}
\emph{Suppose\sss that\oss \textup{(\ref{covering})}\pss is a system of\dss closed subsets of\trs
a compact\sss set\qss $X$\nnsp.\oss
Then\dss there is a number\pss $\varepsilon\qff >\qff 0$\qss with\dss
the following\dss property\fff:\oss
if\trs there is a point\dss of\qss $X$\dss whose distance from\dss several\dss sets
of\dss the system\oss \textup{(\ref{covering})}\pss is\pss $<\qff \varepsilon$\nnsp,\oss
then\dss these sets have non-empty\dss intersection.\oss
Every such\dss number\qss $\varepsilon\qff >\qff 0$\qss
is called a\oss \emph{Lebesgue\dss number}\pss of\dss the system\oss \textup{(\ref{covering})}.\oss}

\proof
Arguing\dss by contradiction,\oss
suppose\sss that\dss for every natural\dss number\dss $m$\dss
there is a point\qss $x_{\dff m}\qff \in\qff X$\qss
and a subsystem\dss $\mathcal{F}_{\fff m}$\sss of\dss the system of\dss sets\qss (\ref{covering})\qss
such\dss that\dss the distance of\dss $x_{\dff m}$\dss from every set\dss of\dss
the system\dss $\mathcal{F}_{\fff m}$\sss is\qss $<\qff 1/m$\nnsp,\oss
but\dss the intersection of\dss sets from\dss $\mathcal{F}_{\fff m}$\sss is empty.\oss
Since\qss (\ref{covering})\qss has only a finite number of\dss subsystems,\oss
we may assume,\oss after passing\dss to a subsequence if\trs necessary,\oss
that\sss all\sss subsystems\dss $\mathcal{F}_{\fff m}$\sss are\sss the same and\dss hence are
equal\dss to\dss $\mathcal{F}_{\dff 1}$\dnsp.\oss
Since\dss $X$\dss is compact\halfff,\oss we can also assume\sss that\dss the points\dss
$x_{\dff m}$\dss converge\dss to a point\qss $x\qff \in\qff X$\nnsp.\oss
Then\dss the distance of\dss $x$\dss from each of\dss the sets of\dss the subsystem\dss
$\mathcal{F}_{\dff 1}$\dss is equal\dss to\dss $0$\nnsp.\oss
Since\sss these sets are closed,\oss $x$\dss belongs\sss to all of\dss them,\oss
and\dss hence\sss the intersection of\dss sets from\dss $\mathcal{F}_{\dff 1}$\sss is non-empty,\oss
contrary\dss to\sss the assumption.\oss
The\dss theorem\dss follows.\oss  \eproof

\myuppar{Subdivisions into small simplices.}
We will\dss need\dss the following elementary\dss result\halfff:\oss
for every\qss $\varepsilon\qff >\qff~0$\qss
there exist\dss triangulations of\dss $\delta$\dss
consisting of\dss simplices of\dss diameter\qss $<\qff \varepsilon$\nnsp.\oss
It\dss is\sss not\dss hard\dss to believe\sss that\sss such\dss
triangulations exist\halfff,\oss and\dss this fact\sss is often accepted
without\sss even explicitly\sss stating\sss it\halfff.\oss
In fact\halfff,\oss even\dss more is\dss true.\oss
For every geometric simplicial complex\dss $S$\dss
there is a subdivision\dss $S'$\dss of\dss $S$\dss
consisting of\dss simplices of\dss diameter\qss $<\qff \varepsilon$\nnsp.\oss
In order\sss not\dss to interrupt\dss the flow of\dss ideas,\oss
the\dss proof\trs is deferred\trs till\qss Appendix\qss \ref{barycent}.\oss

\myuppar{Lebesgue-Sperner\dss theorem.}\oss
\emph{Suppose\dss that\pss \textup{(\ref{other-covering})}\qss
is a closed covering\dss of\qss $\delta$\dss
such\dss that\qss
$v_{\fff i}\qff \in\qff F_{\dff i}$\qss
and\qss
$F_{\dff i}$\dss
is disjoint\dss from\dss $\delta_{\dff i}$\dss for all\qss
$i\qff \in\qff I$\nnsp.\oss
Then\dss the order of\qss the covering\oss \textup{(\ref{other-covering})}\qss
is\oss $\geq\qff n\qff +\qff 1$\nnsp.\oss}

\proof
Let\qss $\varepsilon\qff >\qff 0$\qss be a\dss Lebesgue\dss number of\dss
the covering\qss (\ref{other-covering}),\oss
and\dss let\dss $T$\dss be a triangulation of\dss $\delta$\dss
such\dss that\dss the maximal diameter of\dss a simplex of\qss $T$\dss
is\qss $<\qff \varepsilon$\nnsp.\oss
Suppose\sss that\dss some simplex\dss $\sigma$\dss of\trs $T$\dss
intersects all\sss sets\qss (\ref{other-covering}).\oss
Then every\sss point\sss of\dss $\sigma$\dss has\sss the distance\qss
$<\qff \varepsilon$\dss from each of\dss the sets\qss (\ref{other-covering})\qss
and\dss hence\dss Lebesgue\dss lemma\dss implies\sss that\dss these\qss
$n\qff +\qff 1$\qss sets have non-empty\sss intersection.\oss
It\dss remains\sss to prove\sss that\sss such a simplex\dss $\sigma$\dss exists.\oss

Recall\dss that\dss $\Delta$\dss is\sss the simplicial complex
consisting of\dss the simplex\dss $\delta$\dss and\sss all\dss its faces,\oss
and\dss that\trs $T$\dss is a subdivision of\dss $\Delta$\nnsp.\oss
Let\dss us choose
for every\sss vertex\dss $v$\dss of\trs $T$\dss
a\dss set\dss $F_{\dff i}$\dss containing\dss $v$\dss
and\sss set\qss $\varphi\dff(\dff v\dff)\off =\off v_{\fff i}$\nnsp.\oss
Then\dss $\varphi$\dss is a map from\dss the set\sss of\dss vertices
of\dss $T$\dss to\sss the set\sss of\dss vertices of\dss $\Delta$\nnsp.\oss
Since every set\sss of\dss vertices of\dss $\Delta$\dss is\sss
the set\sss of\dss vertices of\dss a simplex of\dss $\Delta$\nnsp,\oss
the map\dss $\varphi$\dss is a simplicial\dss map\qss
$T\ttoo \Delta$\nnsp.\oss
Moreover\halfff,\oss as we will\sss see in a moment\halfff,\oss $\varphi$\dss is a pseudo-identical\dss
map.\oss
Indeed,\oss suppose\sss that\dss $v$\dss be a vertex of\dss $T$\dss
and\dss $\sigma$\dss is a simplex of\dss $\Delta$\dss such\dss that\qss $v\qff \in\qff \sigma$\nnsp.\oss
It\dss is sufficient\dss to show\sss that\dss in\dss this case\dss $\varphi\dff(\dff v\dff)$\dss
is a vertex of\dss $\sigma$\nnsp.\oss
If\dss not\halfff,\oss
then\qss $v\qff \in\qff F_{\dff i}$\qss
for some\qss $i\qff \in\qff I$\qss
such\dss that\dss $v_{\fff i}$\dss is not\sss a vertex of\dss $\sigma$\nnsp.\oss
In\dss this case\qss $\sigma\qff \subset\qff \delta_{\dff i}$\qss
and\dss hence\qss $v\qff \in\qff \delta_{\dff i}$\nnsp,\oss
contrary\dss to\dss $F_{\dff i}$\dss being disjoint\dss from\dss
$\delta_{\dff i}$\nnsp.\oss
It\dss follows\dss that\dss $\varphi$\dss is\sss pseudo-identical.\oss

Now\qss Alexander\fff's\qss lemma\qss implies\sss that\qss 
$\varphi_*\dff(\dff \fclass{\delta}\dff)
\off =\off
\delta$\qss
and\dss hence\qss
$\varphi\dff(\dff \sigma\dff)\off =\off \delta$\dss
for some\dss $n$\dnsp-simplex\dss $\sigma$\dss of\trs $T$\nnsp.\oss
By\dss the construction of\dss $\varphi$\dss this means for every\qss
$i\qff \in\qff I$\qss
some vertex of\dss $\sigma$\dss belongs\sss to\dss $F_{\dff i}$\dss
and\dss hence\dss $\sigma$\dss intersects all\sss sets\qss (\ref{other-covering}).\oss
This completes\sss the proof\halfff.\oss  \eproof

\myuppar{Lebesgue\dss tiling\qss (covering\fff)\dss theorem.}
\emph{If\pss $\varepsilon$\dss is sufficiently\dss small,\qff\oss
then every\dss $\varepsilon$\dnsp-covering of\pss $\delta$\dss
has\dss order\qss $\geq\qff n\qff +\qff 1$\nnsp.\oss}

\proof
It\dss is sufficient\dss to combine\dss Lebesgue-Sperner\dss theorem
with\dss the first\dss lemma.\oss  \eproof

\myuppar{Remarks.}
After\sss presenting\dss this proof\dss of\dss the\dss Lebesgue\dss tiling\dss theorem,\oss
Alexandroff\qss credits it\dss to\dss Sperner\dss and\dss Hopf\pss
(see\qss \cite{pa1},\oss footnote\qss 40)\fff:\vspace{-9pt}

\begin{quoting}
The above proof\dss of\dss the\sss tiling\dss theorem\dss
is due in essence\sss to\dss Sperner\halfff;\oss
the arrangement\trs given\trs here was communicated\dss to\sss me\dss by\trs
Herr\dss Hopf\halfff.\oss
\end{quoting}

\vspace{-9pt}
In\dss fact\halfff,\pss
the outline and\dss most\sss of\dss the details\sss of\dss this\dss proof\trs
are\sss the same as in\dss Sperner's\dss paper\qss \cite{s}.\oss
Alexandroff--Hopf\qss proof\trs differs from\trs Sperner's\trs one\sss only\dss 
in\dss using\trs Alexander's\trs lemma\dss
instead of\qss Sperner's\trs combinatorial\sss arguments.\oss 
This wouldn't\dss be possible without\dss assigning\dss to a vertex\dss $v$\dss of\trs $T$\dss
a\qss \emph{vertex\dss $v_{\fff i}$\dss of}\dss $\Delta$\dss and\dss
treating\dss the resulting\dss map\dss $\varphi$\dss
as a\qss \emph{simplicial\dss map}\qss $T\ttoo \Delta$\nnsp.\oss 
In contrast\halfff,\oss Sperner\dss assigns\sss to a vertex\dss $v$\dss
of\trs $T$\dss a\qss \emph{number}\pss $i\qff \in\qff I$\qss
(in\dss both versions\sss the as\-sign\-ment\dss is subject\dss to\sss 
the same condition\qss $v\qff \in\qff F_{\dff i}$\nsp).\oss

Nowadays\sss it\dss is only\dss natural\dss to\sss turn\dss
the set\qss
$I
\off =\off
\{\qff 0\fff,\pff 1\fff,\pff 2\fff,\pff \ldots\fff,\pff n\qff\}$\dss
of\dss subscripts enumerating\dss the sets\dss $F_{\dff i}$\qss
into an abstract\sss simplex and\dss then\dss identify\dss it\dss
with\dss the set\sss of\dss the vertices of\dss $\Delta$\nnsp.\oss
This was hardly\dss the case around\qss 1930,\oss
when\qss Alexandroff\qss wrote\sss his book\qss \cite{pa1}.\oss
But\dss the discovery of\trs this\sss idea was certainly\dss 
facilitated\dss by\dss the\sss notion of\dss the\qss \emph{nerve}\qss
of\dss a system of\dss sets,\oss introduced\dss by\trs
Alexandroff\pss \cite{pa0}\qss only a little earlier\halfff.\oss
One may speculate\sss that\qss Alexandroff's\trs contribution\dss
to\sss the above proof\dss is more significant\dss than
writing\sss down\dss Hopf's\dss version\sss of\trs Sperner's\dss proof\halfff.

The beautiful\dss idea of\dss
turning various sets and\dss maps into simplicial complexes and maps 
is well\sss established\dss by\dss now,\qss at\dss least\dss in some quarters.\oss
In skilled\dss hands it\dss is very\dss powerful.\oss

\myuppar{Canonical coverings of\dss simplicial\sss complexes by\dss closed\dss barycentric stars.}
Now we need a converse of\qss Lebesgue\dss tiling\dss theorem\qss
(see\dss Theorem\qss 3\qss below\fff).\oss
To begin with,\oss
for every\dss geometric $m$\dnsp-simplex\dss $\sigma$\dss
we will construct\dss
a canonical closed
covering\sss of\dss $\sigma$\dss by\qss $m\qff +\qff 1$\qss sets.\oss
Let\qss
$x_{\fff 0}\fff,\pff x_{\fff 1}\fff,\pff \ldots\fff,\pff x_{\fff m}$\qss
be\dss the vertices of\dss $\sigma$\nnsp,\oss
and\dss let\qss
$a_{\fff 0}\fff,\pff a_{\fff 1}\fff,\pff \ldots\fff,\pff a_{\fff m}$\qss
be\dss the barycentric coordinates of\trs points in\dss $\sigma$\qss
(see\dss Section\qss \ref{simplicial-complexes}).\oss
Recall\dss that\dss $a_{\fff i}$\dss are non-negative real\dss numbers with\dss the sum\dss $1$\nnsp.\oss
For every\qss
$k\off =\off 0\fff,\pff 1\fff,\pff 2\fff,\pff \ldots\fff,\pff m$\pss
let\dss $B_{\fff k}$\dss be\dss the set\sss of\dss points\dss $x$\dss
\vspace{4.5pt}
\[
\quad
x
\off =\off\dff
\sum\nolimits_{i\qff =\qff 0}^{m}\qff a_{\fff i}\dff x_{\fff i}
\]

\vspace{-7.5pt} 
of\dss $\sigma$\dss such\dss that\dss $a_{\fff k}$\dss is maximal\sss among\dss 
the barycentric coordinates\dss $a_{\fff i}$\dss of\dss $x$\nnsp.\oss
Clearly,\oss the sets\dss $B_{\fff k}$\dss form a closed covering of\dss $\sigma$\nnsp.\oss
The intersection of\dss all\dss these sets consists of\dss one point\halfff,\oss
the\qss \emph{barycenter}\qss of\dss $\sigma$\dnsp,\oss
which is\dss the only\dss point\dss for\sss which all\dss barycentric coordinates are equal.\oss 
Obviously,\oss $x_{\dff k}\qff \in\qff B_{\fff k}$\qss
and\qss $x_{\fff i}\qff \not\in\qff B_{\fff k}$\qss if\qss $i\off \neq\off k$\nnsp.\oss
Moreover\halfff,\oss
$B_{\fff k}$\dss is disjoint\dss from\dss the $(\fff m\dff -\dff 1 \fff)$\dnsp-face\dss of\dss
$\sigma$\dss opposite of\dss $x_{\dff k}$\nnsp,\oss
i.e.\qss from\dss the face having as its vertices\sss the points\dss $x_{\fff i}$\dss
with\qss $i\off \neq\off k$\nnsp.\oss
The sets\dss $B_{\fff k}$\dss naturally correspond\dss to\sss the vertices of\dss $\sigma$\dnsp,\oss
and\dss it\dss is convenient\dss to reflect\dss this in\dss the notations.\oss
Given a vertex\dss $v$\dss of\dss $\sigma$\dnsp,\oss
let\pss $B_{\fff v}\trf(\dff \sigma\dff)\off =\off B_{\fff k}$\nsp,\oss
where\dss $k$\dss is\sss such\dss that\qss $v\off =\off x_{\dff k}$\nnsp.\oss

Now,\oss let\dss $S$\dss be a\sss geometric simplicial complex 
and\dss let\dss $V$\dss be\sss the set\sss of\dss its vertices.\oss
The\qss \emph{(closed)\qss barycentric star}\qss of\dss a vertex\qss
is\dss the union\dss $B_{\dff v}$\dss
of\dss the sets\dss $B_{\fff v}\trf(\dff \sigma\dff)$\dss
over all\sss simplices\dss $\sigma$\dss of\trs $S$\dss having\dss $v$\dss as a vertex.\oss
Clearly,\oss the sets\dss $B_{\dff v}$\dss form a closed covering of\dss $\norm{S}$\nnsp.\oss

We claim\dss that\dss if\dss some point\dss of\dss the set\dss
$B_{\dff v}$\dss
belongs\sss to a simplex\dss $\tau$\dss of\trs $S$\nnsp,\oss
then\dss $v$\dss is a vertex of\dss $\tau$\nnsp.\oss
Suppose\sss that\qss
$x\qff \in\qff B_{\dff v}$\nnsp.\oss 
Then\qss
$x\qff \in\qff B_{\fff v}\trf(\dff \sigma\dff)$\qss
for some simplex\dss $\sigma$\dss having\dss $v$\dss as a vertex.\oss
If\dss $x$\dss belongs\sss to\dss $\tau$\nnsp,\oss
then\qss
$x\qff \in\qff \tau\qff \cap\qff \sigma$\dnsp.\qff\oss
The intersection\qss $\tau\qff \cap\qff \sigma$\qss is a simplex and\dss hence is
a face of\dss $\sigma$\nnsp.\oss
Since\dss $B_{\fff v}\trf(\dff \sigma\dff)$\dss 
is disjoint\dss from\dss the face of\dss $\sigma$\dss opposite\sss to\dss $v$\nnsp,\oss
it\dss follows\sss that\dss $v$\dss is a vertex of\qss
$\tau\qff \cap\qff \sigma$\qss
and\dss hence is a vertex of\dss $\tau$\nnsp.\oss
This proves our claim.\oss

Suppose now\sss that\dss $X$\dss is a subset\sss of\dss $V$\dss such\dss that\dss 
the intersection\vspace{1.25pt}
\[
\quad
\bigcap\nolimits_{\qff v\qff \in\qff X}\qff B_{\dff v} 
\]

\vspace{-10.75pt}
is non-empty.\oss
Let\dss $x$\dss be a point\sss in\dss this intersection and\dss let\dss
$\tau$\dss be some simplex of\trs $S$\dss containing\dss $x$\nnsp.\oss
By\dss the above claim\sss every\qss
$v\qff \in\qff X$\qss is a vertex of\dss $\tau$\nnsp.\oss
It\dss follows\dss that\dss $X$\dss is\sss the set\sss of\dss vertices
of\dss some face of\dss $\tau$\dss and\dss hence of\dss a simplex of\trs $S$\nnsp.\oss
Conversely,\oss if\dss $X$\dss is\dss the set\sss of\dss vertices of\dss
a simplex\dss $\sigma$\dnsp,\oss
then,\pss as we saw,\oss this intersection contains\dss the barycenter\sss of\dss $\sigma$\dnsp.\oss

Therefore\sss the order\sss of\dss the covering of\dss $\norm{S}$\dss
by\dss the sets\dss $B_{\dff v}$\dss is equal\dss to\sss 
the maximal\dss
number of\dss vertices of\dss a simplex of\trs $S$\nnsp,\oss
i.e.\qss is equal\dss to\qss $n\qff +\qff 1$\nnsp,\oss
where\dss $n$\dss is\sss the dimension of\trs $S$\nnsp.\oss
Clear\-ly,\oss the diameter of\dss each\dss
$B_{\dff v}$\dss is less\sss than\dss twice\dss the maximal\sss diameter
of\dss a simplex of\trs $T$\nnsp.\oss

\myuppar{Theorem\qss 3.}
\emph{For every\pss $\varepsilon\qff >\qff 0$\qss
there exists a closed\dss $\varepsilon$\dnsp-covering of\pss $\delta$\dss
of\trs the order\qss $n\qff +\qff 1$\nnsp.\oss}

\proof
Let\dss us apply\dss the above construction\dss to a\sss
triangulation\dss $T$\dss of\dss $\delta$\dss in\dss the role of\trs $S$\nnsp.\oss
If\dss the diameter of\dss simplices of\trs $T$\dss is\qss $<\qff \varepsilon/2$\nnsp,\oss
then\dss the covering\dss by the sets\dss $B_{\dff v}$\dss 
is an\dss $\varepsilon$\dnsp-covering\halfff.\oss  \eproof\vspace{-2.5pt}

\myuppar{Brouwer's\qss invariance\sss of\qss dimension\dss theorem.}
\emph{If\qss $m\qff <\qff n$\nnsp,\oss then\dss
no\dss subset\sss of\oss $\rrr^{\dff n}$\dss
with\dss non-empty\dss interior\dss
is homeomorphic\sss to a subset\sss
of\oss $\rrr^{\dff m}$\dnsp.\oss 
In\dss particular\halfff,\oss
$\rrr^{\dff n}$\dss is not\dss homeomorphic\sss to\pss $\rrr^{\dff m}$\dss
if\pss $n\off \neq\off m$\nnsp.\oss}

\proof
Any\sss subset\sss of\qss $\rrr^{\dff n}$ with\sss non-empty\sss interior 
contains a geometric $n$\dnsp-simplex,\oss
and\dss hence\sss it\dss is sufficient\dss to prove\sss that\sss
the $n$\dnsp-simplex\dss $\delta$\dss is\dss not\dss homeomorphic\sss to a subset\sss
of\qss $\rrr^{\dff m}$\dss if\qss $m\qff <\qff n$\nnsp.\oss
If\dss $h$\sss is a homeomorphism between\dss $\delta$\dss and\qss 
$X\qff \subset\qss \rrr^{\dff m}$\nnsp,\oss
then\dss $X$\dss is compact\dss together\sss with\sss $\delta$\sss 
and\dss hence is contained\sss in
a geometric $m$\dnsp-simplex.\oss
By\dss Theorem\qss 3\qss the latter admits closed $\varepsilon$\dnsp-coverings of\dss
order\qss $m\qff +\qff 1$\qss for every\qss $\varepsilon\qff >\qff 0$\nnsp.\pss
Since\dss $\delta$\dss is compact\sss and\dss hence\dss $h$\dss is uniformly continuous,\pss
transplanting\dss these coverings by\dss $h$\dss to\dss $\delta$\dss 
leads\sss to closed $\varepsilon$\dnsp-coverings of\dss $\delta$\dss
of\dss order\dss $m\qff +\qff 1$\dss
for every\dss $\varepsilon\qff >\qff 0$\dnsp,\qss
contrary\dss to\dss Lebesgue-Sperner\dss theorem.\oss  \eproof

\myuppar{Invariance of\dss dimension\dss for\dss polyhedra.}
\emph{Let\qss $S\fff,\off Q$\qss
be geometric simplicial complexes of\dss dimensions\qss $n\fff,\off m$\qss
respectively.\oss
If\qss $n\off \neq\off m$\nnsp,\oss
then\dss $\norm{S}$\dss and\dss $\norm{Q}$\dss
are not\sss homeomorphic.\oss}

\proof
If\dss $\norm{S}$\dss is homeomorphic\sss to\dss $\norm{Q}$\nnsp,\oss
then an open subset\sss of\dss an $n$\dnsp-simplex of\dss $S$\dss
is homeomorphic\sss to a subset\sss of\dss a $k$\dnsp-simplex\dss of\dss $Q$\dss
with\qss $k\qff \leq\qff m$\nnsp.\oss
By\dss the previous corollary\dss this implies\sss that\qss\
$k\qff \geq\qff n$\qss
and\dss hence\qss
$m\qff \geq\qff n$\nnsp.\oss
Similarly,\oss $n\qff \geq\qff m$\nnsp.\oss  \eproof

\mysection{Brouwer's\qss invariance\qss of\qss domain\qss theorem}{domains}

\myuppar{Brouwer's\dss invariance of\dss domain\dss theorem.}
The goal of\dss this section is\sss to apply\dss the\sss tools developed\dss
in\dss Section\qss \ref{dimension}\qss to prove another\dss famous\sss theorem
of\trs Brouwer\halfff,\oss namely,\oss to prove\sss that\dss if\dss two subsets of\dss $\rrr^{\dff n}$\dss
are homeomorphic and one of\dss them\sss is open,\oss then\dss the other is also open.\oss
The proof\dss is based only\sss on\dss Lebesgue-Sperner\dss theorem,\oss which is a version of\qss
Lebesgue\dss tiling\dss theorem,\oss
and elementary\dss constructions of\dss closed coverings\qss
(including\dss the existence of\dss subdivisions into small simplices).\oss
It\dss is\dss largely\dss due\sss to\dss Lebesgue\pss \cite{l-one},\pss \cite{l-two}.\oss

Lebesgue\qss endeavor\sss to prove\dss the invariance of\dss domain\dss theorem using only\dss
properties of\dss coverings\sss by closed sets was quite audacious,\pss
and\dss it\dss is not\sss very\sss surprising\dss that\dss his arguments contained a gap.\oss
The gap was filled\dss by\dss Sperner\qss \cite{s},\oss
but\dss not\dss without\sss a price:\oss
the resulting\dss proof\qss is dominated\dss by\dss the\sss technical\dss details.\oss
In\dss our\dss presentation\dss the arguments\dss filling\dss that\dss gap are separated\dss
from\dss the rest\sss of\dss the proof\dss as\dss the\sss technical\dss lemma below.\oss
Its proof\dss is based\dss on\dss Sperner\qss \cite{s},\oss 
but\sss differs in details,\oss
in\dss particular\halfff,\oss
in\dss the way\dss the compactness is used.\oss
The rest\dss is\dss based on a modernized\dss
version of\qss Lebesgue\dss ideas\pss \cite{l-one},\pss \cite{l-two}.\oss

\myuppar{Systems of\dss sets\sss differing\dss in a\sss subset\halfff.}
Let\qss (\ref{covering})\qss be a system of\dss sets in\dss $\rrr^{\dff d}$\nsp\dnsp,\oss
and\dss let\dss $Z$\dss be a subset\sss of\trs  $\rrr^{\dff d}$\nsp\dnsp.\dff\oss
A\dss system of\dss sets\qss 
\emph{differs\dss from}\qss (\ref{covering})\qss
\emph{only\dss in}\qss $Z$\qss 
if\qss it\dss consists of\dss subsets of\trs $Z$\dss
and sets\vspace{3pt}
\[
\quad
E_{\dff 0}\fff,\qff\off E_{\dff 1}\fff,\qff\off E_{\dff 2}\fff,\qff\off \ldots\fff,\qff\off E_{\dff s}
\]

\vspace{-9pt}
such\dss that\qss
$E_{\dff i}\off =\off F_{\dff i}$\qss
if\qss $F_{\dff i}$\dss is\dss disjoint\dss from\dss $Z$\dss
and\qss
$E_{\dff i}\qff \smallsetminus\pff Z
\off =\off
F_{\dff i}\qff \smallsetminus\pff Z$\qss
otherwise.\oss

\myuppar{Technical\qss lemma.}
\emph{Suppose\sss that\oss \textup{(\ref{covering})}\qss is a covering of\dss
the order\qss $\leq\qff n$\qss of\trs
a compact\dss set\oss $X\qff \subset\qff \rrr^{\dff d}$\qss 
and\dss that\pss $S$\dss is a geometric simplicial complex of\dss dimension\qss
$\leq\qff n\qff -\qff 1$\nnsp.\oss
Then\dss there is a covering of\qss the union\qss
$X\qff \cup\qff \norm{S}$\qss
which\dss has order\qss $\leq\qff n$\qss and\dss 
differs\dss from\oss \textup{(\ref{covering})}\qss only\dss in\qss $\norm{S}$\nnsp.\oss}

\proof
For every\qss
$i\off =\off 1\fff,\off 2\fff,\off \ldots\fff,\off s$\pss let\pss
$D_{\fff i}
\off =\off
F_{\dff i}\off \cap\off \norm{S}
$\nnsp.\oss 
Some of\dss the sets\dss $D_{\fff i}$\dss may\dss be empty,\oss
but\dss they\dss 
form a closed covering\sss
of\qss $X\qff \cap\qff \norm{S}$\qss with\dss the order\qss $\leq\qff n$\nnsp.\oss
Let\qss $e\qff >\qff 0$\dss be 
smaller\dss than\sss a\dss Lebesgue\dss number of\dss this covering\halfff,\oss
and\dss let\dss $D_{\fff i}\dff(\dff e\dff)$\dss
be\sss the closed $e$\dnsp-neigh\-bor\-hood of\qss $D_{\fff i}$\dss in\dss 
$\norm{S}$\nnsp,\oss
i.e.\qss the set\sss of\dss all\dss points of\dss 
$\norm{S}$\dss
with\dss the distance\qss $\leq\qff e$\qss from\qss $D_{\fff i}$\nnsp.\oss
If\dss $D_{\fff i}$\dss is empty,\oss
then\dss $D_{\fff i}\dff(\dff e\dff)$\dss is also empty.\oss
The sets\dss $D_{\fff i}\dff(\dff e\dff)$\dss form a closed covering of\qss 
$X\qff \cap\qff \norm{S}$\dnsp.\oss
If\qss 
$x\qff \in\qff X$\qss 
and\dss $x$\dss
belongs\sss to\sss the intersection
of\dss several\sss sets\dss $D_{\fff i}\dff(\dff e\dff)$\nnsp,\oss
then\dss the distance of\dss $x$\dss from each of\dss the corresponding
sets\dss $D_{\fff i}$\dss is\qss $\leq\qff e$\nnsp.\oss
By\trs Lebesgue\dss lemma\dss these sets\dss $D_{\fff i}$\dss
have non-empty\dss intersection.\oss
It\dss follows\sss the order of\dss the covering\dss by\dss 
the sets\dss $D_{\fff i}\dff(\dff e\dff)$\dss is\qss $\leq\qff n$\nnsp.\oss

There exists a subdivision\dss $S'$\dss of\dss $S$\dss 
consisting of\dss simplices of\dss diameter\qss $<\qff e/2$\qss
(see\dss Section\qss \ref{dimension}).\oss
By applying\dss the construction of\dss coverings\dss from\dss
Section\qss \ref{dimension}\qss to\dss $S'$\dss in\dss the role of\dss $S$\nnsp,\oss
we get a closed $e$\dnsp-covering\dss of\qss 
$\norm{S}\off =\off \norm{S'}$\dnsp.\oss
Let\qss 
$A_{\qff 0}\dff,\qff\off 
A_{\dff 1}\dff,\qff\off 
\ldots\dff,\qff\off 
A_{\dff p}$\qss
be\sss the sets of\dss this covering\dss intersecting\qss $X\qff \cap\qff \norm{S}$\qss
and\dss let\qss 
$B_{\dff 0}\dff,\qff\off 
B_{\dff 1}\dff,\qff\off 
\ldots\dff,\qff\off 
B_{\dff q}$\qss
be\sss the sets 
disjoint\dss from\qss $X\qff \cap\qff \norm{S}$\nnsp.\oss

Every\dss $A_{\dff k}$\dss intersects\dss $X\qff \cap\qff \norm{S}$\nnsp,\oss
and\dss hence intersects at\dss least\sss one of\dss the sets\dss $D_{\fff i}$\nnsp.\oss
Let\qss $D_{\fff i\dff(\fff k\dff)}$\qss be one of\dss these sets intersecting\dss $A_{\dff k}$\nnsp,\oss 
and\dss let\vspace{3.5pt}
\[
\quad
\overline{D}_{\fff i}
\off\off =\off\off
D_{\fff i}\off \cup \bigcup_{\off
i\qff =\qff i\dff(\fff k\dff)}\off A_{\dff k}
\off.
\]

\vspace{-8.5pt}
For every\qss $k\qff \leq\qff p$\dss
the diameter of\dss $A_{\dff k}$\dss is\qss $\leq\qff e$\qss
and\dss hence\dss $A_{\dff k}$\dss is contained\dss in\qss
$D_{\fff i\dff(\fff k\dff)}\dff(\dff e\dff)$\nnsp.\oss
It\dss follows\dss that\dss
$\overline{D}_{\fff i}$\dss
is contained\dss in\dss 
$D_{\fff i}\dff(\dff e\dff)$\dss
for every\qss $i\qff \leq\qff s$\nnsp.\oss
The sets\vspace{3.5pt}
\begin{equation}
\label{extended-covering}
\quad
F_{\dff 0}\qff \cup\qff \overline{D}_{\fff 0}\dff,\qff\off 
F_{\dff 1}\qff \cup\qff \overline{D}_{\fff 1}\dff,\qff\off  
\ldots\fff,\qff\off 
F_{\dff s}\qff \cup\qff \overline{D}_{\fff s}\fff,\off\off\off
B_{\dff 0}\dff,\qff\off 
B_{\dff 1}\dff,\qff\off 
\ldots\dff,\qff\off 
B_{\dff q}
\end{equation}

\vspace{-8.5pt}
form a closed covering of\qss $X\qff \cup\qff \norm{S}$\qss
which differs from\qss (\ref{covering})\qss only\dss in\dss $\norm{S}$\nnsp.\oss
It\dss remains\sss to prove\sss
that\dss the order of\dss this covering\dss  
is\qss $\leq\qff n$\nnsp,\qff\oss
i.e.\qss that\sss every\dss point\qss
$x\qff \in\qff X\qff \cup\qff \norm{S}$\qss
belongs\sss to\qss $\leq\qff n$\qss sets\qss
(\ref{extended-covering}).\oss 
There are several cases\sss to consider\halfff.\oss

If\qss $x\qff \not\in\qff X$\nnsp,\oss then none of\dss the sets\dss
$D_{\dff i}$\dss contains\dss $x$\nnsp.\oss
At\dss the same\sss time\dss $x$\dss belongs\sss to no more\sss than\dss $n$\dss
sets\qss $A_{\dff k}\fff,\pff B_{\dff l}$\nnsp.\oss
In\qss (\ref{extended-covering})\qss some sets\dss $A_{\dff k}$\dss are merged\dss together into one set\dss
$\overline{D}_{\fff i}$\nnsp,\oss but\sss are never split\halfff.\oss
Therefore in\dss this case\dss $x$\dss also belongs\sss to no more\sss than\dss $n$\dss
sets\qss (\ref{extended-covering}).\oss

If\qss 
$x\qff \in\qff X\qff \smallsetminus \qff \norm{S}$\nnsp,\oss 
then none of\dss the sets\dss
$B_{\dff l}$\dss contains\dss $x$\nnsp,\oss
and\qss $x\qff \in\qff F_{\dff i}\qff \cup\qff \overline{D}_{\fff i}$\qss
if\trs and only\dss if\qss $x\qff \in\qff F_{\dff i}$\nnsp.\oss
Since\sss the order of\pss (\ref{covering})\qss is\qss
$\leq\qff n$\nnsp,\oss
in\dss this case\dss $x$\dss belongs\sss to\qss
$\leq\qff n$\dss sets\qss (\ref{extended-covering}).\oss

If\qss $x\qff \in\qff X\qff \cap\qff \norm{S}$\nnsp,\oss
then\dss still\dss
none of\dss the sets\dss
$B_{\dff l}$\dss contains\dss $x$\nnsp,\oss
and\qss $x\qff \in\qff F_{\dff i}\qff \cup\qff \overline{D}_{\fff i}$\qss
if\trs and only\dss if\qss 
$x\qff \in\qff \overline{D}_{\fff i}\qff \subset\qff D_{\fff i}\dff(\dff e\dff)$\nnsp.\oss
Since\sss the order of\dss the covering\dss by\dss the sets\dss $D_{\fff i}\dff(\dff e\dff)$\dss
is\qss $\leq\qff n$\nnsp,\oss
in\dss this case\dss $x$\dss belongs\sss to\qss $\leq\qff n$\qss
sets\dss $D_{\fff i}\dff(\dff e\dff)$\nnsp,\oss
and\dss hence\sss to\qss 
$\leq\qff n$\qss sets\qss (\ref{extended-covering}).\oss  \eproof

\myuppar{Theorem\qss 4.}
\emph{Let\qss $X$\dss be a compact\sss subset\sss of\pss $\rrr^{\dff n}$\nnsp.\oss 
Sup\-pose that}\dss\vspace{3pt} 
\begin{equation}
\label{e-covering}
\quad
F_{\dff 0}\fff,\qff\off F_{\dff 1}\fff,\qff\off F_{\dff 2}\fff,\qff\off \ldots\fff,\qff\off F_{\fff s}
\end{equation}

\vspace{-9pt}
\emph{is closed\sss covering of\pss $X$\qss
such\dss that\dss its order\dss is equal\dss to\qss
$n\qff +\qff 1$\qss
and\dss only\sss one point\dss $y$\dss belongs\dss
to\qss $n\qff +\qff 1$\qss sets\oss \textup{(\ref{e-covering})}.\oss
If\qss $y$\dss belongs\sss to\sss the boundary of\oss $X$\nnsp,\oss
then\dss for every open set\oss $U\qff \subset\qff \rrr^{\dff n}$\qss containing\dss $y$\dss
there\dss exists a covering\dss of\oss $X$\qss of\qss order\qss $\leq\qff n$\qss
which\dss differs\dss from\oss \textup{(\ref{e-covering})}\pss only\dss in\dss $U$\nnsp.\oss}

\proof
Let\dss $\sigma$\dss be a geometric\dss $n$\dnsp-simplex
contained\dss in\dss $U$\dss and\sss
containing\dss $y$\dss in\dss its\dss
interior\qss 
$\int\dff \sigma\off =\off \sigma\qff \smallsetminus\qff \bd \sigma$\nnsp.\oss
Since\sss $y$\sss is\dss a\sss boundary\dss point\sss of\trs $X$\nnsp,\oss
there\dss exist\sss a point\dss $y'$\sss contained\dss in\dss $\int\dff \sigma$\dss
but\dss not\sss contained\dss in\dss $X$\nnsp.\oss 
Let\qss
$X'\off =\off X\qff \smallsetminus\qff \int\dff \sigma$\qss 
and\dss let\vspace{3pt}
\[
\quad
r\dff \colon\dff
X\qff \ttoo\qff X'\qff \cup\qff \bd \sigma
\]

\vspace{-9pt}
be\sss the map
equal\dss to\dss the identity on\dss $X'$\dss
and\dss to\sss the radial\dss projection\dss
from\dss $y'$\dss to\dss $\bd \sigma$\dss on\qss $X\qff \cap\qff \sigma$\dnsp.\oss
It\dss is well defined\dss because\qss $y'\qff \not\in\qff X$\nnsp.\oss
In\dss more details,\oss
if\qss $x\qff \in\qff X\qff \cap\qff \sigma$\dnsp,\oss
then\dss $r\dff(\dff x \dff)$\dss is\sss the point\sss of\dss
the intersection with\dss $\bd \sigma$\sss of\trs the ray\dss
starting at\dss $y'$\dss and\dss passing\dss through\dss $x$\nnsp.\oss

For every\qss
$i\off =\off 1\fff,\off 2\fff,\off \ldots\fff,\off s$\pss
let\qss
$E_{\dff i}
\off =\off
F_{\dff i}\off \smallsetminus\off \int\dff \sigma$\dnsp.\oss
The sets\dss $E_{\dff i}$\dss form a closed\sss covering\sss of\dss 
$X'$\dss
with\dss the order\qss $\leq\qff n$\qss
(since\sss the only\dss point\dss belonging\dss to\qss $n\qff +\qff 1$\qss
sets is removed\halfff).\oss
By\dss the\sss technical\dss lemma,\oss
there exists a closed covering\dss $\mathcal{G}$\dss  
of\dss the union\qss
$X'\qff \cup\qff \bd \sigma$\qss
with\dss the order\qss $\leq\qff n$\qss
which\dss differs\dss from\dss the system of\trs sets\dss $E_{\dff i}$\dss
only\dss in\dss $\bd \sigma$\nnsp.\oss
The collection of\dss preimages\dss $r^{\fff -\dff 1}\fff (\dff G\dff)$\dss
of\dss the sets\qss $G\qff \in\qff \mathcal{G}$\qss is a closed
covering\sss of\dss $X$\nnsp.\oss
Clearly,\oss its order\dss is\qss $\leq\qff n$\nnsp.\oss

It\dss remains\sss to check\dss that\dss this collection of\dss preimages differs
from\qss (\ref{e-covering})\qss only\dss in\dss 
$U$\dnsp.\oss
Actually,\oss it\dss differs
from\qss (\ref{e-covering})\qss only\dss in\dss 
$\sigma$\dnsp.\oss
If\qss $F_{\dff i}$\dss is disjoint\dss from\dss $\sigma$\dnsp,\oss
then\qss
$F_{\dff i}\off =\off E_{\dff i}$\qss
belongs\dss to\dss $\mathcal{G}$\dss
and\qss
$r^{\fff -\dff 1}\fff (\dff F_{\dff i}\dff)
\off =\off
F_{\dff i}$\nnsp.\oss
If\qss $F_{\dff i}$\dss intersects\dss $\sigma$\dnsp,\oss
then\dss $\mathcal{G}$\dss includes a set\dss
$G_{\dff i}$\dss
such\dss that\qss\vspace{1.75pt}
\[
\quad
G_{\dff i}\qff \smallsetminus\qff \bd \sigma
\off\qff =\off\qff
E_{\dff i}\qff \smallsetminus\qff \bd \sigma
\off\qff =\off\qff
F_{\dff i}\qff \smallsetminus\qff \sigma
\]

\vspace{-10.25pt}
and\dss hence\qss
$r^{\fff -\dff 1}\fff (\dff G_{\dff i}\dff)\qff \smallsetminus\qff \sigma
\off =\off
F_{\dff i}\qff \smallsetminus\qff \sigma$\dnsp.\oss
All other sets\dss $G$\dss from\dss the covering\dss $\mathcal{G}$\dss
are contained\dss in\dss $\bd \sigma$\dnsp,\oss
and\dss hence\sss their\dss preimages\dss $r^{\fff -\dff 1}\fff (\dff G\dff)$\dss
are contained\dss in\dss $\sigma$\nnsp.\oss
It\dss follows\dss that\dss the collection of\dss preimages
indeed\dss differs from\qss (\ref{e-covering})\qss only\dss in\dss 
$\sigma$\dnsp.\oss  \eproof

\myuppar{Brouwer's\qss invariance\sss of\qss domains\dss theorem.}
\emph{Suppose\sss that\pss
$Y\dff,\off X\qff \subset\pff \rrr^{\dff n}$\pss
and\qss
$h\dff \colon\dff Y\qff \ttoo\qff X$\qss
is a homeomorphism.\oss
If\qss $z$\dss belongs\sss to\sss the interior of\qss $Y$\nnsp,\oss
then\dss $h\dff(\dff z\dff)$\dss belongs\sss to\sss the interior of\oss $X$\dnsp.\oss
In\dss particular\halfff,\oss if\oss $Y$\dss is an open\dss subset\sss of\pss $\rrr^{\dff n}$\nnsp,\oss
then\qss $X$\dss is also an open subset\halfff.\oss}

\proof
Suppose\sss that\qss 
$y\off =\off h\dff(\dff z\trf)$\qss 
belongs\sss to\sss the boundary of\qss $X$\dnsp.\oss
Let\dss us choose a geometric $n$\dnsp-simplex\dss $\tau$\dss in\dss $\rrr^{\dff n}$\dss
such\dss that\dss $\tau$\dss is contained\dss in\dss the interior of\dss $Y$\dss
and\dss the point\dss $z$\dss is\dss the barycenter\qss
(see\dss Section\qss \ref{dimension})\qss of\dss $\tau$\nnsp.\oss
In\dss Section\qss \ref{dimension}\qss
we constructed a canonical closed covering\dss 
$\mathcal{F}$\dss of\dss $\tau$\dss by\qss $n\qff +\qff 1$\qss sets.\oss
The covering\dss $\mathcal{F}$\dss has order\qss $n\qff +\qff 1$\qss and\dss there is only one point\halfff,\oss
namely,\oss the barycenter\dss $z$\dss of\dss $\tau$\dnsp,\oss
which\dss belongs\dss to\qss $n\qff +\qff 1$\qss sets of\dss $\mathcal{F}$\nnsp.\oss
In addition,\oss the covering\dss $\mathcal{F}$\dss obviously satisfies\sss the assumptions of\qss
Lebesgue-Sperner\dss theorem.\oss
The image of\dss $\mathcal{F}$\dss under\dss $h$\dss is a closed covering\dss $\mathcal{F}'$\dss of\dss
$h\dff(\dff \tau\dff)$\dss of\dss order\qss $n\qff +\qff 1$\qss
and\dss  
$y\off =\off h\dff(\dff z\dff)$\dss
is\dss the only\dss point\dss
belonging\dss to\qss $n\qff +\qff 1$\qss sets of\dss $\mathcal{F}'$\nnsp.\oss
Since\sss the point\dss $y$\dss belongs\dss to\sss the boundary of\qss $X$\dnsp,\oss
it\dss belongs\dss to\sss the boundary of\qss $h\dff(\dff \tau\dff)$\dss also.\oss
Let\dss $U$\dss be an open\dss neighborhood\dss of\dss $h\dff(\dff z\dff)$\dss
such\dss that\dss the preimage\qss 
$h^{\dff -\dff 1}\dff(\dff U\dff)$\qss
is contained\dss in\dss the
interior\qss $\int\dff \tau\off =\off \tau\qff \smallsetminus\qff \bd \tau$\nnsp.\oss
Since\qss $\bd \tau\off =\off \norm{S}$\dss for an $(\fff n\dff -\dff 1 \fff)$\dnsp-dimensional
complex\dss $S$\nnsp,\oss
Theorem\qss 4\qss implies\sss that\dss there is a closed covering\sss 
of\dss $h\dff(\dff \tau\dff)$\dss 
with\dss the order\qss $\leq\qff n$\qss
which differs from\dss $\mathcal{F}'$\dss only\dss in\dss $U$\nnsp.\oss
The preimage $\mathcal{G}$ of\dss this covering\dss under\dss $h$\dss is
a closed covering of\dss $\tau$\dnsp.\oss
Clearly,\oss its order\dss is\qss $\leq\qff n$\qss
and\dss it\dss differs from\dss 
$\mathcal{F}$\dss
only\dss in\qss $h^{\dff -\dff 1}\dff(\dff U\dff)$\nnsp.\oss

Let\trs $G$\sss be an element\sss of\dss $\mathcal{G}$\dss 
containing some vertex of\dss $\tau$\nnsp.\oss
Let\dss us replace\dss $G$\dss by\dss the union of\trs $G$\dss and all\sss elements
of\dss $\mathcal{G}$\dss contained\dss in\dss $h^{\dff -\dff 1}\dff(\dff U\dff)$\nnsp,\oss
remove\dss from\dss $\mathcal{G}$\dss the latter elements,\oss
and\dss denote by\trs $\mathcal{G}'$\trs the resulting covering of\dss $\tau$\nnsp.\oss
Since\qss
$h^{\dff -\dff 1}\dff(\dff U\dff)\qff \subset\qff \int\dff \tau$\qss
and\dss $\mathcal{F}$\dss satisfies\sss the assumptions of\qss
Lebesgue-Sperner\dss theorem,\oss 
$\mathcal{G}'$\qss
also satisfies\sss these assumptions and\dss hence\sss its
order\dss is\qss $\geq\off n\qff +\qff 1$\nnsp.\oss
But\dss the previous paragraph implies\sss that\dss the order of\dss $\mathcal{G}'$\dss
is\dss $\leq\qff n$\nnsp.\oss  \eproof

\myuppar{Remark.}
Lebesgue\dss overlook\dss the need\dss to modify
and expand\dss the covering of\qss $X\qff \cap\qff \norm{S}$\qss
by\dss the sets\dss $D_{\fff i}$\dss in\dss the situation of\trs Theorem\dss 4,\oss
i.e.\qss when\qss $S\off =\off \bd \tau$\qss and\qss (\ref{covering})\qss
is\sss the image of\dss $\mathcal{F}$\nnsp.\oss

\mysection{Brouwer's\qss fixed-point\qss theorem}{simplicial-app}

\myuppar{Open stars.}
Let\dss $\sigma$\dss be a geometric $m$\dnsp-simplex and\dss 
let\dss $v$\dss be a vertex of\dss $\sigma$\dnsp.\oss
The $(\fff m\dff -\dff 1 \fff)$\dnsp-face of\dss $\sigma$\dss
having as its vertices all\sss vertices of\dss $\sigma$\dss
except\dss $v$\dss is called\dss the face\qss 
\emph{opposite}\pss to\dss $v$\nnsp.\oss
We already\dss informally\dss used\dss this notion\dss
and called\dss $\delta_{\dff i}$\dss the face of\dss $\delta$\dss opposite\sss to\dss $v_{\fff i}$\nnsp.\oss
If\dss $\sigma$\dss has dimension $0$\nnsp,\oss
i.e.\qss consists of\dss the point\dss $v$\nnsp,\oss
then\dss $v$\dss has no opposite face\qss
(alternatively,\oss one can consider\dss the empty set\sss as\dss the face opposite\sss to\dss $v$\nsp).\oss

Suppose now\dss that\dss $S$\dss is a geometric simplicial complex
and\dss $v$\dss is a vertex of\dss $S$\nnsp.\oss
The\qss \emph{open\dss star}\qss
of\dss $v$\dss in\dss $S$\dss
is\dss the subset\dss $\st(\dff v\fff,\qff S\dff)$\dss of\dss $\norm{S}$\dss
obtained\dss by\dss taking\dss the union of\dss all simplices of\dss $S$\dss
having\dss $v$\dss as a vertex with\dss the\sss faces opposite\sss to\dss $v$\dss removed.\oss
An easy exercise shows\sss that\sss every\sss open star\dss $\st(\dff v\fff,\qff S\dff)$\dss 
is indeed an open
subset\sss of\dss $\norm{S}$\nnsp.\oss
The following\dss lemma immediately\dss implies\sss
that\sss open stars form an open covering of\dss $\norm{S}$\nnsp.\oss

\myuppar{Lemma.}
\emph{Let\qss $x\qff \in\qff \norm{S}$\qss and\dss let\dss $\sigma$\dss be\dss 
the carrier\dss of\qss $x$\nnsp,\oss
i.e.\qss the minimal\dss simplex\sss of\pss $S$\sss containing\dss $x$\nnsp.\oss
Then\qss 
$x\qff \in\qff \st(\dff v\fff,\qff S\dff)$\qss
if\qss and\dss only\trs if\qss
$v$\dss is\dss a\sss vertex\dss of\qss $\sigma$\dnsp.\oss}

\proof
Since $\sigma$ is\sss the carrier of\dss $x$\nnsp,\oss
no proper\dss face of\dss $\sigma$ contains\dss $x$\nnsp.\oss
It\dss follows\dss that\qss
$x$\dss belongs\sss to\sss $\st(\dff v\fff,\qff S\dff)$\dss 
for every vertex\dss $v$\dss of\dss $\sigma$\dnsp.\oss
Conversely,\oss if\qss
$x\qff \in\qff \st(\dff v\fff,\qff S\dff)$\nnsp,\oss
then\dss $x$\dss belongs\sss to
some simplex\dss $\tau$\dss having\dss $v$\dss as a vertex\halfff,\oss
but\dss not\dss to its face opposite\sss to\dss $v$\nnsp.\oss
Since\dss $\sigma$\dss is\sss the minimal\sss simplex containing\dss $x$\nnsp,\oss
it\dss is contained\dss in\dss $\tau$\nnsp,\oss
and since\qss $x\qff \in\qff \sigma$\nnsp,\oss
it\sss cannot\dss be contained\dss in\dss the face opposite\sss to\dss $v$\nnsp.\oss
It\dss follows\sss that\dss $v$\dss is a vertex of\dss $\sigma$\dnsp.\oss  \eproof

\myuppar{Corollary.}
\emph{Let\qss
$w_{\fff 0}\fff,\pff w_{\fff 1}\fff,\pff \ldots\fff,\pff w_{\fff m}$\qss
be several\sss vertices of\pss $S$\dnsp.\oss
The intersection}\vspace{0.5pt}
\[
\quad
\bigcap\nolimits_{\qff i\qff =\qff 0}^{\qff m}\pff
\st(\fff w_{\fff i}\fff,\pff S\dff)
\]

\vspace{-10.75pt}
\emph{is non-empty\dss if\qss and\dss only\dss if\qss
$w_{\fff 0}\fff,\pff w_{\fff 1}\fff,\pff \ldots\fff,\pff w_{\fff m}$\qss 
are vertices of\qss a simplex of\pss $S$\dnsp.\oss}

\proof
If\dss $x$\dss belongs\sss to\sss this intersection and\dss $\sigma$\dss
is\dss the carrier of\dss $x$\nnsp,\oss
then all\dss $w_{\fff i}$\dss are vertices of\dss $\sigma$\dnsp.\oss
Conversely,\oss
if\dss $\sigma$\dss is a simplex with\dss the set\sss of\dss vertices\qss
$\{\qff
w_{\fff 0}\fff,\pff w_{\fff 1}\fff,\pff \ldots\fff,\pff w_{\fff m}
\qff\}$\nnsp,\oss
then\dss $\sigma$\dss is\sss the carrier of\dss every\qss
$x\qff \in\qff \sigma\qff \smallsetminus\qff \bd\dff \sigma$\qss
and\dss hence\qss
$\sigma\qff \smallsetminus\qff \bd\dff \sigma$\qss
is contained\dss in\dss this intersection.\oss  \eproof

\myuppar{Simplicial approximations.}
Let\qss $S\fff,\qff Q$\qss be simplicial complexes and\dss $f\dff \colon\dff
\norm{S} \ttoo \norm{Q}$\qss
be a continuous map.\oss
A simplicial\dss map\qss 
$\varphi\dff \colon\dff
S \ttoo Q$\qss
is called 
a\qss \emph{simplicial\dss approximation}\qss of\dss $f$\pss
if\vspace{1.5pt}
\[
\quad
f\qff \bigl(\dff \st(\dff v\fff,\qff S\dff) \dff\bigr)
\off \subset\off\qff
\st(\dff \varphi\dff (\dff v\dff)\fff,\qff Q\dff)
\]

\vspace{-10.5pt}
for every vertex\dss $v$\dss of\dss $S$\dnsp.\oss
Usually\dss $f$\qss admits no simplicial approximations.\oss
But\dss if\dss we allow\dss to replace\dss $S$\dss by\dss its subdivisions,\oss
the simplicial approximations always exist\halfff.\oss
The proof\dss is based on\dss a\sss version of\qss
Lebesgue\dss lemma from\dss Section\qss \ref{dimension}.\oss

\myuppar{Lebesgue\dss lemma\sss for open\sss coverings.}
\emph{Suppose\sss that\oss 
$\dis
\mathcal{U}
\off =\off 
\{\off U_{\dff i}\off |\bigr.\off i\qff \in\qff \mathcal{I} \off\}$\oss
is an open covering of\qss
a compact\sss set\qss $X$\nnsp.\oss
Then\dss there is a number\pss $\varepsilon\qff >\qff 0$\qss with\dss
the following\dss property\fff:\oss
if\qss the diameter\dss of\trs a subset\qss 
$Y$\dss of\oss $X$\qss
is\pss $<\qff \varepsilon$\nnsp,\oss
then\dss $Y$\dss is contained\dss in\qss $U_{\dff i}$\qss for some\qss 
$i\qff \in\qff \mathcal{I}$\nnsp.\oss 
Every such\dss number\qss $\varepsilon\qff >\qff 0$\qss
is called a\oss \emph{Lebesgue\dss number}\pss of\dss 
the covering\qss $\mathcal{U}$\nnsp.\oss}

\proof
Let\qss $F_{\dff i}\off =\off X\qff \smallsetminus\qff U_{\dff i}$\qss
for every\qss $i\qff \in\qff \mathcal{I}$\nnsp.\oss
The sets\dss $F_{\dff i}$\dss are closed.\oss
Arguing\dss by contradiction,\oss
suppose\sss that\dss for every natural\dss number\dss $m$\dss
there is a set\oss 
$Y_{\dff m}\qff \subset\off X$\oss 
of\dss diameter\qss $\leq\qff 1/m$\qss
not\sss contained\dss in any\dss $U_{\dff i}$\nnsp.\oss
Then\dss $Y_{\dff m}$\dss intersects every\dss $F_{\dff i}$\nnsp.\oss
Let\dss us choose some points\qss $y_{\dff m}\qff \in\qff Y_{\dff m}$\nnsp.\oss
Then\dss the distance of\dss $y_{\dff m}$\dss from every\dss $F_{\dff i}$\dss
is\qss $\leq\qff 1/m$\nnsp.\oss
Since\dss $X$\dss is compact\halfff,\oss
we can assume\sss that\dss the points\dss $y_{\dff m}$\dss
converge\sss to a point\qss $y\qff \in\qff X$\nnsp.\oss
Then\dss the distance of\dss $y$\dss from each set\dss $F_{\dff i}$\dss
is equal\dss to\dss $0$\nnsp.\oss
Since\sss these sets are closed,\pss $y$\dss belongs to every\dss $F_{\dff i}$\nnsp,\oss
and\dss hence does not\dss belong\dss to any\dss $U_{\dff i}$\nnsp.\oss
But\dss this contradicts\sss to\sss the assumption\dss
that\dss $\mathcal{U}$\dss is a covering\halfff.\oss  \eproof

\myuppar{The\dss simplicial approximation\dss theorem.}
\emph{Suppose\dss that\qss $S\fff,\qff Q$\qss are\sss two simplicial complexes
and\qss $f\dff \colon\dff
\norm{S}\qff \ttoo\qff \norm{Q}$\dss
is a continuous map.\oss
Let\qss $S'$\dss be a subdivision of\qss $S$\dnsp.\oss
If\qss the diameter of\dss the simplices of\oss $S'$\dss 
is\sss sufficiently\dss small,\oss
then\qss $f$\dss 
admits a simplicial approximation\qss
$S'\qff \ttoo Q$\nnsp.\oss
In\dss particular\halfff,\pss there exists\qss
$S'$\dss such\dss that\qss $f$\dss admits a simplicial approximation\qss
$S'\qff \ttoo Q$\nnsp.\oss}

\proof
The family of\dss open stars\dss 
$\st(\fff w\fff,\qff Q\dff)$\nnsp,\oss
where\dss $w$\dss runs over\dss the vertices of\dss $Q$\nnsp,\oss 
is an open covering of\dss $\norm{Q}$\dss
and\dss hence\sss the family of\dss their\dss preimages\vspace{2.5pt}
\[
\quad
f^{\qff -\dff 1}\dff \bigl(\dff \st(\fff w\fff,\qff Q\dff) \dff\bigr)
\]

\vspace{-9.5pt}
is an open covering of\dss $\norm{S}$\nnsp.\oss
Let\qss $\varepsilon\qff >\qff 0$\qss be\sss the\dss Lebesgue\dss number
of\dss this covering\halfff.\oss
Suppose\sss that\dss
the maximal\sss diameter of\dss a simplex
of\qss $S'$\dss is\qss $<\qff \varepsilon/2$\nnsp.\oss
Then\dss the diameter of\dss each open star\dss is\qss
$<\qff \varepsilon$\qss
and\dss hence every open star\dss $\st(\dff v\fff,\qff S'\dff)$\dss
is contained\sss in one of\dss the preimages.\oss
Equivalently,\oss for every vertex\dss $v$\dss of\dss $S'$\dss
there is a vertex\dss $w$\dss of\dss $Q$\dss such\dss that\vspace{1.5pt}
\[
\quad
f\qff \bigl(\dff \st(\dff v\fff,\qff S'\dff) \dff\bigr)
\off \subset\off\qff
\st(\dff w\fff,\qff Q\dff)
\qff.
\]

\vspace{-10.5pt}
For each vertex $v$ of\qss $S'$\dss
let\dss us choose one of\dss such vertices $w$
and denote it\dss by\dss $\varphi\dff(\fff v\dff)$\nnsp.\oss 
It\dss is sufficient\dss to check\dss that\dss the resulting\dss map
$\varphi$ is a simplicial\dss map\qss
$S'\qff \ttoo\qss Q$\nnsp,\oss
i.e.\qss that\dss it\dss maps simplices of\dss $S'$\dss to simplices of\dss $Q$\nnsp.\oss
Suppose\sss that\qss
$v_{\fff 0}\fff,\pff v_{\fff 1}\fff,\pff \ldots\fff,\pff v_{\fff m}$\qss
is\dss the set\sss of\dss vertices of\dss a simplex of\dss $S'$\dnsp.\oss
By\dss the above corollary\dss the intersection of\dss the open stars\qss
$\st(\dff v_{\fff i}\fff,\qff S'\dff)$\qss in non-empty.\oss
The image of\dss this intersection under\dss the map\dss $f$\qss
is contained\dss in\dss the intersection of\dss the open stars\qss
$\st(\dff \varphi\dff(\dff v_{\fff i}\dff)\fff,\qff Q\dff)$\nnsp,\oss
which is\dss therefore non-empty.\oss
Now\dss the same corollary\dss implies\dss that\qss 
$\varphi\dff(\dff v_{\fff 0}\dff)\fff,\pff
\varphi\dff(\dff v_{\fff 1}\dff)\fff,\pff
\ldots\fff,\pff
\varphi\dff(\dff v_{\fff m}\dff)$\qss
are vertices of\dss some simplex of\dss $Q$\nnsp.\oss  \eproof

\myuppar{Simplicial approximations and compositions.}
Let\qss $S\fff,\qff P$\nnsp,\oss and\dss $Q$\dss be simplicial complexes.\oss
Let\qss $f\dff \colon\dff
\norm{S} \ttoo \norm{P}$\qss
and\qss
$g\dff \colon\dff
\norm{P} \ttoo \norm{Q}$\qss
be continuous maps.\oss
If\qss
$\varphi\dff \colon\dff
S \ttoo P$\qss
and\qss
$\psi\dff \colon\dff
P \ttoo Q$\qss
are simplicial approximations of\dss the maps\dss
$f$\dss and\dss $g$\dss respectively,\oss
then,\oss obviously,\oss
$\psi\fff \circ\dff \varphi\dff \colon\dff
S \ttoo Q$\qss
is a simplicial\sss approximation of\qss
$g\dff \circ\fff f\dff \colon\dff
\norm{S} \ttoo \norm{Q}$\nnsp.\oss

\myuppar{Lemma.}
\emph{Suppose\dss that\trs $S'$\dss is a subdivision of\qss $S$\dss
and\pss
$\norm{S'} \ttoo \norm{S}$\qss
is\sss the identity\dss map.\oss
Every\dss simplicial\dss approximation\qss 
$\varphi\dff \colon\dff
S'\qff \ttoo\qff S$\qss
of\trs this identity\dss map\dss
is a pseudo-identical\dss map.}

\proof
By\dss the definition,\oss $\varphi$\dss is a simplicial approximation of\dss the
identity\dss if\qss and\dss only\trs if\qss\vspace{2.5pt}
\begin{equation}
\label{app-id}
\quad
\st(\dff v\fff,\qff S'\dff)
\off \subset\off\qff
\st(\dff \varphi\dff (\dff v\dff)\fff,\qff S\dff)
\end{equation}

\vspace{-9.5pt}
for every\sss vertex $v$ of\qss $S'$\dnsp.\oss
Clearly,\oss (\ref{app-id})\qss implies\dss
that\qss
$v\qff \in\pff
\st(\dff \varphi\dff (\dff v\dff)\fff,\qff S\dff)$\nnsp.\oss
This means\sss that\dss there is a simplex $\sigma$ of\trs $S$\dss
having\dss $\varphi\dff (\dff v\dff)$\dss as a vertex and
such\dss that\sss $v$ belongs\sss to $\sigma$
but\dss not\dss to\sss the face of\sss $\sigma$\ opposite\sss to\dss
$\varphi\dff (\dff v\dff)$\nnsp.\oss
It\dss follows\dss that\dss the carrier of\dss $v$ in\dss $S$\dss
has\dss $\varphi\dff (\dff v\dff)$\dss as a vertex.\oss
Since $v$ is an arbitrary vertex of\qss $S'$\dnsp,\oss
this means\sss that\dss $\varphi$\dss is a pseudo-identical\dss map.\oss  \eproof

\myuppar{Remark.}
The converse is also\sss true,\oss but\dss is\dss less useful\dss and\dss we do not\dss need\dss it\halfff.\oss

\myuppar{The no-retraction\dss theorem.}
\emph{There exists no retraction\qss
$\delta\qff \ttoo\qff \bd\dff \delta$\nnsp,\oss
i.e.\qss no continuous map\qss
$r\dff \colon\dff
\delta\qff \ttoo\qff \bd\dff \delta$\qss
which is equal\dss to\sss the identity on\dss $\bd\dff \delta$\nnsp.\oss}

\proof
Suppose\sss that\qss
$r\dff \colon\dff
\delta\qff \ttoo\qff \bd\dff \delta$\qss
is a continuous map equal\dss to\sss the identity on\dss $\bd\dff \delta$\nnsp.\oss
The sets\dss $\delta$\dss and\dss $\bd\dff \delta$\dss
are equal\dss to polyhedra\dss $\norm{\Delta}$\dss
and\dss $\norm{\nsp\bd\fff \Delta}$\dss respectively.\oss
By\dss the simplicial\sss approximation\dss theorem\dss
there exists a\sss subdivision\dss $T$\dss of\dss the complex\dss $\Delta$\dss
such\dss that\dss the continuous map\dss $r$\dss admits simplicial\sss approximation\qss
$\varphi\dff \colon\dff
T\qff \ttoo\qff \bd\fff \Delta$\nnsp.\oss

The restriction\qss 
$\bd\fff \varphi\dff \colon\dff
\bd\fff T\qff \ttoo\qff \bd\fff \Delta$\qss of\dss $\varphi$\dss
is a simplicial approximation of\dss
the restriction of\dss $r$\dss
to\sss the boundary\qss
$\bd\dff \delta$\nnsp,\oss
i.e.\qss a simplicial\sss approximation of\dss the identity\dss map of\dss
$\bd\dff \delta$\nnsp.\oss
By\dss the last\dss lemma\dss
the restriction\dss $\bd\fff \varphi$\dss is a pseudo-identical\dss map.\oss

We claim\dss that\dss $\varphi$\dss considered as a simplicial\dss map\qss
$T\qff \ttoo\qff \Delta$\qss
is a pseudo-identical\dss map.\oss
Let\sss $v$ be a vertex of\trs $T$\nnsp.\oss
If\qss $v\qff \in\qff \bd\dff \delta$\nnsp,\oss
then\dss $\varphi\dff(\fff v\dff)$\dss is a vertex of\dss carrier of\dss $v$\dss
because\dss $\bd\fff \varphi$\dss is a pseudo-identical\dss map.\oss
If\dss $v\qff \in\qff  
\delta\qff \smallsetminus\qff \bd\dff \delta$\nnsp,\oss
then\dss the carrier of\dss $v$\dss is\dss $\delta$\dss and\dss 
$\varphi\dff(\fff v\dff)$\dss is\dss tautologically\sss a vertex of\dss $\delta$\nnsp.\oss
This proves our claim.\oss
 
Now\dss Alexander\fff's\qss lemma\qss implies\dss that\qss
$\varphi\dff(\dff \sigma\dff)\off =\off \delta$\qss
for some simplex\dss $\sigma$\dss of\trs $T$\nnsp.\oss
But\dss $\varphi\dff(\dff \sigma\dff)$\dss
is a simplex of\dss $\bd\fff \Delta$\nnsp,\oss
i.e.\qss proper\dss face of\dss $\delta$\nnsp.\oss
The contradiction completes\sss the proof\halfff.\oss  \eproof

\myuppar{Brouwer's\qss fixed-point\dss theorem.}
\emph{Every continuous map\qss 
$\delta\qff \ttoo\qff \delta$\qss has a fixed\dss point\halfff.\oss}

\proof
As is well\dss known,\oss the geometric $n$\dnsp-simplex\dss $\delta$\dss is
homeomorphic\dss to an $n$\dnsp-dimensional\dss ball\dss $B$\dss
by a homeomorphism\dss taking\dss the boundary\dss $\bd\dff \delta$\dss
to\dss the boundary\dss $\partial\dff B$\dss of\dss the ball.\oss
The no-retraction\dss theorem\dss implies\sss that\dss there exists no retraction\qss
$B\qff \ttoo\qff \partial\dff B$\nnsp.\oss
On\dss the other\dss hand,\oss a continuous map\qss
$\delta\qff \ttoo\qff \delta$\qss without\dss fixed\dss points leads\sss to
a continuous map\qss 
$f\dff \colon\dff
B\qff \ttoo\qff B$\qss 
without\dss fixed\dss points,\oss
and\dss $f$\dss defines a retraction\qss
$r\dff\colon\dff
B\qff \ttoo\qff \partial\dff B$\qss
by\dss the following\dss rule\fff:\oss
if\qss $x\qff \in\qff B$\nnsp,\oss
then\dss $r\dff(\dff x\dff)$\dss
is\dss the point\sss of\dss intersection with\dss $\partial\dff B$\dss
of\dss the ray\dss going\dss from\dss $f\dff(\dff x\dff)$\dss to\dss $x$\nnsp.\oss  \eproof

\myuppar{The\qss Knaster-Kuratowski-Mazurkiewich\qss
argument\halfff.}
Now\sss we\sss turn\sss to\sss the celebrated\qss
Knaster-Kuratowski-Mazurkiewich\qss 
proof\dss
of\dss Brouwer's\dss fixed-point\dss theorem.\oss
The heart\sss of\trs this proof\trs is\dss an\sss ingenious use of\trs
of\trs barycentric coordinates,\oss  the\qss
Knaster-Kuratowski-Mazurkiewich\qss
ar\-gu\-ment\halfff.\oss
The rest\sss of\trs this section\dss is\dss devoted\dss to\sss this proof\dss
and an explanation of\trs its relation\dss with\dss the above proof\trs
based on simplicial\sss approximations.

\myuppar{KKM\dss theorem.}
\emph{Suppose\dss that\pss \textup{(\ref{other-covering})}\qss
are\qss $n\qff +\qff 1$\qss closed subsets of\oss $\delta$\nnsp.\qff\oss
If\pss
for\dss every\qss $J\qff \subset\qff I$}\vspace{3pt}\vspace{-4.25pt}
\[
\quad
\bigcup\nolimits_{\qff i\dff \in\dff J}\qff F_{\dff i}
\]

\vspace{-9pt}
\emph{contains\sss the face of\qss $\delta$\dss having\qss 
$\{\qff v_{\fff i}\qff \mid\qff
i\qff \in\qff J \pff\}$\qss as its set\sss of\dss vertices,\oss
then\oss
$\dis
\bigcap\nolimits_{\qff i\dff \in\dff I}\qff F_{\dff i}
\off \neq\off
\varnothing$\nsp.}

\proof
It\dss is a minor\dss modification of\dss the proof\dss of\qss
Lebesgue-Sperner\dss theorem\qss (see\dss Section\qss \ref{dimension}).\oss
The only difference is\dss in\dss the construction of\dss
simplicial\dss maps\qss
$\varphi\dff \colon\dff
T\qff \ttoo\qff \Delta$\nnsp.\oss
Let\dss $v$ be a vertex of\trs $T$\dss and\dss let\dss $\sigma$\dss
be\sss the smallest\dss face of\dss $\delta$\dss containing $v$\nnsp,\oss
i.e\qss the carrier of\dss $v$ in\dss $\Delta$\nnsp.\oss
Then\qss $v\qff \in\qff \sigma$\qss 
and\dss hence\sss there exists\qss $i\qff \in\qff I$\qss
such\dss that\qss
$v\qff \in\qff F_{\dff i}$\qss
and\dss 
$v_{\fff i}$\dss is a vertex of\dss $\sigma$\dnsp.\oss
Let\dss us choose an arbitrary\sss such\dss $i$\dss and\sss set\qss
$\varphi\dff(\dff v\dff)
\off =\off
v_{\fff i}$\nnsp.\oss
By\dss the construction,\oss $\varphi$\dss is a pseudo-identical\sss
simplicial\dss map.\oss
The rest\sss of\dss the proof\dss is\dss the same as before.\oss  \eproof

\myuppar{Notations.}
Let\dss us fix a continuous map\qss
$f\dff \colon\dff
\delta\qff \ttoo\qff \delta$\nnsp.\qff\oss
For a point\qss $a\qff \in\qff \delta$\qss we will\sss denote by\qss\vspace{1.5pt}
\[
\quad
(\qff
a_{\trf 0}\fff,\pff a_{\dff 1}\fff,\pff \ldots\fff,\pff a_{\fff n}
\qff)
\hspace*{1.5em}\mbox{and}\hspace*{1.5em}
(\qff
b_{\dff 0}\fff,\pff b_{\dff 1}\fff,\pff \ldots\fff,\pff b_{\fff n}
\qff)
\]

\vspace{-10.5pt}
be\sss the\sss barycentric coordinates of\dss $a$\dss and\qss
$b\off =\off f\dff(\dff a \dff)$\qss 
respectively.\oss

\prooftitle{A\dss proof\qss of\pss Brouwer's\qss fixed-point\trs theorem\dss
by\qss KKM\qss argument\halfff}
For\qss $i\qff \in\qff I$\qss let\dss $F_{\dff i}$\dss
be\sss the set\sss of\dss points\qss $a\qff \in\qff \delta$\qss
such\dss that\qss
$b_{\fff i}\qff \leq\qff a_{\dff i}$\nnsp.\oss
Clearly,\oss
the sets\dss $F_{\dff i}$\dss are closed.\oss

Suppose\sss that\qss $J\qff \subset\qff I$\qss and 
$a$ belongs\sss to\sss the face of\dss $\delta$\dss
having\qss 
$\{\qff v_{\fff i}\qff \mid\qff
i\qff \in\qff J \pff\}$\qss as its set\sss of\dss vertices.\oss
Then\dss the barycentric coordinate $a_{\dff i}$ can\dss be non-zero only\dss
for\qss $i\qff \in\qff J$\qss and\dss hence\vspace{4pt}
\[
\quad
\sum\nolimits_{\qff i\qff \in\qff J}\qff a_{\dff i}
\off\off =\off\off
1
\qff.
\]

\vspace{-8pt}
If\dss $a$ does not\dss belong\dss to any\dss $F_{\dff i}$\dss with\qss $i\qff \in\qff J$\nnsp,\oss
then\qss $b_{\fff i}\qff >\qff a_{\dff i}$\qss
for every\qss $i\qff \in\qff J$\qss and\dss hence\vspace{4.25pt}
\[
\quad
1
\off\off =\off\off
\sum\nolimits_{\qff i\qff \in\qff I}\qff b_{\fff i}
\off\off \geq\off\off
\sum\nolimits_{\qff i\qff \in\qff J}\qff b_{\fff i}
\off\off >\off\off
\sum\nolimits_{\qff i\qff \in\qff J}\qff a_{\dff i}
\off\off =\off\off
1
\qff.
\]

\vspace{-7.75pt}
The contradiction shows\sss that\qss
$a\qff \in\qff F_{\dff i}$\qss for some\qss $i\qff \in\qff J$\nnsp.\oss
It\dss follows\sss that\dss the sets\dss $F_{\dff i}$\dss satisfy\dss
the assumptions of\pss KKM\dss theorem
and\dss hence\sss the intersection of\dss all\dss these sets is non-empty.\oss
If\dss $a$\dss belongs\sss to\sss this intersection,\oss
then\qss
$b_{\fff i}\qff \leq\qff a_{\dff i}$\qss
for all\qss $i\qff \in\qff I$\nnsp.\oss
Since\sss the sum of\dss barycentric coordinates of\dss every\dss point\dss is\sss equal\dss to $1$\nnsp,\oss
this implies\sss that\qss
$a\off =\off b\off =\off f\dff(\dff a\dff)$\nnsp,\oss
i.e.\qss that\sss $a$ is\dss a\sss fixed\dss point\dss of\trs the map\dss $f$\dnsp.\oss  \eproof

\myuppar{Remarks.}
This\dss proof\qss is\dss  due\sss to\pss
Knaster-Kuratowski-Mazurkiewich\qss \cite{kkm}.\oss 
Thanks\dss to\sss this\sss proof\trs
a\sss minor\sss strengthening\dss of\pss Lebesgue-Sperner\qss theorem\dss
became known as\dss
Knaster-Kura\-towski-Mazurkiewich\dss or\dss KKM\dss theorem.\oss
There\dss is\dss no 
doubt\dss that\trs Sperner\dss would\dss easily\dss
prove\sss this\sss theorem,\pss
would\dss a need arise.\pss 
But\dss he was interested\sss in\sss the invariance of\dss dimension
and\dss domain\dss theorems,\qss
and\sss his version of\qss Lebesgue\dss ideas,\qss
which we called\qss Lebesgue-Sperner\qss theorem,\pss
is\dss precisely\sss tailored\dss for\sss his\sss goals. 

The\qss Knaster-Kuratowski-Mazurkiewich\qss argument\halfff,\pss
used\dss to deduce\dss Brouwer's\dss fixed-point\dss theorem\dss
from\dss KKM\dss theorem,\pss
is striking\halfff.\oss
The dependence on\dss KKM\dss theorem could\dss be easily
eliminated\qss (see another\dss proof\dss below),\oss
but\dss KKM\dss argument\sss
strongly contrasts with\dss conventional\dss proofs,\oss
invariably\dss reducing\dss the fixed-point\dss theorem\dss
to\dss the no-retraction\dss theorem.\oss

It\dss turns out\dss the standard construction of\dss a retraction\dss
from a map without\dss fixed\dss points\dss is\dss still\dss present\dss in\dss KKM\dss
argument\halfff,\oss but\dss in a veiled\dss form.\oss
Moreover\halfff,\oss KKM\dss argument\dss turns out\dss
to be similar\dss to\sss the proof\dss of\dss the\dss no-retraction\dss
theorem\dss based on\dss the simplicial\sss approximation\dss theorem.\oss
In order\dss to unveil\dss the retraction and explain\dss this similarity,\oss
we need\dss to rewrite\dss KKM\dss argument\sss
as a proof\trs by\sss contradiction starting\dss with\dss
assuming\dss that\dss $f$\dss has no fixed\dss points.\oss
As a\sss byproduct\halfff,\pss this\sss will\sss eliminate\sss the dependence on\dss
KKM\dss theorem.\oss\vspace{1.5pt}

\prooftitle{A\dss version of\pss the\qss KKM\qss proof\halfff}
Arguing\dss by contradiction,\oss suppose\sss that\dss $f$\dss has no fixed\dss points.\oss
Then\qss
$b\off \neq\off a$\qss for every\qss $a\qff \in\qff \delta$\nnsp.\oss
For\qss $i\qff \in\qff I$\qss let\dss $U_{\fff i}$\dss
be\sss the set\sss of\dss points\qss $a\qff \in\qff \delta$\qss
such\dss that\qss
$b_{\fff i}\qff <\qff a_{\dff i}$\nnsp.\oss
The sets\dss $U_{\fff i}$\dss are open,\oss but\dss this\sss is\sss irrelevant\dss for\dss the proof\halfff.\oss
Since\qss $b\off \neq\off a$\qss and\vspace{4.5pt}
\[
\quad
\sum\nolimits_{\qff i\qff \in\qff I}\qff b_{\fff i}
\off\off =\off\off
\sum\nolimits_{\qff i\qff \in\qff I}\qff a_{\dff i}
\off\off =\off\off
1
\qff,
\]

\vspace{-7.5pt}
there exists\qss $i\qff \in\qff I$\qss such\dss that\qss
$b_{\fff i}\qff <\qff a_{\dff i}$\nnsp.\oss
Therefore\sss the sets\dss $U_{\fff i}$\dss form a covering of\dss $\delta$\nnsp.\oss
For every\dss $i$\dss the $(\fff n\dff -\dff 1 \fff)$\dnsp-face\dss 
$\delta_{\dff i}$\dss of\dss $\delta$\dss is defined\dss by\dss the equation\qss
$a_{\dff i}\off =\off 0$\nnsp.\oss
Since\sss the barycentric coordinates are\qss $\geq\qff 0$\nnsp,\oss
if\qss $a\qff \in\qff \delta_{\dff i}$\nnsp,\oss
then\qss $a\qff \not\in\qff U_{\fff i}$\nnsp.\oss
In other words,\oss $U_{\fff i}$\dss is disjoint\dss from\dss $\delta_{\dff i}$\nnsp.\oss

Let\trs $T$\dss be a\sss triangulation of\dss $\delta$\nnsp.\oss
For every vertex $v$ of\trs $T$\dss there is at\dss least\sss one vertex\dss $v_{\fff i}$\dss
of\dss $\delta$\dss
such\dss that\qss $v\qff \in\qff U_{\fff i}$\qss
({\fff}because\sss the sets\dss $U_{\fff i}$\dss form a covering\fff).\oss
Let\dss
$\varphi\dff(\dff v\dff)$\dss be one of\dss them.\oss
If\dss $v$\dss belongs\sss to a face\dss $\sigma$\dss of\dss $\delta$\dss and\dss $v_{\fff k}$\dss
is not\sss a vertex of\dss $\sigma$\nnsp,\oss
then\qss $\sigma\qff \subset\qff \delta_{\dff k}$\qss
and\dss hence\dss $\sigma$\dss is disjoint\dss from\dss $U_{\fff k}$\nnsp.\oss
Therefore\qss $v\qff \not\in\qff U_{\fff k}$\qss and\dss hence\qss
$\varphi\dff(\dff v\dff)\off \neq\off v_{\fff k}$\nnsp.\oss
It\dss follows\dss that\dss $\varphi$\dss is a pseudo-identical\dss map.\oss
By\dss Alexander's\dss lemma\dss there exists a simplex\dss $\tau$\dss of\trs $T$\dss
such\dss that\qss
$\varphi\dff(\dff \tau\dff)\off =\off \delta$\nnsp.\oss
By\dss the construction of\dss $\varphi$\nnsp,\oss
the simplex\dss $\tau$\dss intersects every set\dss $U_{\fff i}$\nnsp.\oss

Let\qss $\varepsilon\qff >\qff 0$\nnsp.\oss
By\dss taking as\dss $T$\dss a\dss triangulation of\dss $\delta$\dss consisting
of\dss simplices of\dss diameter\qss $<\qff \varepsilon$\qss
we conclude\sss that\dss there exists a set\sss of\dss diameter\qss 
$<\qff \varepsilon$\qss
intersecting every set\dss $U_{\fff i}$\nnsp.\oss
Such a set\dss tautologically\dss intersects\sss also\sss the closures of\dss these sets.\oss
By applying\dss to\sss these closures\dss Lebesgue\dss lemma\dss for closed\sss sets
we see\dss that\dss these closures have a common\dss point\halfff.\oss
Clearly,\oss if\dss $a$\dss is such a common\dss point\halfff,\oss
then\qss
$b_{\fff i}\qff \leq\qff a_{\dff i}$\qss
for all\qss $i\qff \in\qff I$\nnsp.\oss
Since\sss the sum of\dss the barycentric coordinates of\dss a point\dss is always $1$\nnsp,\oss
it\dss follows\dss that\qss $a\off =\off b\off =\off f\dff(\dff a\dff)$\nnsp,\oss
i.e.\qss $a$\dss is a fixed\dss point\sss of\dss $f$\dnsp,\oss
contrary\dss to\sss the assumption\dss that\dss $f$\dss has no fixed\dss points.\oss
The contradiction completes\sss the proof\halfff.\oss  \eproof

\myuppar{Adjusting\dss the map\dss $f$\nsp\dnsp.}
Assuming\dss that\dss $f$\dss has no fixed\dss points,\oss
we would\dss like\sss  
to construct\sss a retraction\qss
$r\dff \colon\dff 
\delta\qff \ttoo\qff \bd\dff \delta$\nnsp.\oss
For every\qss  $x\qff \in\qff \delta$\qss
there is a unique ray starting at\sss $f\dff(\dff x \dff)$\sss and passing\dss through\dss $x$\nnsp,\oss
and\dss we would\dss like\sss to define\dss $r\dff(\dff x\dff)$\dss as\sss the point\sss of\dss intersection
of\dss this ray with\dss $\bd\dff \delta$\nnsp.\oss
But\dss the point\sss of\dss intersection\dss 
is\sss not\dss well\sss defined\sss 
when $x$ and\sss $f\dff(\dff x \dff)$\sss
belong\dss to\sss the same proper\dss face of\dss $\delta$\nnsp.\oss
This difficulty\sss is easy\dss to deal\sss with.\oss
If\dss  $f\dff(\dff \delta\dff)$\dss is contained\dss in\dss the interiour\qss
$\int\dff \delta\off =\off \delta\qff \smallsetminus\qff \bd\dff \delta$\nnsp,\oss
then\dss $r$\dss is well\sss defined.\oss
Since\sss there are maps arbitrarily close\sss to\dss $f$\dss with\dss the image contained\dss
in $\int\dff \delta$\nnsp,\oss
and every map sufficiently close\sss to\dss $f$\dss has no fixed\dss points,\oss
one can assume\sss that\qss
$f\dff(\dff \delta\dff)\qff \subset\qff \int\pff \delta$\qss
and\dss hence\dss $r$\dss is well\dss defined.\oss\vspace{1.25pt}

\myuppar{Bringing\dss forward\dss the retraction\dss $r$\dnsp.}
The central\sss element\sss of\dss the proof\dss is\dss the choice of\dss
a simplicial\dss map\qss
$\varphi\dff \colon\dff 
T\qff \ttoo\qff \Delta$\nnsp.\oss
The key\sss and\dss the most\sss original\dss feature of\dss the\dss KKM\dss argument\dss
is\dss the condition\dss imposed on\dss this choice\fff:\oss
given a vertex $v$ of\trs $T$\nnsp,\oss
one can\dss take as\dss $\varphi\dff(\dff v\dff)$\dss
any\dss vertex\dss $v_{\fff i}$\dss  
such\dss that\qss
$v\qff \in\qff F_{\dff i}$\nnsp,\oss
or\qss
$v\qff \in\qff U_{\dff i}$\dss in\dss the second version of\dss the proof\halfff.\oss

These conditions can\dss be easily\dss restated\dss in\dss terms of\dss $r$\nnsp.\oss
We\dss limit\sss ourselves\sss by\qss 
$v\qff \in\qff U_{\fff i}$\nnsp.\oss
Let\dss $\Lambda_{\dff i}\dff(\dff a\trf)$\dss
be\sss the set\sss of\dss points\qss $z\qff \in\qff \bd\dff \delta$\qss 
such\dss that\qss
$z_{\dff i}\qff >\qff a_{\dff i}$\nsp,\qss
where\dss
$(\trf
z_{\dff 0}\fff,\pff z_{\dff 1}\fff,\pff \ldots\fff,\pff z_{\fff n}
\qff)$\dss
are\sss the barycentric coordinates of\dss $z$\nnsp.\oss
Then\qss
$a\qff \in\qff U_{\fff i}$\pss
if\qss and\dss only\qss if\qss
$r\dff(\dff a \trf)
\qff \in\qff
\Lambda_{\dff i}\dff(\dff a\trf)$\nnsp.\oss
Therefore,\oss
in\dss the\qss (second version of\dss the)\qss 
KKM\dss argument\sss 
one can\dss take as\dss $\varphi\dff(\dff v\dff)$\dss
any\dss vertex\dss $v_{\fff i}$\dss 
such\dss that\qss
\begin{equation}
\label{kkm-r}
\quad
r\dff(\dff v \trf)
\qff \in\qff
\Lambda_{\dff i}\dff(\dff v\trf)
\qff.
\end{equation}
Let\qss
$\Lambda_{\dff i}
\off =\off
\Lambda_{\dff i}\dff(\dff v_{\fff i}\trf)
\off =\off
\bd\dff \delta\qff \smallsetminus\qss \delta_{\dff i}$\nsp.\oss
Then\qss
$\Lambda_{\dff i}\dff(\dff a\trf)
\pff \subset\pff
\Lambda_{\dff i}$\nnsp.\oss
The\sss set\sss $\Lambda_{\dff i}$\sss has\sss the advantage of\dss being\sss
independent\sss of\dss $a$\nnsp,\oss
and\dss one may\dss try\dss to
allow as\dss $\varphi\dff(\dff v\dff)$\dss any\dss vertex $v_{\fff i}$ such\dss that
\begin{equation}
\label{kkm-naive}
\quad
r\dff(\dff v \trf)
\qff \in\qff
\Lambda_{\dff i}
\qff.
\end{equation}
This condition\dss turns out\dss to be\sss too weak\dss to complete\sss
the proof\dss along\dss the lines of\dss the\dss KKM\dss argument\halfff.\oss
Namely,\oss by\trs
Lebesgue\dss lemma for closed sets there is\dss
a point\qss $a\qff \in\qff \delta$\qss such\dss that\dss
$r\dff(\dff a \trf)$\dss belongs\sss to\sss the closure of\dss $\Lambda_{\dff i}$\dss
for every\qss $i\qff \in\qff I$\nnsp.\oss
But\dss the closures of\dss the sets\dss $\Lambda_{\dff i}$\dss
have non-empty\dss intersection\qss
(equal\dss to\dss the union of\trs all\sss
$\bd\dff \delta_{\dff i}$\nsp),\oss
and we fail\dss to reach a contradiction.\oss
One can deal\dss with\dss this problem\dss by\dss
requiring\dss that\dss not\sss only\qss
$r\dff(\dff v \trf)
\qff \in\qff
\Lambda_{\dff i}$\nsp,\qff\oss
but\sss also
\begin{equation}
\label{kkm-simplicial}
\quad
r\dff(\dff \st (\dff v\dff)\dff)
\qff \subset\qff
\Lambda_{\dff i}
\qff.
\end{equation}
If\trs the diameters of\trs the simplices of\trs $T$\dss 
are less\sss than\dss the\dss Lebesgue\dss number of\dss 
the open covering of\dss $\delta$\dss
by\dss the preimages\dss
$r^{\dff -\dff 1}\dff(\dff \Lambda_{\dff i} \dff)$\nnsp,\qff\oss
then\sss such\dss $i\qff \in\qff I$\dss always exists.\oss
In\dss this case\dss Alexander\fff's\qss lemma implies\sss that\qss 
$\varphi\dff(\dff \sigma\dff)\off =\off \delta$\qss
for some simplex\dss $\sigma$\dss of\trs $T$\nnsp.\oss
Since\dss 
$\int\dff \sigma$\dss
is contained\dss in\dss the open star of\dss every vertex of\dss $\sigma$\dnsp,\oss
the image\sss $r\dff(\trf \int\dff \sigma \trf)$\sss is contained\dss in\dss
$\Lambda_{\dff i}$\dss for every\qss $i\qff \in\qff I$\nnsp.\oss
But\dss the intersection of\dss the sets\dss $\Lambda_{\dff i}$\dss
is obviously empty\halfff.\oss
We reached\dss the desired contradiction.\oss

The sets\dss $\Lambda_{\dff i}$\dss are\sss 
the open stars\dss $\st(\fff v_{\fff i}\fff,\pff \bd\dff \Delta\dff)$\nnsp,\oss
and\dss hence\sss the condition\qss (\ref{kkm-simplicial})\qss 
is equivalent\dss to\dss $\varphi$\dss being\sss a simplicial approximation of\dss $r$\nnsp.\oss
We are
back\dss to\sss the proof\dss from\dss Section\qss \ref{simplicial-app}.\oss
So,\oss the only essential\sss difference
of\dss the\dss KKM\dss argument\dss from\dss that\dss proof\dss
is\dss the\dss choice of\dss the condition\qss (\ref{kkm-r})\qss
instead of\qss (\ref{kkm-simplicial})\qss
as\sss the  
strengthening\dss of\qss the naive condition\qss (\ref{kkm-naive}).\oss

\mysection{Homology\qss groups\qss and\qss their\qss topological\qss invariance}{homology}

\myuppar{\nnsp$\partial\partial$\dnsp-theorem.}\oss
$\partial\dff \circ\qff \partial
\off =\off
0$\dnsp.\off\oss
\emph{In\dss more details,\oss
the composition of\qss the boundary operators}\vspace{0pt}
\begin{equation*}
\quad
\begin{tikzcd}[column sep=large]\dis
C_{\fff m}\dff(\dff S \dff)\dff\off
\arrow[r, 
"\dis \partial"]
& 
\off C_{\fff m\dff -\dff 1}\dff(\dff S \dff)
\arrow[r,   
"\dis \partial"]
& 
\off C_{\fff m\dff -\dff 2}\dff(\dff S \dff)
\end{tikzcd}
\end{equation*}

\vspace{-9.5pt}
\emph{is equal\dss to\dss $0$\dss for every simplicial\sss complex\qss $S$\dss
and\dss every\dss $m$\nnsp.\oss}

\proof
Clearly,\oss the versions dealing\dss with\dss geometric and\dss with\dss
abstract\sss simplicial\sss complexes are equivalent\halfff,\oss
and we can assume\sss that $S$ is\dss an abstract\dss simplicial\sss complex.\pss
For an $n$\dnsp-simplex\dss $\sigma$\dss and\sss a vertex\qss
$v\qff \in\qff \sigma$\qss let\dss us\dss denote\dss by\qss
$\sigma\qff -\qff v$\qss
the\dss $(\fff n\dff -\dff 1\fff)$\dnsp-face\qss
$\sigma\qff \smallsetminus\qff \{\dff v \trf\}$\qss
of\trs $\sigma$\dnsp.\oss
It\dss is\dss sufficient\dss to prove\sss that\qss
$\partial\dff(\trf \partial\fff\sigma \trf)\off =\off 0$\qss
for every simplex\dss $\sigma$\dnsp.\oss
Clearly,\oss\vspace{3pt}
\[
\quad
\partial\fff\sigma
\off\off =\off\off
\sum_{v\qff \in\qff \sigma\vphantom{X}}\qff 
\sigma\qff -\qff v
\]

\vspace{-12.5pt}
and\dss hence\vspace{2.5pt}
\[
\quad
\partial\dff(\trf \partial\fff\sigma \trf)
\off\off =\off\off
\sum_{v\qff \in\qff \sigma\vphantom{X}}\off
\sum_{w\qff \in\qff \sigma\qff -\qff v\vphantom{X}}\qff
\sigma\qff -\qff v\qff -\qff w
\off\off =\off\off
\sum_{v\fff,\qff w\vphantom{X}}\qff
\sigma\qff -\qff v\qff -\qff w
\qff,
\]

\vspace{-9pt}
where the last sum is taken over all\sss ordered pairs\dss $(\dff v\fff,\qff w \dff)$\dss
of\dss distinct vertices of\dss $\sigma$\nnsp.\oss
Every face of\dss $\sigma$\sss of\trs the form\qss 
$\sigma\qff \smallsetminus\qff \{\dff v\fff,\qff w \dff\}$\nnsp,\oss
where\qss $v\off \neq\off w$\nnsp,\oss
enters into this sum twice.\oss
Since we are working over\dss $\ftwo$\nnsp,\oss
this\dss implies\dss that\qss
$\partial\dff(\trf \partial\fff\sigma \trf)
\off =\off
0$\nnsp,\oss
completing the proof\halfff.\oss  \eproof

\myuppar{Cycles,\oss boundaries,\oss and homology.}
The identity\qss $\partial\qff \circ\qff \partial
\off =\off
0$\qss 
is the heart of\dss algebraic to\-pol\-o\-gy\halfff.\oss
It\dss is an algebraic form of\dss a simple geometric idea\fff:\oss
the boundary of\dss a geometric figure,\pss
considered as a geometric figure in\dss its own\dss right\halfff,\oss
has no boundary.\oss
Informally,\oss ``cycles''\qss are geometric figures with\sss no boundary.\oss
The\sss boundaries are cycles\sss by\dss a\sss trivial\dss reason,\pss
and one may ask\dss if\trs there are other cycles.\oss
This naturally\dss leads\sss to\sss treating\dss two cycles as equivalent\dss
when\halfff,\pss taken\dss together\halfff,\pss
they\dss form a boundary.\oss
These vague ideas can\dss be made precise only\dss if\dss
we specify\dss the class of\pss ``geometric figures''\qss considered.\oss
In our context\dss an appropriate class of\pss
``geometric figures''\qss
is\dss formed\dss by\dss
chains of\dss a simplicial\sss complex\qss (compare\dss Section\qss \ref{simplicial-complexes}).\oss

Let\dss $S$\dss be a simplicial complex and $m$ be a non-negative integer\halfff.\oss
An $m$\dnsp-chain $\alpha$ of\dss $S$\dss is called a\qss \emph{cycle}\oss
if\pss
$\partial\dff(\dff \alpha\dff)
\off =\off
0$\pss
or\qss $m\off =\off 0$\nnsp,\oss
and\dss a\qss \emph{boundary}\oss
if\pss
$\alpha
\off =\off
\partial\dff(\trf \beta\dff)$\qss
for\dss some $(\fff m\dff +\dff 1 \fff)$\dnsp-chain\dss $\beta$\nnsp,\oss
or\qss $\alpha\off =\off 0$\qss
and\sss $m$\sss is\dss the dimension of\qss $S$\dnsp.\oss
Let\vspace{3pt}
\[
\quad
Z_{\fff m}\fff(\dff S \dff)
\hspace*{1.0em}\mbox{and}\hspace*{1.0em}
B_{\fff m}\fff(\dff S \dff)
\]

\vspace{-9pt}
be\sss the spaces of\sss $m$\dnsp-chains which are cycles and\dss boundaries respectively.\oss
The $\partial\partial$\dnsp-theorem\sss implies\sss that\sss 
every $m$\dnsp-boundary\sss is\sss an $m$\dnsp-cycle,\oss
i.e.\qss
$B_{\fff m}\fff(\dff S \dff)
\qff \subset\off
Z_{\fff m}\fff(\dff S \dff)$\nnsp.\oss

The quotient\sss space\vspace{2pt}
\[
\quad
H^{\fff m}\fff(\dff S \dff)
\off =\off\qff
Z^{\fff m}\fff(\dff S \dff)\dff
\bigl/\bigr.\qff
B^{\fff m}\fff(\dff S \dff)
\]

\vspace{-10.375pt}
it\dss called\dss the
\emph{$m$\dnsp-dimensional\dss homology\dss group of}\qss
$S$\dnsp.\oss
The\sss term\qss \emph{``homology\dss group''}\oss
is standard even when\dss it\dss is a vector space.\oss
The image of\dss a cycle\qss
$\alpha\qff \in\qff Z_{\fff m}\fff(\dff S \dff)$\qss
in\dss the homology\dss group\qss
$H_{\fff m}\fff(\dff S \dff)$\qss
is called\dss the\qss \emph{homology\dss class}\qss of\dss $\alpha$\dss
and\dss is denoted\dss by\dss $\hclass{\alpha}$\nnsp.\oss

Let\pss
$\varphi\dff \colon\dff S\ttoo Q$\pss 
be a simplicial\dss map.\oss
Theorem\qss 1\qss implies\sss that\dss $\varphi_*$\dss
maps\dss $Z^{\fff m}\fff(\dff S \dff)$\dss to\dss $Z^{\fff m}\fff(\dff Q \dff)$\dss
and\dss maps\dss $B_{\fff m}\fff(\dff S \dff)$\dss to\dss $B_{\fff m}\fff(\dff Q \dff)$\nnsp.\oss
Therefore\dss $\varphi_*$\dss leads\sss to maps\vspace{2pt}
\[
\quad
\varphi_{*\nsp *}
\dff \colon\dff
H_{\fff m}\fff(\dff S \dff)
\qff \ttoo\qff
H_{\fff m}\fff(\dff Q \dff)
\]

\vspace{-10.375pt}
of\dss homology groups,\oss called\dss the\qss \emph{induced\dss maps}\qss in\dss homology.\oss

\myuppar{The\dss invariance\sss of\dss homology\dss groups under homeomorphisms and subdivisions.}
The rest\sss of\dss this section is devoted\dss to a proof\dss
of\dss the\sss topological\dss invariance of\dss 
the homology\dss groups\dss
$H_{\fff m}\fff(\dff S \dff)$\nnsp,\oss 
i.e.\qss to a proof\dss that\dss homology\dss groups\dss 
$H_{\fff m}\fff(\dff S \dff)$\dss is isomorphic\sss to\dss $H_{\fff m}\fff(\dff Q \dff)$\dss
if\trs the polyhedra
$\norm{S}$ and
$\norm{Q}$ are homeomorphic.\oss
The proof\trs is\dss based\dss on\dss Alexander's\dss ideas\qss \cite{a2},\pss
especially\dss on\trs Alexander's\qss lemma,\pss and\dss to some extent\dss
follows\trs Alexandroff's\dss exposition\qss \cite{pa1}.\oss

An\dss important\sss special case of\dss the\sss topological\dss invariance of\dss
homology\dss groups is\dss the following\halfff.\oss
Let\dss $S'$\dss be a subdivision of\dss a simplicial\sss complex\dss $S$\nnsp.\oss
Then\dss
$\norm{S}\off =\off \norm{S'}$\dss
and\dss hence\sss the homology\dss groups $H_{\fff m}\fff(\dff S \dff)$ and
$H_{\fff m}\fff(\dff S' \dff)$
are isomorphic.\oss
In\dss fact\halfff,\pss
it\dss is\sss easy\dss to describe\sss a canonical\dss isomorphism\qss
$H_{\fff m}\fff(\dff S \dff)\qff \ttoo\qff H_{\fff m}\fff(\dff S' \dff)$.\oss
Recall\dss the subdivision of\dss chains\sss map\qss
$\alpha\qff \longmapsto\qff \fclass{\alpha}$\qss
from\dss Section\qss \ref{alexander-lemma}.\oss
It\dss maps\dss $C_{\fff m}\fff(\dff S \dff)$\dss to\dss $C_{\fff m}\fff(\dff S' \dff)$\nnsp.\oss
The corollary\dss from\dss Section\qss \ref{alexander-lemma}\qss 
({\fff}preceding\dss the discussion of\dss pseudo-identical\dss maps)\qss
immediately\dss implies\dss that\dss the map\qss
$\alpha\qff \longmapsto\qff \fclass{\alpha}$\qss
maps cycles\sss to cycles and\dss boundaries\sss to boundaries.\oss
Hence\sss the subdivision\dss of\dss chains induces a homomorphism\qss
$H_{\fff m}\fff(\dff S \dff)\qff \ttoo\qff H_{\fff m}\fff(\dff S' \dff)$\dnsp,\oss
which we will call\dss the\qss \emph{homology subdivision\dss map}\pss
and\dss denote\dss by\qss
$h\qff \longmapsto\qff \fclass{h\fff}$\nnsp.\oss
The first\sss step\dss toward\dss the proof\dss of\dss 
the\sss topological\dss invariance of\dss homology groups is\dss to
prove\sss that\dss this map is an\dss isomorphism\dss when\dss $S'$\dss is\dss the so-called\dss
barycentric\dss subdivision.\oss

\myuppar{Barycentric\dss subdivisions.}
For every\dss geometric simplicial\sss complex\dss $S$\dss
there is a canonical subdivision $\bs S$ of\dss $S$\dnsp,\oss
called\dss the\qss \emph{barycentric subdivision}.\oss
See\dss Appendix\qss \ref{barycent}.\oss
Here we only\dss review\sss the main\dss properties of\dss $\bs S$\dnsp.\oss
Passing\dss
from\dss $S$\dss to $\bs S$ decreases\sss
the maximal diameter of\dss simplices by\dss the factor\qss $n/(\dff n\qff +\qff 1 \dff)$\nnsp,\oss
where $n$ is\dss the dimension of\dss $S$\dnsp.\pss
Therefore,\oss 
by\dss iterating\dss the construction of\dss the barycentric subdivision
one can construct\dss for every\qss $\varepsilon\qff >\qff 0$\dss
a subdivision\dss $S'$\dss of\dss $S$\dss such\dss that\dss the diameter of\dss
every simplex of\dss $S'$\dss is\qss $<\qff \varepsilon$\nnsp.\oss
See\dss Appendix\qss \ref{barycent}.\oss
Alexander\qss  \cite{a2}\qss observed\dss that\dss $\bs S$\dss
itself\dss can\dss be constructed\dss by\dss iterating\sss
a simpler\sss operation of\trs taking a\qss \emph{stellar\sss subdivision}\qss
and\dss proved\dss that\dss for stellar\sss subdivisions\sss the homology\sss subdivision\dss map\dss
is\dss an\sss isomorphism.\oss
This implies\sss that\dss the homology\sss subdivision\dss maps are isomorphisms
for simple and\dss iterated\dss barycentric subdivisions.\oss
Alexander's\dss proof\trs is\sss purely\sss combinatorial,\oss
but\dss fairly\dss technical.\oss
It\dss is\dss presented\dss in\dss Appendix\qss \ref{stellar}\halfff.\oss

\myuppar{The\sss topological\dss invariance\sss theorem\dss for\dss homology\dss groups.\off}
\emph{Let\pss $S$\dss and\pss $Q$\qss be\dss geometric\dss simplicial complexes.\qff\oss
If\pss the polyhedra\qss $\norm{S}$\dss and\pss $\norm{Q}$\dss
are homeomorphic,\oss then\dss the homology\dss groups\trs
$H_{\fff m}\fff(\dff S \dff)$\sss and\pss $H_{\fff m}\fff(\dff Q \dff)$\dss
are\sss isomorphic\dss for\dss every\sss $m$\nnsp.\oss}

\emph{Suppose\sss that\pss 
$f\dff \colon\dff
\norm{S}\qff \ttoo \norm{Q} $\qss
is a homeomorphism.\oss 
Then\dss there is an\dss iterated\dss barycentric\dss subdivision\qss $S'$\dss of\oss $S$\dss
such\dss that\dss $f$\dss admits\sss a\sss simplicial\sss approximation\qss
$\varphi\dff \colon\dff
S'\qff \ttoo\qff Q$\dss and\dss the induced\dss map\sss in\dss homology\oss 
$\varphi_{*\nsp *}\dff \colon\dff 
H_{\fff m}\fff(\dff S' \dff)
\qff \ttoo\qff 
H_{\fff m}\fff(\dff Q \dff)$\pss
is\dss an\dss isomorphism.}

\proof
Let\dss
$g\dff \colon\dff
\norm{Q}\qff \ttoo \norm{S}$\dss
be\sss the inverse of\dss the homeomorphism\dss $f$\dnsp.\oss 
By\dss the simplicial\sss approximation\dss theorem\dss there are\sss
subdivisions\qss
$S'\fff,\pff Q'\fff,\pff S''$\qss of\dss complexes\qss
$S\fff,\pff Q\fff,\pff S'$\qss respectively\sss and\sss simplicial\dss maps\qss
$\varphi\fff,\pff \vartheta\fff,\pff \psi$\qss
as on\dss the diagram\vspace{3pt}
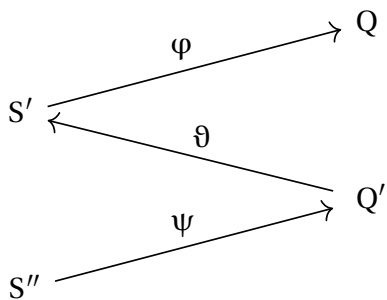
\[
\quad
\begin{tikzcd}[column sep=speccc, row sep=scriptsize]\dis 
&
Q
\\
S' 
\arrow[ur,   
"\dis \varphi\dnsp"]
&
\\
&
\pff Q'
\arrow[ul, 
"\dis \dnsp\vartheta"'] 
\\ 
\pff S''
\arrow[ur,   
"\dis \psi\dnsp\nsp"]
&
\end{tikzcd}
\]

\vspace{-9pt}
such\dss that\dss the maps\qss $\varphi\fff,\off \vartheta\fff,\off \psi$\qss 
are simplicial\sss approximations of\qss $f\dnsp,\off g\fff,\off f$\qss respectively.\oss

By\dss the discussion\dss preceding\dss the\sss theorem\dss
we can assume\sss that\dss
$S'\fff,\pff Q'\fff,\pff S''$\dss are\dss iterated\dss barycentric\sss subdivision of\qss 
$S\fff,\pff Q\fff,\pff S'$\qss respectively.\oss
Then\dss the homology subdivision\dss maps\qss\vspace{3pt}
\[
\quad
H_{\fff m}\fff(\dff S \dff)\qff \ttoo\qff H_{\fff m}\fff(\dff S' \dff)\dff,
\hspace{1.5em}
H_{\fff m}\fff(\dff S' \dff)\qff \ttoo\qff H_{\fff m}\fff(\dff S'' \dff)\dff,
\hspace{1.5em}\mbox{and}\hspace{1.5em}
H_{\fff m}\fff(\dff Q \dff)\qff \ttoo\qff H_{\fff m}\fff(\dff Q' \dff)
\]

\vspace{-9pt}
are isomorphisms.\oss
The above diagram 
leads\sss to\sss the diagram of\dss homology groups\vspace{3pt}
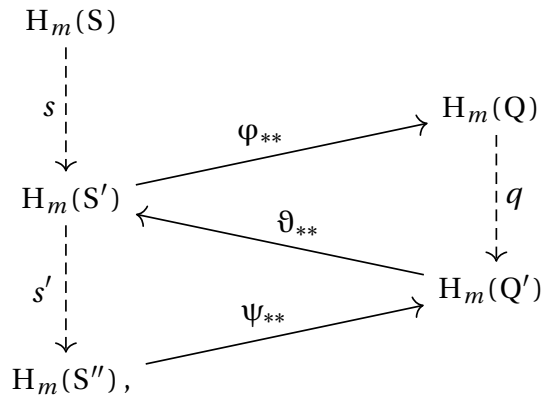
\[
\quad
\begin{tikzcd}[column sep=speccc, row sep=scriptsize]\dis
\off H_{\fff m}\fff(\dff S \dff)
\arrow[dashed]{dd}[']{\dis s\dff}
&
\\
&
H_{\fff m}\fff(\dff Q \dff)\off 
\arrow[dashed]{dd}{\dis q\qff}
\\
\off H_{\fff m}\fff(\dff S' \dff)
\arrow[ur,   
"\dis \varphi_{*\nsp *}\dnsp\dnsp"]
\arrow[dashed]{dd}[']{\dis s'\qff}
&
\\
&
H_{\fff m}\fff(\dff Q' \dff)\off 
\arrow[ul, 
"\dis \dnsp\dnsp\vartheta_{*\nsp *}"'] 
\\
\off H_{\fff m}\fff(\dff S'' \dff)\qff,
\arrow[ur,   
"\dis \psi_{\nsp*\nsp *}\dnsp\dnsp"]
&
\end{tikzcd}
\]

\vspace{-9pt}
where\dss the dashed\dss vertical arrows\dss
$s\fff,\pff s'\fff,\pff q$\dss are\sss the homology\dss subdivision\dss homomorphisms.\oss

By\dss Section\qss \ref{simplicial-app}\qss the composition\dss
$\varphi\dff \circ\dff \vartheta$\dss
is a simplicial\sss approximation of\dss the composition\dss $f\dff \circ\dff g$\dnsp,\oss
i.e.\qss of\dss the identity\dss map of\trs $\norm{Q}$\nnsp.\oss
Similarly,\oss $\vartheta\dff \circ\dff \psi$\dss
is a simplicial\sss approximation of\dss $g\dff \circ\dff f$\nnsp,\oss
i.e.\qss of\dss the identity\dss map of\trs $\norm{S}$\nnsp.\oss
Now\dss the second\dss lemma\dss from\dss Section\qss \ref{simplicial-app}\qss
implies\sss that\dss $\varphi\dff \circ\dff \vartheta$\dss
and\dss $\vartheta\dff \circ\dff \psi$\dss
are pseudo-identical\dss maps.\oss

Clearly,\oss
$\varphi_*\dff \circ\qff \vartheta_*
\off =\off
(\dff \varphi\dff \circ\dff \vartheta \dff)_*$\dss 
and\dss hence\qss Alexander's\dss lemma\dss implies\sss that\vspace{2.125pt}
\[
\quad
\varphi_*\dff \circ\qff \vartheta_*\dff \bigl(\dff \fclass{\alpha}\dff\bigr)
\off =\off
(\dff \varphi\dff \circ\dff \vartheta \dff)_*\dff \bigl(\dff \fclass{\alpha}\dff\bigr)
\off =\off\qff
\alpha
\]

\vspace{-9.875pt}
for every chain\dss $\alpha$\dss of\dss $Q$\dnsp.\oss
It\dss follows\sss that\dss if\dss $\alpha$\dss is a cycle and\qss 
$h$\dss is\sss its homology class,\oss then\vspace{2.125pt}
\[
\quad
\varphi_{*\nsp *}\dff \circ\qff \vartheta_{*\nsp *}\qff \bigl(\dff \fclass{h\fff}\dff\bigr)
\off =\off\qff
h
\qff.
\]

\vspace{-9.875pt}
Since\qss $q\dff(\dff h\dff)\off =\off \fclass{h\fff}$\nnsp,\oss
this means\sss that\qss
$\varphi_{*\nsp *}\dff \circ\qff \vartheta_{*\nsp *}\dff \circ\dff q$\qss
is\dss the identity\dss homomorphism of\qss
$H_{\fff m}\fff(\dff Q \dff)$\nnsp.\oss
Since $q$ is an\dss isomorphisms,\oss
this implies,\oss
in\dss particular\halfff,\oss
that\dss $\vartheta_{*\nsp *}$\dss is\dss injective.\oss

A completely similar argument\dss shows\dss that\qss
$\vartheta_{*\nsp *}\dff \circ\qff \psi_{\nsp*\nsp *}\dff \circ\dff s\fff'$\qss
is\dss the identity\dss homomorphism of\qss
$H_{\fff m}\fff(\dff S' \dff)$\nnsp.\oss
Since $s\fff'$ is an\dss isomorphisms,\oss
this implies,\oss
in\dss particular\halfff,\oss
that\dss $\vartheta_{*\nsp *}$\dss is\dss surjective.\oss

It\dss follows\dss that\dss $\vartheta_{*\nsp *}$\dss is\dss bijective
and\dss hence\sss is\sss an\dss isomorphism\qss
$H_{\fff m}\fff(\dff S' \dff)\qff \ttoo\qff H_{\fff m}\fff(\dff Q' \dff)$\nnsp.\oss
This already\dss implies\sss that\dss the homology\dss groups\dss
$H_{\fff m}\fff(\dff S \dff)$\dss and\dss $H_{\fff m}\fff(\dff Q \dff)$\dss
are isomorphic.\oss
Since\sss the composition\qss
$\varphi_{*\nsp *}\dff \circ\qff \vartheta_{*\nsp *}\dff \circ\dff q$\qss
is\dss the identity\dss homomorphism
and\dss $q$\dss is an\dss isomorphism,\oss
this implies\sss that\dss $\varphi_{*\nsp *}$\dss is also isomorphism.\oss  \eproof

\myuppar{The isomorphism of\qss homology\dss groups induced\dss
by\sss a homeomorphism.}
The proof\dss of\dss the\sss topological\dss invariance\sss theorem\dss
for\dss homology\dss groups proves\dss more\sss than\dss the isomorphism of\dss the\sss
groups\dss 
$H_{\fff m}\fff(\dff S \dff)$\dss and\dss $H_{\fff m}\fff(\dff Q \dff)$\nnsp.\oss
It\dss shows\dss that\sss a homeomorphism\dss
$f\dff \colon\dff
\norm{S}\qff \ttoo \norm{\hnsp Q}$\qss leads\sss to an\dss isomorphism\qss
$H_{\fff m}\fff(\dff S \dff)\qff \ttoo\qff H_{\fff m}\fff(\dff Q \dff)$\nnsp,\oss
namely,\pss to\sss the isomorphism\qss
$\varphi_{*\nsp *}\trf \circ\qff s$\nnsp.\oss
Moreover\halfff,\oss this isomorphisms depends only on\dss $f$\nsp\dnsp.\qff\oss
It\trs is\sss called\dss the\qss \emph{isomorphism\dss induced\dss by}\pss $f$\nnsp\dnsp.

The proof\dss is based on\dss the same ideas and\dss to a big extent\dss
is contained\dss in\dss the proof\dss of\dss the\sss topological\dss
invariance\sss theorem\dss itself\halfff.\oss
We need\dss the following\dss two lemmas.\oss
The first\sss one is almost\dss contained\dss in\dss the proof\dss of\dss
the\sss topological\dss invariance\sss theorem.\oss
The second one is a basic fact\dss about\dss the existence of\dss
simplicial\dss approximations.\oss

\myuppar{Lemma.}
\emph{Under\dss the assumptions of\qss the topological\dss invariance\sss theorem\dss
for\dss homology\dss groups,\oss
the\sss map\oss $\varphi_{*\nsp *}$\dss 
does not\dss depend\dss on\dss the choice\dss of\qss the\sss simplicial\dss
approximation\dss $\varphi$\nnsp.\oss}

\proof
It\dss is a continuation of\dss the proof\dss of\dss the\dss theorem.\oss
Since\dss
$\varphi_{*\nsp *}\dff \circ\qff \vartheta_{*\nsp *}\dff \circ\dff q$\dss
is\dss the identity\dss map of\qss
$H_{\fff m}\fff(\dff Q \dff)$\dss
and $q$\dss is\sss an\dss isomorphism,\pss
$q\trf \circ\trf \varphi_{*\nsp *}\dff \circ\qff \vartheta_{*\nsp *}$\dss
is\dss the identity\dss map of\qss
$H_{\fff m}\fff(\dff Q' \dff)$\nnsp.\oss
Therefore\dss $q\trf \circ\trf \varphi_{*\nsp *}$\dss
is\dss the inverse of\dss $\vartheta_{*\nsp *}$\nnsp.\oss
But\dss $\vartheta_{*\nsp *}$\dss is independent\sss on\dss the choice of\dss $\varphi$\dss
and\dss $\varphi_{*\nsp *}$\dss is independent\sss on\dss the choice of\dss $\vartheta$\nnsp.\oss
Therefore\sss
both\dss maps\dss $\vartheta_{*\nsp *}$\dss
and\dss $\varphi_{*\nsp *}$\dss are independent\sss on\dss the choice of\dss
simplicial\sss approximations\qss
$\vartheta\fff,\pff \varphi$\nnsp.\oss  \eproof

\myuppar{Lemma.}
\emph{If\oss $S$\dss is a geometric simplicial\sss complex\sss
and\pss $S'$\dss is a subdivision of\oss $S$\dnsp,\oss
then\dss the identity\dss map\qss
$\norm{S'}\qff \ttoo\qff \norm{S}$\qss
admits simplicial\sss approximation\qss
$S'\qff \ttoo\qff S$\dnsp.\oss}

\proof
Let\dss $w$\dss be a vertex of\qss $S'$\dss
and\dss let\dss $\sigma$\dss be\sss the carrier of\dss $w$\dss in\dss $S$\dnsp.\oss
We claim\dss that\qss\vspace{3pt}
\begin{equation}
\label{subdivision-stars}
\quad
\st(\dff w\fff,\qff S'\dff)\qff \subset\qff \st(\dff v\fff,\qff S\dff)
\qff.
\end{equation}

\vspace{-9pt}
for every vertex\dss $v$\dss of\dss $\sigma$\dnsp.\oss
Let\dss us consider some simplex\dss $\tau$\dss of\qss $S'$\dss
having\dss $w$\dss
as a vertex\halfff.\oss
Since\dss $S'$\dss is a subdivision of\trs $S$\dnsp,\oss
the simplex\dss $\tau$\dss is contained\dss in some simplex\dss $\rho$\dss of\qss $S$\dnsp.\oss
The simplex\dss $\rho$\dss contains $w$ and\dss hence has\sss the carrier\dss $\sigma$\dss
of\dss $w$\dss as a face.\oss
It\dss follows\sss that\dss $\sigma$\dss has $v$ as a vertex.\oss
The intersection of\dss $\tau$\dss with\dss the face of\dss $\rho$\dss
opposite\sss to\dss $v$\dss is a face of\dss $\tau$\qss
(perhaps,\oss empty)\qss 
not\sss containing\dss $w$\nnsp.\oss
Therefore\sss this intersection\dss is\sss contained\dss in\dss the face of\dss
$\tau$\dss opposite\sss to\dss $w$\nnsp.\oss
The inclusion\qss (\ref{subdivision-stars})\qss follows.\oss
One gets a  simplicial\sss approximation of\dss the identity
by choosing\dss for every\sss vertex\dss $w$\dss of\qss $S'$\dss
some vertex\dss $v$\dss of\qss $S$\dss as above.\oss  \eproof

\myuppar{Theorem.}
\emph{Under\dss the assumptions of\qss the topological\dss invariance\sss theorem\dss
for\dss homology\dss groups,\oss 
the composition\qss
$\varphi_{*\nsp *}\dff \circ\qff s\dff \colon\dff
H_{\fff m}\fff(\dff S \dff)\qff \ttoo\qff H_{\fff m}\fff(\dff Q \dff)$\nnsp,\oss
where\qss
$s\dff \colon\dff
H_{\fff m}\fff(\dff S \dff)\qff \ttoo\qff H_{\fff m}\fff(\dff S' \dff)$\qss
is\dss the homology\dss subdivision\dss map,\oss
depends only\sss on\qss $f$\dnsp.\oss}

\proof
Since\dss $S'$\dss is an\dss iterated\dss barycentric
subdivision of\dss $S$\dnsp,\oss
it\dss is sufficient\dss to prove\sss that\dss this composition does not\dss
change if\qss $S'$\dss is replaced\dss by an
arbitrary\dss iterated\dss barycentric subdivision\dss
$S''$\dss of\qss $S'$\dnsp.\oss
In\dss more details,\oss
let\dss $s\fff'$\dss be\sss the homology\sss subdivision\dss map\qss
$H_{\fff m}\fff(\dff S' \dff)\qff \ttoo\qff H_{\fff m}\fff(\dff S'' \dff)$\nnsp.\oss
Then\sss $s\fff'\dff \circ\dff s$\sss 
is\dss the homology subdivision\sss map\dss
$H_{\fff m}\fff(\dff S \dff)\qff \ttoo\qff H_{\fff m}\fff(\dff S'' \dff)$\nnsp,\pss
and we need\dss to show\sss that\vspace{3pt}
\[
\quad
\varphi_{*\nsp *}'\dff \circ\qff (\dff s\fff'\dff \circ\dff s \dff)
\off =\off
\varphi_{*\nsp *}\dff \circ\qff s
\]

\vspace{-9pt}
if\qss
$\varphi'\dff \colon\dff
S''\qff \ttoo\qff Q$\qss
is a simplicial\sss approximation\sss
of\pss $f$\nsp\dnsp.\oss
By\dss the first\dss lemma\dss $\varphi_{*\nsp *}'$ does not\sss depend on\dss
the choice of\dss the simplicial\sss approximation\dss $\varphi'$\dss
and\dss hence we are free\sss to choose\dss $\varphi'$\dnsp.\oss 

By\dss the second\dss lemma\sss 
there exists a simplicial\sss approximation\qss
$\lambda\dff \colon\dff
S''\qff \ttoo\qff S'$\qss
of\dss the identity map\sss
$\norm{S''}\qff \ttoo\qff \norm{S'}$\nnsp.\oss
By\dss Section\qss \ref{simplicial-app}\pss
the composition\qss 
$\varphi\dff \circ\dff \lambda
\dff \colon\dff
S'\qff \ttoo\qff Q$\qss
is\sss a\sss simplicial\dss approximation\sss of\qss $f\dff \circ\qff \id
\off =\off 
f$\nsp\dnsp.\oss
Hence we can\dss take\qss 
$\varphi'
\off =\off
\varphi\dff \circ\dff \lambda$\nnsp.\oss
By\dss Alexander's\dss lemma\qss\vspace{2.625pt}
\[
\quad
\lambda_{*\nsp *}\dff \circ\qff s\fff'
\off =\off
\id
\qff,
\]

\vspace{-9.375pt}
and\dss hence\oss
$\varphi_{*\nsp *}'\dff \circ\qff (\dff s\fff'\dff \circ\dff s \dff)
\off =\off
\varphi_{*\nsp *}\dff \circ\qff \lambda_{*\nsp *}\dff \circ\qff s\fff'\dff \circ\dff s
\off =\off
\varphi_{*\nsp *}\dff \circ\qff s
$\nnsp.\oss  \eproof

\myuppar{Remarks.}
The same construction applies if\dss $f$\dss is only a continuous map
and\dss leads\sss to a homomorphism\qss
$H_{\fff m}\fff(\dff S \dff)\qff \ttoo\qff H_{\fff m}\fff(\dff Q \dff)$\nnsp,\oss
called\dss the\qss \emph{homomorphism\dss induced}\pss by\dss $f$\dss
and\dss usually denoted\dss by\dss $f_*$\nnsp.\oss
It\dss is also independent\sss on\dss the choices involved\dss in\dss its construction,\oss
but\dss the proof\qss is\dss more\sss technical\dss than\dss the above one and\sss involves other\sss ideas.\oss
The above proof\dss is an apparently\dss unintended application of\qss
Alexander's\dss methods\qss \cite{a2}.\oss

\mysection{Sperner's\pss lemma\qss and\qss its\qss combinatorial\qss proof}{sperner-lemma}

\myuppar{Sperner\dss colorings.}
We continue\sss to use\sss the notations introduced at\dss
the beginning of\trs Section\dss \ref{alexander-lemma}.\oss
Let\qss $V\off =\off v\dff(\dff T\dff)$\qss
be\sss the set\sss of\dss vertices of\dss the\sss triangulation\dss $T$\nnsp,\oss
and\dss let\qss $V_i\off =\off v\dff(\dff T_{\dff i}\dff)$\qss
be\sss the set\sss of\dss vertices of\dss $T_{\dff i}$\nnsp.\oss
In other words,\oss $V_i\off =\off V\qff \cap\pff \delta_{\dff i}$\nnsp.\oss
A\dss map\vspace{3pt}
\[
\quad
\varphi\qff \colon\qff 
v\dff(\dff T\dff)\off \ttoo\off 
I
\off =\off
\{\qff 0\fff,\pff 1\fff,\pff \ldots\fff,\pff n \qff\}
\]

\vspace{-9pt}
is\dss said\dss to\dss be\dss a\pss \emph{Sperner\dss coloring}\dff\oss 
if\oss
$\varphi\fff(v\fff)\off \neq\off i$\oss
for every\qss $i\qff \in\qff I$\qss and\oss 
$v\qff \in\qff V_i$\nnsp.

\myuppar{Theorem\qss (Sperner\halfff's\qss lemma).}\oss
\emph{If\oss
$\varphi\dff \colon\dff 
v\dff(\dff T\dff)\ttoo I$\oss
is\dss a\dss Sperner coloring\halfff,\oss 
then the number of\trs $n$\dnsp-sim\-plices $\sigma$ of\oss $a\dff(\dff T\dff)$\qss 
such\dss that\oss $\varphi\fff(\sigma)\off =\off\qff I$\oss is\dss odd\halfff.\oss
In particular\halfff,\oss it\dss is non-zero.}

\prooftitle{A combinatorial\dss proof\halfff}
It\dss is due\sss to\dss Sperner\dss and\dss is based on a celebrated\sss
double counting\sss argument\halfff.\oss
Let us\dss fix\dss an element\qss $i\qff \in\qff I$\qss
and\dss let\dss $N$\dss be the number of\dss 
pairs\dss
$(\fff \tau\fff,\pff \sigma \fff)$\dss
such that $\tau$ is an 
$(\fff n\dff -\dff 1 \fff)$\dnsp-simplex of\dss $a\dff(\dff T\dff)$\dss 
and\qss $\varphi\dff(\dff \tau\dff)\off =\off I\qff -\qff i$\nnsp,\oss 
and\dss $\sigma$\dss is\dss an $n$\dnsp-simplex\dss of\qss $a\dff(\dff T\dff)$\dss having $\tau$ as a face.\oss
Let\dss us count\dss such pairs in\dss two ways.

The first way\dss 
is based on the 
non-branching\dss property\halfff.\oss 
Since\sss the map $\varphi$ is\sss a\dss Sperner coloring\halfff,\oss
if\oss
$\tau\off \subset\off V_k$\oss for some $k$ and\pss
$\varphi\dff(\dff \tau\dff)
\off =\off 
I\qff  -\qff i$\nnsp,\oss
then\pss $k\off =\off i$\pss and\pss
$\tau\off \subset\off V_i$\nsp.\oss
Let $h$ be the number of\dss $(\fff n\dff -\dff 1 \fff)$\dnsp-simplices\qss 
$\tau\off \subset\off   
V_i$\qss such that\pss
$\varphi\dff(\dff \tau\dff)
\off =\off 
I\qff  -\qff i$\pss
and\dss let $g$ be the number of\dss the other 
$(\fff n\dff -\dff 1 \fff)$\dnsp-simplices $\tau$
such that\pss
$\varphi\dff(\dff \tau\dff)
\off =\off 
I\qff  -\qff i$\nnsp.\oss
By the non-branching\dss property\pss\vspace*{0pt}
\[
\quad
N\off =\off h\qff +\qff 2\dff g\dff. 
\]

\vspace{-12pt}
The second way\dss is\dss independent of\dss the non-branching\dss property.\oss
If\trs an $n$\dnsp-simplex $\sigma$ has a face $\tau$ such that\pss
$\varphi\dff(\dff \tau\dff)
\off =\off 
I\qff  -\qff i$\nnsp,\oss
then either\pss 
$\varphi\dff(\dff \sigma\dff)
\off =\off
I$\nnsp,\oss
or\pss 
$\varphi\dff(\dff \sigma\dff)
\off =\off
I\qff  -\qff i$\nnsp.\qff\oss
Let\qss $e$\qss be the number of\dss 
$n$\dnsp-simplices $\sigma$ such\dss that\pss 
$\varphi\fff(\fff\sigma\fff)\off =\off I$\oss
and\dss let\dss $f$\dss be the number of\dss 
$n$\dnsp-simplices $\sigma$ such\dss that\pss 
$\varphi\dff(\dff\sigma\dff)
\off =\off 
I\qff  -\qff i$\nnsp.\qff\oss
If\pss 
$\varphi\dff(\dff\sigma\dff)
\off =\off
I$\nnsp,\oss 
then\pss
$\varphi\dff(\dff \tau\dff)
\off =\off 
I\qff  -\qff i$\pss
for exactly one face\dss $\tau$\dss if\dss $\sigma$\nnsp.\oss 
Suppose now that\pss 
$\varphi\dff(\dff\sigma\dff)
\off =\off 
I\qff  -\qff i$\nnsp.\oss
Since\pss $\num{\sigma}\off =\off n\qff +\qff 1$\pss 
and\pss $\num{\dff I\qff  -\qff i \dff}\off =\off n$\nnsp,\oss
there is a unique pair\qss $a\fff,\pff b\qff \in\qff \sigma$
such that\qss
$\varphi\dff(\dff a\dff)
\off =\off
\varphi\dff(\dff b\dff)$\qss
and\qss
$a\qff \neq\qff b$\nnsp.\oss
Clearly\halfff,\oss\vspace{2.75pt}
\[
\quad
\varphi\dff(\trf\sigma\qff \smallsetminus\qff \{\dff a \qff\}\trf)
\off =\off
\varphi\dff(\trf\sigma\qff \smallsetminus\qff \{\trf b \qff\}\trf)
\off =\off
I\qff  -\qff i
\]

\vspace{-9.25pt}
and\pss
$\num{\varphi\fff(\fff\tau\fff)}
\off <\off 
n\qff -\qff 1$\pss  
for all $(\fff n\dff -\dff 1 \fff)$\dnsp-faces $\tau$ of\dss $\sigma$ different\dss from\qss
$\sigma\qff \smallsetminus\qff  \{\dff a \qff\}$\nnsp,\pss
$\sigma\qff \smallsetminus\qff \{\trf b \qff\}$\nnsp.\oss
Hence $\sigma$ has exactly two
$(\fff n\dff -\dff 1 \fff)$\dnsp-faces such that\oss
$\varphi\dff(\dff \tau\dff)
\off =\off 
I\qff  -\qff i$\nnsp.\oss 
It\dss follows that\pss\vspace{1.5pt} 
\[
\quad
N\off =\off e\qff +\qff 2\fff f\dff. 
\]

\vspace{-12pt}
By comparing the two expressions for $N$\nnsp,\oss
we\dss see\dss that\oss\vspace{0pt} 
\begin{equation}
\label{sperner-equation}
\quad
h\qff +\qff 2\dff g
\off =\off
e\qff +\qff 2\fff f
\off.
\end{equation}

\vspace{-12pt}
Let us now use an induction by $n$\nnsp.\oss
Sperner's\dss lemma is trivially true for\qss $n\qff =\qff 0$\nnsp.\oss
Suppose that\qss $n\qff >\qff 0$\nnsp.\oss
The simplices of\dss
$a\dff (\dff T_i \dff)$\dss
are nothing else but\dss the simplices\dss of\qss
$a\dff (\dff T \dff)$\qss
contained\dss in\dss $V_i$\nsp.\qff\oss
The map $\varphi$ induces a map\qss
$\varphi_{\fff i}
\qff \colon\qff
V_i
\off \ttoo\off
I\qff  -\qff i$\nnsp.\oss
Renumbering the elements of\dss the set\pss 
$I\qff  -\qff i$\pss
by\pss
$0\fff,\pff 1\fff,\pff \ldots\fff,\pff n\qff -\qff 1$\pss
turns\dss $\varphi_{\fff i}$\dss into a Sperner coloring\halfff.\oss
Hence  the inductive assumption implies\dss that\dss the\dss number\dss of\dss
$(\fff n\dff -\dff 1 \fff)$\dnsp-simplices $\tau$ of\qss
$a\dff (\dff T_i \dff)$\qss
such that\oss
$\varphi_{\fff i}\dff(\dff \tau\dff)
\off =\off 
I\qff  -\qff i$\oss 
is odd.\qff\oss
But\dss this number is equal\dss to $h$\nnsp.\qff\oss
Therefore the equality\qss (\ref{sperner-equation})\qss
implies that\dss $e$\dss is odd.\qff\oss
This completes the induction step and\dss hence the proof\halfff.\oss  \eproof

\myuppar{Remark.}
Strictly\sss speaking\halfff,\oss
there is no such statement\dss in\dss Sperner's\dss paper\qss \cite{s},\oss
but\dss there is\sss the same proof\halfff,\oss up\sss to\sss the language and\dss notations.\oss
Sperner\dss starts with a closed\sss covering of\dss $\delta$\dss by\qss
$n\qff +\qff 1$\qss sets satisfying\dss the assumptions of\qss Lebesgue-Sperner\dss theorem\qss 
(see\dss Section\qss \ref{dimension})\qss
and chooses\dss $\varphi$\dss in\dss the same way as in\dss its\sss proof\halfff.\oss
\emph{``Sperner's\dss lemma''}\oss appeared\dss for\dss the first\dss time in\dss
Knaster-Kuratowski-Mazurkiewicz\dss
paper\qss \cite{kkm}\qss 
as\dss the\qss \emph{``combinatorial core of\qss Sperner's\dss new\sss proof\dss
of\dss the invariance of\dss dimension''}.\oss
The authors\dss of\pss \cite{kkm}\qss also modified\dss the proof\halfff.\oss
They considered only\dss the numbers\dss $e$\dss and\dss $h$\dss
and\dss worked\dss modulo\dss $2$\nnsp,\oss
in contrast\dss with\dss Sperner\halfff.\oss

\myuppar{Remark.}
The second\dss method\sss of\dss counting\sss is parallel\dss to\sss the proof\dss
of\trs Theorem\qss 1.\oss
The case\pss 
$\varphi\fff(\fff\sigma\fff)\off =\off I$\oss
corresponds\sss to\dss Case\dss 1\dss of\trs that\dss proof\halfff,\oss
and\dss the case\pss 
$\varphi\dff(\dff\sigma\dff)
\off =\off 
I\qff  -\qff i$\oss
corresponds\sss to\dss Case\qss 3.\oss
The case when\dss the dimension of\dss $\varphi\dff(\dff \sigma\dff)$\dss
is\qss $\leq\qff n\qff -\qff 2$\qss
would correspond\dss to\dss Case\qss 2,\oss
but\sss such simplices\dss $\sigma$\dss
do not\sss occur\dss in\dss pairs\dss
$(\fff \tau\fff,\pff \sigma \fff)$\dss
such\dss that\qss $\varphi\dff(\dff \tau\dff)\off =\off I\qff -\qff i$\nnsp.\oss

\myuppar{Sperner\dss colorings as simplicial\dss maps.}
Recall\dss that\dss $\Delta$\dss is\dss the simplicial complex consisting of\dss
simplex\dss $\delta$\dss and all\dss its faces.\oss
Let us identify the vertices\dss $v_{\fff i}$\dss of\dss $\delta$\dss 
with\dss their subscripts\qss $i\qff \in\qff I$\nsp.\oss
This\sss turns a\dss Sperner\dss colorings\dss $\varphi$\dss into a map\qss
$v\dff(\dff T\dff)\ttoo
v\dff(\dff \Delta\dff)\off =\off I$\nnsp.\oss
Since every subset\sss of\qss $v\dff(\dff \Delta\dff)$\dss
is a simplex\dss of\qss $a\dff(\dff \Delta\dff)$\nnsp,\oss
this\dss is\sss a simplicial\dss map\qss
$T\ttoo \Delta$\nnsp.\oss
A simplicial\dss map\qss
$\varphi
\dff \colon\dff
T
\ttoo
\Delta$\qss
is\sss a\dss Sperner\dss coloring\trs if\trs and\sss only\dss if\qss
for every\qss $i\qff \in\qff I$\qss it\dss 
takes each vertex of\qss $T$\dss 
belonging\dss to\dss $\delta_{\dff i}$\dss
into an element\dss of\qss $I\qff -\qff i$\nnsp,\oss
i.e.\qss into\sss a\dss vertex of\dss $\delta_{\dff i}$\nnsp.\oss

Since every\dss proper face of\dss $\delta$\dss is equal\dss to\sss the intersection of\dss
several\dss $(\fff n\dff -\dff 1\fff)$\dnsp-dimensional\dss 
faces\dss $\delta_{\dff i}$\nsp,\oss 
this condition implies\sss that\dss $\varphi$\dss
maps\sss the set\sss of\dss vertices of\dss $T$\dss belonging\dss to a face\dss $\tau$\dss of\dss $\delta$\dss
into\sss the set\sss of\dss vertices of\dss $\tau$\nnsp.\oss
It\dss follows\sss that\qss 
$\varphi\dff \colon\dff
v\dff(\dff T\dff)\ttoo
v\dff(\dff \Delta\dff)$\qss 
is a\dss Sperner\dss coloring\dss
if\dss and only\dss if\dss $\varphi$\dss is a pseudo-identical simplicial\dss map\qss
$T\ttoo \Delta$\nnsp.\oss

Now\qss Sperner's\trs lemma\dss takes\sss the following\sss form.\oss

\myitpar{The simplicial\dss form\sss of\qss Sperner's\dss lemma.}
\emph{If\pss
$\varphi\dff \colon\dff 
T\qff \ttoo\qff \Delta$\qss
is\dss a\dss pseudo-identical\dss simplicial\dss map,\oss 
then the\dss number of\dss $n$\dnsp-simplices\dss $\sigma$\dss of\oss
$a\dff(\dff T\dff)$\dss
such\dss that\pss
$\varphi\dff(\dff \sigma\dff)\off =\off v\dff(\dff \Delta \dff)$\qss 
is\dss odd\halfff.}\vspace{6pt}

This\sss form of\qss Sperner's\trs lemma\sss is an\sss immediate corollary\dss
of\qss Alexander's\qss lemma,\oss i.e.\qss Theorem\qss 3.\qss
Indeed,\pss Theorem\dss 3\dss implies\sss that\dss
$\varphi_*\dff(\dff \fclass{\delta} \dff)
\off =\off
\delta$\nnsp.\pss
But\dss
by\dss the definition of\qss induced\dss maps\qss
$\varphi_*\dff(\dff \fclass{\delta} \dff)
\off =\off
e\trf \delta$\nnsp,\oss
where\dss $e$\dss is\dss the number of\dss $n$\dnsp-simplices\dss $\sigma$\dss
of\dss $T$\dss
such\dss that\qss 
$\varphi\dff(\dff \sigma\dff)\off =\off \delta$\nnsp.\oss
Therefore\dss $e$\dss is\sss equal\dss to\dss $1$\dss in\dss $\ftwo$\nnsp,\qff\oss
i.e.\dss $e$\dss is\sss odd.\oss

\mysection{Cochains\pss and\pss Sperner's\pss lemma}{dualizing}

\myuppar{The combinatorial\dss proof\dss and\sss algebraic\sss topology.}
The goal of\dss this section is\sss to show\dss that\dss not\sss only\dss
Sperner's\dss lemma admits a natural\dss topological\dss interpretation,\oss
but\dss its combinatorial\dss proof\dss is also a\sss topological\dss proof\dss
in disguise.\oss
The proof\dss of\qss Alexander's\dss lemma\qss (i.e.\qss of\qss Theorem\qss 3)\qss
depends on\trs Theorem\qss 1\qss and\dss hence
indirectly\sss contains a part\sss of\dss the combinatorial\dss proof\qss
(see\sss the second\dss remark\sss in\dss Section\qss \ref{sperner-lemma}).\oss
Also,\oss the induction\dss by\dss $n$\dss is used\sss in\dss both\dss proofs\sss
in an essentially\dss the same manner\halfff.\oss
Still,\oss the proofs look quite different\halfff.\oss

The induced\dss maps\dss $\varphi_*$\sss are a good\dss tool\dss to deal
with the images of\dss simplices under a simplicial\dss map\dss $\varphi$\dnsp.\oss
But\dss the combinatorial\dss proof\dss of\qss Sperner's\dss lemma
operates not with\dss the images but with\dss the\qss \emph{preimages},\oss
the sets of\dss simplices\dss mapped\dss by\dss $\varphi$\sss 
to\dss particular\dss simplices\dss
in\dss the\dss target\dss complex\halfff,\oss
namely\halfff,\pss  
to the simplices $I$ and\qss $I\qff -\qff i$\qss
of\dss $a\dff(\dff \Delta\dff)$\nnsp.\oss

This suggests to dualize the notions of\dss induced\dss maps
and\dss boundary operators in the sense of\dss the linear algebra
over the field\dss $\ftwo$\dss and\dss leads\sss
to\sss the notion of\dss cochains.\oss
The reader should\dss keep in\sss mind\dss that\dss this motivation\sss
is\sss an\sss artificial\sss one.\oss
The cochains were introduced\sss in\dss 1935\dss
by\dss completely different\dss reasons independently\dss
by\dss Alexander\qss \cite{a3},\pss \cite{a4}\qss and\dss Kolmogoroff\qss \cite{ko}.

\myuppar{Cochains.}
For a simplicial complex\dss $S$\dss and 
a non-negative integer\dss $m$\dss
let\oss\vspace{1.5pt} 
\[
\quad
C^{\fff m}\fff(\dff S \dff)
\off\qff =\off\qff
C_{\fff m}\fff(\dff S \dff)^*
\] 

\vspace{-12pt}
be the vector space dual\dss over\sss $\ftwo$\sss 
to\sss $C_{\fff m}\fff(\dff S \dff)$\dnsp.\qff\oss
Its elements
are called \emph{$m$\dnsp-cochains}\qss of\dss $S$\nnsp.\qff\oss
Let\vspace{3pt}
\[
\quad
\partial^*
\qff \colon\qff
C^{\fff m\dff -\dff 1}\fff(\dff S \dff)
\off \ttoo\off
C^{\fff m}\fff(\dff S \dff)
\]

\vspace{-10.5pt}
be the\sss linear map dual\dss to the boundary operator\pss
$\dis
\partial
\dff \colon\dff
C_{\fff m}\fff(\dff S \dff)
\ttoo
C_{\fff m\dff -\dff 1}\fff(\dff S \dff) $\dnsp.\oss
The map\dss $\partial^*$\dss is called\dss the\qss \emph{coboundary\sss operator\halfff.}\oss
For a simplicial\dss map\pss 
$\varphi\dff \colon\dff S\ttoo S'$\pss
the\qss \emph{induced\dss map}\qss\vspace{3pt}
\[
\quad
\varphi^*
\qff \colon\qff
C^{\fff m}\fff(\dff S' \dff)
\off \ttoo\off
C^{\fff m}\fff(\dff S \trf)
\]

\vspace{-10.5pt}
is defined as the linear map dual\dss to\sss the induced map\pss
$\dis
\varphi_*
\dff \colon\dff
C_{\fff m}\fff(\dff S \dff)
\ttoo
C_{\fff m}\fff(\dff S' \trf)$\dnsp.\oss

\myuppar{Cochains as formal sums of\dss simplices.}
Since\dss $C_{\fff m}\fff(\dff S \dff)$\dss
is a vector space over\dss $\ftwo$\dss having a canonical\dss basis 
consisting\sss of\dss $m$\dnsp-simplices of\dss $S$\nnsp,\oss
the $m$\dnsp-cochains can\dss be identified\dss with\dss
$\ftwo$\dnsp-valued\dss functions on\dss the set\sss of\dss
$m$\dnsp-simplices of\dss $S$\nnsp.\oss
Since all\sss our complexes are assumed\dss to be finite,\oss
this basis is finite and can\sss be used\dss to identify\dss
the vector space\dss $C_{\fff m}\fff(\dff S \dff)$\dss
with\sss its dual\dss $C^{\fff m}\fff(\dff S \dff)$\dss
and\dss interpret\sss cochains,\pss
like chains,\pss 
as formal sums of\dss simplices.\oss

In what\dss follows,\pss we\sss write
cochains as formal sums of\dss simplices,\oss
but\dss keep\sss the notation\dss
$C^{\fff m}\fff(\dff S \dff)$\dss
as an indicator showing\dss that\dss we\sss treat\dss these formal sums as cochains.\oss
The identification of\dss cochains with formal sums of\dss simplices\sss
turns\sss the maps\dss
$\partial^*$\dss and\dss $\varphi^*$\dss into the\qss 
\emph{adjoint\dss operators}\pss of\qss $\partial$\dss and\dss $\varphi_*$\dss
respectively with respect\dss to the pairings\qss
$\langle\dff \bullet\fff,\qff \bullet\dff\rangle$\qss such\dss that\vspace{3pt}
\[
\quad
\langle\dff \sigma\fff,\qff \tau\dff\rangle
\off =\off
1
\hspace*{0.8em}\mbox{if}\hspace*{0.9em}
\sigma
\off =\off
\tau
\hspace*{1.5em}\mbox{and}\hspace*{1.5em}
\langle\dff \sigma\fff,\qff \tau\dff\rangle
\off =\off
0
\hspace*{0.8em}\mbox{if}\hspace*{0.9em}
\sigma
\off \neq\off
\tau
\qff.
\]

\vspace{-9pt}
A\sss trivial\sss verification shows\sss that\dss
if\dss $\tau$\sss is\dss an 
$(\fff m\dff -\dff 1 \fff)$\dnsp-simplex\sss of\dss $S$\dnsp,\oss
then\dss\vspace{3pt}
\begin{equation}
\label{coboundary}
\quad
\partial^*\dff(\dff \tau \trf)
\off =\qff\off
\sum\qff \sigma\qff,
\end{equation}

\vspace{-9pt}
where the sum is taken over\dss all $m$\dnsp-simplices $\sigma$ having $\tau$ as a face.\oss
Hence\dss the coboundary operator\dss $\partial^*$\nnsp,\oss
like\sss the boundary operator\dss $\partial$\nnsp,\oss
encodes\sss the relation\oss 
\emph{``\nsp$\tau$\qss is\dss a\dss face\dss of\qss $\sigma$\dnsp''}\oss
between\sss simplices\qss $\tau\fff,\pff \sigma$\qss
such\dss that\dss the dimension of\dss $\tau$\dss is less by\dss $1$\dss
than\dss the dimension of\dss $\tau$\nnsp.\oss

Similarly\halfff,\oss
if\pss 
$\varphi\dff \colon\dff S\ttoo S'$\pss
is a simplicial\dss map and\dss
$\rho$\sss is an $m$\dnsp-simplex of\trs $S'$\nnsp,\oss
then\dss\vspace{3pt} 
\begin{equation}
\label{coinduced}
\quad
\varphi^*\dff(\dff \rho \dff)
\off =\qff\off
\sum\qff \tau\qff,
\end{equation}

\vspace{-9pt}
where the sum is over\dss all $m$\dnsp-simplices\dss $\tau$\dss of\trs $S$\dss  
such that\oss
$\varphi\dff(\dff \tau \dff)
\off =\off
\rho$\nnsp,\qff\oss
as another\dss trivial verification shows.\oss
In other words,\oss
$\varphi^*\dff(\dff \rho \dff)$\dss
indeed encodes
the preimage of\dss $\rho$\dnsp.\oss\vspace{0pt}

\myuppar{Theorem\qss 1*.}\oss
$\partial^*\qff \circ\qff \varphi^*
\off =\off
\varphi^*\qff \circ\qff \partial^*$\dnsp.\oss\vspace{0pt}

\proof
This\sss immediately\dss follows from\trs Theorem\qss 1\qss by dualizing\halfff.\oss
One can also give a direct\dss proof\qss
based on\qss (\ref{coboundary})\qss and\qss (\ref{coinduced}).\oss
We leave\sss this\sss task\sss to\sss the interested\dss
readers as an exercise.\oss  \eproof\vspace{0pt}

\prooftitle{A cochains-based\dss proof\qss  
of\dss the simplicial\qss form\dss of\pss Sperner's\trs lemma\fff}
We are going\dss to partially\sss dualize\qss
the\qss proof\trs of\qss Alexander's\trs lemma.\oss
The latter\dss is based on\dss Theorems\qss 1\qss and\qss 2,\oss 
the equality\qss (\ref{boundary-delta}),\oss
and\qss
Lemma\dss from\dss Section\qss \ref{alexander-lemma},\oss
a compressed\dss form of\dss the non-branching\dss property.\oss
The dualization of\qss Theorem\qss 1\qss is\qss Theorem\qss 1*.\oss
Theorem\qss 2\qss and\dss
the equality\qss (\ref{boundary-delta})\qss cannot\dss be straightforwardly\sss dualized,\oss
but\trs if\trs we\sss fix some\qss $i\qff \in\qff I$\nnsp,\oss
then\dss the obvious equality\vspace{2.875pt}
\begin{equation}
\label{delta-coboundary}
\quad
\partial^*\dff \bigl(\dff \delta_{\dff i}\dff\bigr)
\off =\off\qff 
\delta
\dff
\end{equation}

\vspace{-9.125pt}
turns out\dss to\sss be a reasonable substitution for\qss (\ref{boundary-delta}).\oss
By\qss Theorem\qss 1*\vspace{3pt}
\[
\quad
\varphi^*\left(\qff \partial^*(\dff \delta_{\dff i}\dff) \dff\right)
\off =\off
\partial^*\left(\qff \varphi^*(\dff \delta_{\dff i}\dff) \dff\right)
\qff,
\]

\vspace{-9pt}
and\dss together\sss with\qss (\ref{delta-coboundary})\qss 
this implies\sss that\vspace{3pt} 
\begin{equation}
\label{pull-back-delta}
\quad
\varphi^*\fff (\dff \delta\dff)
\off =\off
\partial^*\dff \bigl(\qff \varphi^*\dff (\dff \delta_{\dff i} \dff) \qff\bigr)
\qff.
\end{equation}

\vspace{-9pt}
Instead of\qss Lemma\dss from\dss Section\qss \ref{alexander-lemma}\qss
we will\dss use\sss the non-branching\dss property\sss directly.\oss
Let\dss us
explicitly compute\sss the cochains\dss in\qss (\ref{pull-back-delta})\qss
and\dss relate\sss them\sss to\sss the numbers\qss
$e\fff,\pff f\fff,\pff g\fff,\pff h$\pss
from\dss the combinatorial\dss proof\halfff.\oss
We need\dss the following\dss four sets of\dss simplices.\oss\vspace{0pt}

Let\dss $E$\dss and\dss $F$\dss be\sss the sets of\dss
$n$\dnsp-simplices\dss $\sigma$\dss of\dss $T$\dss
such\dss that\qss
$\varphi\dff(\dff \sigma\dff)\off =\off \delta$\qss
and\qss
$\varphi\dff(\dff \sigma\dff)\off =\off \delta_{\dff i}$\qss
respectively.\oss
The sets\dss $E$\dss and\dss $F$\dss
consist\sss of\dss $e$\dss and\dss $f$\dss 
elements respectively.\oss

Let\dss $H$\dss be\sss the set of\dss
$(\fff n\dff -\dff 1\fff)$\dnsp-sim\-plices\dss $\tau$\dss of\trs $\bd\dff T$\dss
such\dss that\qss
$\varphi\dff(\dff \sigma\dff)\off =\off \delta_{\dff i}$\nsp.\oss
Since\dss $\varphi$\dss is a pseudo-identical\sss simplicial\dss map,\pss
every such simplex\dss $\tau$\dss is actually a simplex of\dss $T_{\dff i}$\nnsp.\oss
In\dss particular\halfff,\pss $H$\dss consists of\dss $h$\dss elements.\oss

Finally,\oss let\dss $G$\dss be\sss the set of\dss 
$(\fff n\dff -\dff 1\fff)$\dnsp-simplices\dss $\tau$\dss of\dss $T$\dss
such\dss that\qss
$\varphi\dff(\dff \tau\dff)\off =\off \delta_{\dff i}$\nsp,\oss
but\dss $\tau$\dss is not a simplex of\trs $\bd\dff T$\nnsp,\oss
i.e.\dss $\tau\qff \not\in\qff H$\nnsp.\oss 
The set\dss $G$\dss consists of\dss $g$\dss elements.\oss\vspace{3pt}

In\dss terms of\dss the sets\qss $E\fff,\off F\fff,\off G\fff,\off H$\qss
the cochains\dss $\varphi^*\fff (\dff \delta\dff)$\dss
and\dss $\varphi^*\fff (\dff \delta_{\dff i}\dff)$\dss
can\dss be written as follows\fff:\vspace{4.5pt}
\[
\quad
\varphi^*\fff (\dff \delta\dff)
\off\off =\qff\off\off
\sum\nolimits_{\qff \sigma\qff \in\qff E}\dff \sigma
\]

\vspace{-16.5pt}
and\vspace{-4.5pt}
\[
\quad
\varphi^*\fff (\dff \delta_{\dff i}\dff)
\off\off =\qff\off\off
\sum\nolimits_{\qff \tau\qff \in\qff G}\dff \tau
\off\off +\qff\off\off
\sum\nolimits_{\qff \tau\qff \in\qff H}\dff \tau
\off.
\]

\vspace{-7.5pt}
Therefore we can rewrite\qss (\ref{pull-back-delta})\qss as\vspace{4.5pt}
\begin{equation}
\label{unfolded}
\quad
\sum\nolimits_{\qff \sigma\qff \in\qff E}\dff \sigma
\qff\off\off =\qff\off\off
\sum\nolimits_{\qff \tau\qff \in\qff G}\qff \partial^*\dff(\dff \tau\dff)
\off\off +\qff\off\off
\sum\nolimits_{\qff \tau\qff \in\qff H}\qff \partial^*\dff(\dff \tau\dff)
\off.
\end{equation}

\vspace{-7.5pt}
By\dss the\sss non-branching\trs property\halfff,\qff\oss
if\pss $\tau\qff \in\qff H$\nnsp,\oss
then\dss $\tau$\dss is a face of\dss exactly one $n$\dnsp-simplex of\trs $T$\nnsp,\oss
and\dss if\pss $\tau\qff \in\qff G$\nnsp,\oss
then\dss $\tau$\dss is a face of\dss exactly\dss two $n$\dnsp-simplices of\trs $T$\nnsp.\oss
In\dss terms of\dss the coboundary operator\dss $\partial^*$\dss this means\sss that\dss
if\pss $\tau\qff \in\qff H$\nnsp,\oss
then\dss $\partial^*\dff \tau$\dss is a simplex,\oss
and\qss if\pss $\tau\qff \in\qff G$\nnsp,\oss
then\dss $\partial^*\dff \tau$\dss is a sum of\trs two simplices.\oss 
Hence\dss the\dss right\trs hand side of\pss (\ref{unfolded})\pss
is\dss a\dss sum\sss of\pss $h\qff +\qff 2\dff g$\qss simplices.\oss

If\dss some $n$\dnsp-simplex $\sigma$ occurs in\dss this sum at\dss least\dss twice,\oss
then $\sigma$ has at\dss least\dss two $(\fff n\dff -\dff 1 \fff)$\dnsp-faces\qss
$\tau\fff,\pff \tau'$\qss such\dss that\qss
$\varphi\dff(\dff \tau\dff)
\off =\off
\varphi\dff(\dff \tau'\dff)
\off =\off
\delta_{\dff i}$\nnsp.\oss
In\dss this case\qss
$\varphi\dff(\dff \sigma\dff)
\off =\off
\delta_{\dff i}$\qss
and\qss
$\varphi\dff(\dff \tau''\dff)
\off \neq\off
\delta_{\dff i}$\qss
for any other\dss face\dss $\tau''$\dss of\dss $\sigma$\dnsp.\oss
Therefore,\oss
in\dss this case\qss
$\sigma\qff \in\qff F$\qss
and\dss $\sigma$\dss occurs\sss in\dss the sum exactly\dss two\sss times.\oss
Conversely,\oss if\pss
$\sigma\qff \in\qff F$\nnsp,\oss
then\dss $\sigma$\dss has\sss two such\dss faces
and\dss hence $\sigma$ occurs in\dss this sum\dss twice.\oss
In other words,\oss
pairs of\dss equal simplices at\dss the
right hand side of\qss (\ref{unfolded})\qss 
correspond\dss to elements of\trs $F$\dss 
and\dss
there are $f$ such pairs.\oss
Over\dss $\ftwo$\dss
such\dss pairs cancel\halfff.\oss

There are no other cancellations
and\dss hence\sss the right\dss hand side of\qss (\ref{unfolded})\qss
is equal\dss to\sss a sum of\qss
$h\qff +\qff 2\dff g\qff -\qff 2\fff f$\qss distinct\sss simplices.\oss
At\dss the same\sss time\sss 
the left\dss hand side of\qss (\ref{unfolded})\qss
is obviously a sum of\dss $e$\dss distinct\sss simplices.\oss
Therefore\qss (\ref{unfolded})\qss implies\sss that\vspace{3.75pt} 
\begin{equation}
\label{cochains-efgh-equation}
\quad
e\off =\off h\qff +\qff 2\dff g\qff -\qff 2\fff f
\dff.
\end{equation}

\vspace{-8.25pt}
It\dss follows\dss that\qss
$e\qff \equiv\qff h$\qss modulo $2$\nnsp.\oss
Now one can\sss use
induction\dss by $n$
to complete\sss the proof\halfff.\oss  \eproof\dnsp

\prooftitle{An alternative\dss ending\halfff}
After\dss the equality\qss (\ref{pull-back-delta})\qss is proved,\oss
one can use\sss the non-branching\dss property\dss in\dss the form of\qss
Lemma\dss from\dss Section\qss \ref{alexander-lemma}.\oss
This leads\sss to a proof\dss closer\dss to\dss Alexander's\dss one.\oss
By\dss pairing\dss both sides of\qss (\ref{pull-back-delta})\qss
with\dss $\fclass{\delta}$\dss  
we see\sss that\vspace{4.5pt}
\begin{equation*}
\quad
\left\langle\qff 
\varphi^*\fff (\dff \delta\dff)\fff,\off
\fclass{\delta}
\qff\right\rangle
\off\qff =\off\qff
\left\langle\qff 
\partial^*\dff \bigl(\qff \varphi^*\dff (\dff \delta_{\dff i} \dff) \qff\bigr)\fff,\off 
\fclass{\delta}
\qff\right\rangle
\qff.
\end{equation*}

\vspace{-7.5pt}
Together\dss with\dss the definition of\dss $\partial^*$\dss
this implies\sss that\vspace{4.5pt}
\begin{equation*}
\quad
\left\langle\qff 
\varphi^*\fff (\dff \delta\dff)\fff,\off
\fclass{\delta}
\qff\right\rangle
\off\qff =\off\qff
\left\langle\qff 
\varphi^*\dff (\dff \delta_{\dff i} \dff)\fff,\off 
\partial\dff \fclass{\delta}
\qff\right\rangle
\qff.
\end{equation*}

\vspace{-7.5pt}
Now\dss the equality\qss
$\partial\dff \fclass{\delta}
\off =\off
\fclass{\dff\partial\dff \delta\dff}$\qss
of\pss Lemma\dss from\dss Section\qss \ref{alexander-lemma}\qss
implies\sss that\vspace{4.5pt}
\begin{equation}
\label{dual-basic}
\quad
\left\langle\qff 
\varphi^*\fff (\dff \delta\dff)\fff,\off
\fclass{\delta}
\qff\right\rangle
\off\qff =\off\qff
\left\langle\qff 
\varphi^*\dff (\dff \delta_{\dff i} \dff)\fff,\off 
\fclass{\partial\dff \delta}
\qff\right\rangle
\qff.
\end{equation}

\vspace{-7.5pt}
Pairing cochains with\dss $\fclass{\delta}$\dss and\dss
$\fclass{\partial\dff \delta}$\dss amounts\sss to counting\dss
their simplices modulo $2$\nnsp.\oss
In\dss more details,\oss
since we are working over\dss $\ftwo$\nnsp,\oss
any cochain can\dss be written as a sum of\dss several\qss
\emph{distinct}\qss simplices.\oss 
Obviously,\oss
if\dss $\alpha$\dss is an\dss $n$\dnsp-cochain of\dss $T$\nnsp,\oss
then\dss
$\left\langle\qff \alpha\fff,\pff \fclass{\delta\fff} \qff\right\rangle$\dss
is equal\dss to\sss number of\dss simplices of\trs $T$\dss
in\dss the sum\dss $\alpha$\dss taken \dss modulo\dss $2$\nnsp.\oss
Similarly,\oss
if\qss $\beta$\dss is an\dss $(\fff n\dff -\dff 1 \fff)$\dnsp-cochain,\oss
then\sss
$\left\langle\qff \beta\fff,\pff \fclass{\dff\partial\dff \delta\dff} \qff\right\rangle$\sss
is equal\dss to\sss number of\dss simplices of\trs $\bd\dff T$\dss
in\dss the sum\dss $\beta$\dss taken \dss modulo\dss $2$\nnsp.\oss

It\dss follows\sss that\dss
the right\dss hand side of\qss (\ref{dual-basic})\qss is equal\dss to\sss
the\sss taken\dss modulo\dss $2$\dss
number of\dss simplices\dss $\tau$\dss of\trs $\bd\dff T$\dss
such\dss that\qss
$\varphi\dff(\dff \tau\dff)\off =\off \delta_{\dff i}$\nsp.\oss
But\sss since\dss $\varphi$\dss is a pseudo-identical\dss simplicial\dss map,\oss
every such simplex $\tau$ is a simplex of\trs $T_{\fff i}$\nnsp.\oss
Hence\sss the right\dss hand side of\qss (\ref{dual-basic})\qss is equal\dss to
$h$ modulo\dss $2$\nnsp.\oss
Similarly,\oss
the left\dss hand side of\qss (\ref{dual-basic})\qss is equal\dss to
$e$ modulo\dss $2$\nnsp.\oss
Now\qss (\ref{dual-basic})\qss
implies\sss that\qss
$e\off \equiv\off h$\qss modulo\dss $2$\dss
and\dss one can complete\sss the proof\trs by\sss
using an\dss induction\dss by\dss $n$\nnsp.\oss  \eproof\vspace{3pt}

\myuppar{The cochains-based\dss proofs\dss and\dss the combinatorial\dss proof\halfff.}
The equality\qss
(\ref{cochains-efgh-equation})\qss
from\dss the co\-chains-based\dss proof\dss
is\dss trivially equivalent\dss to\sss the equality\qss (\ref{sperner-equation})\qss
around\dss
which\dss the combinatorial\dss proof\dss is centered.\oss
The equality\qss (\ref{unfolded})\qss between cochains\sss is\sss a\sss
realization\qss (or a\sss lift\dss to\sss the linear algebra)\qss
of\trs the equalities\qss
(\ref{cochains-efgh-equation})\qss and\qss (\ref{sperner-equation})\qss
between\sss numbers,\oss
and\dss the whole\sss cochains-based\dss proof\dss 
is essentially\dss
the combinatorial\dss proof\dss
rewritten in\dss the spirit\sss of\dss the linear algebra methods in combinatorics.\oss
But\dss from\dss the\sss point\sss of\dss view of\dss a\sss topologist\dss
both\dss these proofs are hardly satisfactory,\oss
in contrast\dss with\qss Alexander's\dss one.\oss

On\dss the one hand,\pss the numbers\qss $f\nsp,\pff g$\qss
and\dss the number of\dss cancellations
are irrelevant\dss to\sss the problem at\dss hand.\oss
The alternative version of\trs the cochains-based\dss proof\dss
is\dss better\dss in\dss this respect\qss
(at\dss the cost\sss of\dss being\dss further\dss from\dss Sperner's\dss one).\oss
While\dss the numbers of\dss interest\qss
$e\fff,\pff h$\qss
naturally appear\sss in\dss the proof\halfff,\oss
in\dss this version\dss
the numbers\qss $f\nsp,\pff g$\qss 
are hidden\dss by\dss the equality\qss 
$\partial\dff \fclass{\delta}
\off =\off
\fclass{\dff\partial\dff \delta\dff}$\nnsp.\oss
On\dss the other hand,\pss these proofs ignore a fundamental\dss
property of\trs $T$\dss and\dss $\bd\dff T$\nnsp,\oss 
namely,\oss
the fact\dss that\dss their\dss top-dimensional\qss
\emph{cohomology\dss groups}\pss
are isomorphic\sss to\dss $\ftwo$\nnsp.\oss
Passing\dss from cochains\sss to cohomology classes and\dss using\dss this fact\sss
allows\sss to clarify\dss the proof\dss and
carry out\dss the counting\sss in a more natural\sss way\dss
than\dss pairing\sss cochains with\dss $\fclass{\delta}$\dss and\dss
$\fclass{\dff\partial\dff \delta\dff}$\nnsp.\oss

\mysection{Graphs\qss and\pss path-following\qss algorithms}{graphs}

\myuppar{Graph-theoretical\dss interpretation of\trs the 
cochains-based\dss proof\halfff.}
There\dss is\dss a\sss widespread opinion\dss that\dss the classical\dss 
proofs\dss of\pss Sperner's\trs lemma\sss are\sss pure\sss existence\sss proofs.\oss
In\dss fact\halfff,\pss an analysis of\trs these proofs
naturally\dss leads\sss to an algorithm\dss leading\dss to a simplex $\sigma$
such\dss that\qss
$\varphi\dff(\dff \sigma\dff)\off =\off \delta$\nnsp.\oss
Let\dss us\dss begin\dss with such an analysis of\trs one step of\dss induction\dss in\dss the 
cochains-based\dss proof\halfff.\oss
In\dss this analysis 
we will\dss freely\dss use\sss the notations introduced\dss
in\dss this proof\halfff.\oss\vspace{-0.375pt}

The equality\qss (\ref{cochains-efgh-equation})\qss was proved\dss by\sss
counting\dss the number of\dss cancellations in\dss the equality\qss (\ref{unfolded}).\oss
But\dss this counting\dss was\sss based on determining\dss what\sss simplices do
actually\sss cancel\halfff.\oss
A convenient\dss way\dss to record\dss this more detailed\dss information\dss is\dss
to introduce an appropriate graph\dss $\mathbb{G}_{\dff i}$\nsp,\oss
where\qss $i\qff \in\qff I$\qss is\dss the element\dss fixed at\dss the beginning of\trs
the proof\halfff.\oss
This\sss graph\dss 
has\sss two kinds of\dss vertices.\oss
The vertices of\trs the first\dss kind are
$(\fff n\dff -\dff 1\fff)$\dnsp-simplices belonging\dss to\sss the union\qss
$G\qff \cup\qff H$\nnsp.\oss
The vertices of\trs the second\dss kind are
$n$\dnsp-simplices belonging\dss to\sss the union\qss
$E\qff \cup\qff F$\nnsp.\oss
A vertex\qss $\tau\qff \in\qff G\qff \cup\qff H$\qss
is\dss connected\dss to a vertex\qss $\sigma\qff \in\qff E\qff \cup\qff F$\qss
if\dss $\tau$\dss is\dss an
$(\fff n\dff -\dff 1\fff)$\dnsp-face of\trs $\sigma$\nnsp.\oss
There are no other edges.\oss

The graph\dss $\mathbb{G}_{\dff i}$\dss encodes all\dss relevant\dss information about\dss
the equality\qss (\ref{unfolded}).\oss
Indeed,\oss the elements of\qss $G\qff \cup\qff H$\qss
correspond\dss to\sss the summands at\dss the right\dss hand side of\qss
(\ref{unfolded}).\oss
The elements of\trs $E$\dss correspond\dss to\sss 
the summands at\dss the\sss left\dss hand side of\qss
(\ref{unfolded}),\oss
and\dss the elements of\trs $F$\dss correspond\dss to\sss 
the cancellations at\dss the right\dss hand side.\oss
Finally\halfff,\oss a\sss vertex\qss $\tau\qff \in\qff G\qff \cup\qff H$\qss
is\dss connected\dss to a vertex\qss $\sigma\qff \in\qff E\qff \cup\qff F$\qss
if\trs and\dss only\trs if\dss
$\sigma$\dss is\dss a\sss summand\sss of\trs the coboundary\dss
$\partial^*\dff(\dff \tau\dff)$\nnsp.\oss
 
The main\dss properties of\trs $\mathbb{G}_{\dff i}$\sss are\sss the following\halfff.\oss
By\dss the non-branching\dss property\halfff,\oss
every\dss vertex\sss in\sss $G$\sss is\dss an endpoint\sss of\dss exactly\dss two edges,\oss
and every\dss vertex\dss in\sss $H$\sss is an endpoint\sss of\dss exactly\sss one edge.\oss
Clearly\halfff,\oss every $n$\dnsp-simplex\dss in\sss $E$\sss has exactly\sss one face
belonging\dss to\qss $G\qff \cup\qff H$\qss and\dss hence\dss is\dss an endpoint\sss
of\dss exactly\dss one edge.\oss
Also,\oss 
every $n$\dnsp-simplex\dss in\sss $F$\sss has exactly\dss two faces
belonging\dss to\qss $G\qff \cup\qff H$\qss and\dss hence\dss is\dss an endpoint\sss
of\dss exactly\dss two edges.\oss
In\dss particular\halfff,\oss every\dss vertex\sss of\trs $\mathbb{G}_{\dff i}$\sss
is\dss an endpoint\sss of\dss either one or\dss two edges.\oss
It\dss follows\dss that\trs $\mathbb{G}_{\dff i}$\sss
consists of\dss several\sss disjoint\dss paths and cycles.\oss
Clearly\halfff,\oss every\dss path connects\sss two vertices in\dss the union\qss
$E\qff \cup\qff H$\nnsp,\oss
and every\sss vertex\dss in\dss this union\dss is\dss an endpoints of\dss a\sss path.\oss
Therefore\dss the number of\dss elements of\qss
$E\qff \cup\qff H$\nnsp,\oss i.e.\qss $e\qff +\qff h$\nnsp,\oss is\dss even
and\dss hence\qss
$e\off \equiv\off h$\qss modulo $2$\nnsp.\oss
We see\sss that\dss the inductive step 
in\dss the cochains-based\dss proof\dss can\dss be rephrased\sss
in\dss terms of\trs the graph\dss $\mathbb{G}_{\dff i}$\nsp.\oss\vspace{-0.375pt}

\myuppar{The\dss graphs\trs $\mathbb{G}_{\dff i}$\dss in\dss the combinatorial\dss proof\halfff.}
In order\dss to match\dss the above discussion,\oss
let\dss us\dss switch\dss from\dss the abstract\sss
complex $a\dff(\dff T\dff)$\sss to\sss the geometric complex\dss $T$\dnsp.
The combinatorial\dss proof\trs is\dss based\dss on counting\dss pairs\dss
$(\fff \tau\fff,\pff \sigma \fff)$\dss
such that $\tau$ is\dss an 
$(\fff n\dff -\dff 1 \fff)$\dnsp-simplex of\trs $T$\sss 
and\qss $\varphi\dff(\dff \tau\dff)\off =\off \delta_{\dff i}$\nsp,\oss 
and\dss $\sigma$\dss is\dss an $n$\dnsp-simplex\dss of\qss $T$\dss having $\tau$ as a face.\oss
Clearly\halfff,\oss the simplices $\tau$ occurring\dss in such\sss pairs are
exactly\dss the elements of\qss $G\qff \cup\qff H$\nnsp,\oss
i.e.\qss the vertices of\trs the first\dss kind of\trs the graph\dss $\mathbb{G}_{\dff i}$\nsp.\oss
The simplices $\sigma$ occurring\dss in such\sss pairs are
exactly\dss the elements of\qss $E\qff \cup\qff F$\nnsp,\oss
i.e.\qss the vertices of\trs the second\dss kind of\trs the graph\dss $\mathbb{G}_{\dff i}$\nsp.\oss
The pair\dss
$(\fff \tau\fff,\pff \sigma \fff)$\dss
is\dss among\dss the counted\dss pairs\dss if\trs and\dss only\trs if\dss
$\tau$ and $\sigma$ are connected\dss by\sss an edge in\dss $\mathbb{G}_{\dff i}$\nsp.\oss
We see\sss that\dss the graph\dss $\mathbb{G}_{\dff i}$\dss
is\dss present\dss in\dss the combinatorial\dss proof\dss
even more explicitly\dss than\dss in\dss the cochains-based one.\oss

\myuppar{Searching\dss for\sss elements of\dss $E$\nnsp.}
The graph-theoretical\dss version of\trs the proof\dss is\dss very\sss attractive,\oss
especially\dss because\sss it\sss suggests a way\sss of\dss finding\sss
elements of\trs $E$\nnsp,\oss i.e.\qss the $n$\dnsp-simplices $\sigma$
such\dss that\qss
$\varphi\dff(\dff \sigma\dff)\off =\off \delta$\nnsp.\oss
Indeed,\oss since both sets\sss $E$\sss and\dss $H$\sss have an odd\dss number
of\dss elements,\oss
there\dss is\dss at\dss least\sss one path starting\dss in\sss $H$\sss
and ending\dss in\trs $E$\nnsp.\oss
Following\dss a\sss path\dss in\sss $\mathbb{G}_{\dff i}$\sss can be easily\dss turned\dss into
an algorithm,\oss
but\sss at\dss the first\sss sight\dss such an algorithm\dss is\dss hardly\sss satisfactory\fff:\oss
it\sss seems\sss that\dss 
in\sss order\dss to find\dss in\dss this\sss way\sss  
even one element\sss of\trs $E$\sss
one needs\sss to know\sss all\sss elements of\trs $H$\nnsp.\oss
Still,\oss at\dss the very\dss least\sss one can\dss replace an exhaustive search among\dss
the $n$\dnsp-simplices of\trs $T$\dss by\sss 
an exhaustive search 
among\dss the $(\fff n\dff -\dff 1\fff)$\dnsp-simplices\dss of\trs $T_{\dff i}$\dss
plus\dss following\sss several\dss paths.\oss\vspace{1pt}

A moment\sss of\trs thought\dss leads\sss to\sss the conclusion\dss
that\sss one shouldn't\sss expect\sss that\dss simply\dss following\dss paths 
starting\dss in\sss $H$\sss would\dss be\sss a satisfactory\sss
search strategy\halfff.\oss
Indeed,\oss this method\sss 
does not\dss fully\sss reflect\sss even\dss the inductive
step:\oss 
the fact\dss that\sss elements of\trs $E$\sss not\dss
reachable in\dss this way\sss are pairwise connected\dss by\sss paths of\trs
$\mathbb{G}_{\dff i}$\sss 
is\dss equally\dss important\halfff.\oss
Even\dss more importantly\halfff,\oss
the graph\sss $\mathbb{G}_{\dff i}$\sss
encodes only\sss one step of\trs the induction,\oss
and one step\dss is\dss not\sss sufficient\sss
even\dss to establish\dss that $n$\dnsp-simplices $\sigma$
such\dss that\qss
$\varphi\dff(\dff \sigma\dff)\off =\off \delta$\qss
exist\halfff.\oss
The proof\dss of\trs the existence implicitly\dss involves
a similar graph\sss related\dss to\sss the 
$(\fff n\dff -\dff 1\fff)$\dnsp-face\dss $\delta_{\dff i}$ of\qss $\delta$\nnsp,\oss
a\sss graph\sss related\dss to an $(\fff n\dff -\dff 2\fff)$\dnsp-face
of\trs $\delta_{\dff i}$\nsp,\oss etc.\oss
One can combine\sss the corresponding\dss path-following algorithms\qss
(including\dss paths connecting\sss one element\sss of\trs $E$\sss with another\dss
for\dss proper\dss faces of\dss $\delta$ in\dss the role of\dss $\delta$\nsp),\oss
but\dss there\dss is\dss a\sss better\sss approach.\oss
Namely\halfff,\oss
one can concatenate all\dss relevant\sss graphs into a single graph.\pss\vspace{1pt}

\myuppar{A\dss graph\sss $\mathbb{G}$\sss encoding\dss all\dss steps of\qss induction.}
The\sss graph\sss $\mathbb{G}$\sss depends\sss
not\sss only\sss on\dss the choice of\dss $i$\nnsp,\oss
but\sss also on\dss the corresponding choices in\dss lower dimensions.\oss
So,\oss let\dss\vspace{3pt}
\[
\quad
\delta^{\fff n}\off \supset\off
\delta^{\fff n\dff -\dff 1}\off \supset\off
\ldots\off \supset\off
\delta^{\fff 1}\off \supset\off
\delta^{\fff 0}
\] 

\vspace{-9pt}
be a sequence of\trs faces of\dss $\delta$\dss starting\dss with\dss
$\delta^{\fff n}\off =\dff\off \delta$\dss and such\dss that\dss the dimension
of\trs $\delta^{\fff m}$\dss is\sss $m$\nnsp.\oss
Let\trs $T^{\dff m}$\trs be\sss the\sss triangulation of\trs $\delta^{\fff m}$\dss
consisting of\dss simplices of\trs $T$\sss contained\dss in\dss $\delta^{\fff m}$\nnsp.\oss
Now\dss we are ready\dss to define\dss $\mathbb{G}$\nnsp.\oss
For\sss every\qss $m\qff \leq\qff n$\qss
every\sss $m$\dnsp-simplex\sss $\sigma$ of\qss $T^{\dff m}$\sss
such\dss that\qss\vspace{3pt}
\[
\quad
\varphi\dff(\dff \sigma\dff)
\off =\off\qff 
\delta^{\fff m}
\hspace*{1.2em}\mbox{or}\hspace*{1.2em}
\delta^{\dff m\dff -\dff 1}
\] 

\vspace{-9pt}
is\dss a\sss vertex\sss of\trs $\mathbb{G}$\nnsp.\oss
Also,\oss if\qss $1\qff \leq\qff m\qff \leq\qff n$\nnsp,\oss
then every $(\fff m\dff -\dff 1\fff)$\dnsp-simplex\sss $\tau$
of\trs $T^{\dff m}$\sss such\dss that\vspace{3pt}
\[
\quad
\varphi\dff(\dff \tau\trf)
\off =\off\qff 
\delta^{\dff m\dff -\dff 1}
\] 

\vspace{-9pt}
is\dss a vertex\sss of\trs $\mathbb{G}$\nnsp.\oss
There are no other vertices.\oss
Two vertices\dss $\sigma\fff,\pff \tau$\dss as above are connected\dss
by\sss an edge if\trs $\tau$\sss is\dss an $(\fff m\dff -\dff 1\fff)$\dnsp-face of\dss $\sigma$\nnsp.\oss
There are no other edges.\pss
Note\sss that\qss if\dss an $(\fff m\dff -\dff 1\fff)$\dnsp-simplex $\tau$
as above\dss is\dss contained\dss in\sss $\delta^{\fff m\dff -\dff 1}$\dnsp,\oss
then $\tau$ is\dss connected\dss by\sss an edge\sss to some $(\fff m\dff -\dff 2\fff)$\dnsp-simplex
of\trs $T^{\dff m\dff -\dff 1}$\nnsp.\qff\oss
If\pss $\delta^{\fff n\dff -\dff 1}\off =\off\dff \delta_{\dff i}$\nsp,\oss
then,\oss obviously\halfff,\pss $\mathbb{G}_{\dff i}$\sss 
is\dss a subgraph of\trs $\mathbb{G}$\nnsp.\oss\vspace{1pt}

\myuppar{Theorem.}
\emph{Every\sss vertex\sss of\trs $\mathbb{G}$\sss
is\dss an endpoint\sss of\dss one or\dss two edges.\oss
A vertex $\sigma$ is\dss an endpoint\sss of\dss only\sss one edge\dss
if\trs and\dss only\trs if\trs either\qss 
$\sigma\off =\off \delta^{\dff 0}$\nnsp,\qff\oss
or $\sigma$ is\dss an $n$\dnsp-simplex\dss and\qss
$\varphi\dff(\dff \sigma\dff)
\off =\off 
\delta^{\fff n}\off =\off \delta$\nnsp.\oss}

\proof
Suppose\sss first\dss that\dss $\sigma$ is\dss 
an $m$\dnsp-simplex\sss of\trs $T^{\dff m}$\sss
such\dss that\qss
$\varphi\dff(\dff \sigma\dff)
\off =\off 
\delta^{\fff m}$\dnsp.\oss
If\qss $m\off =\off 0$\nnsp,\oss then\qss $\sigma\off =\off \delta^{\dff 0}$\nnsp.\oss
Otherwise\dss there\dss is\dss exactly\sss one $(\fff m\dff -\dff 1\fff)$\dnsp-face $\tau$
of\dss $\sigma$ such\dss that\qss
$\varphi\dff(\dff \tau\trf)
\off =\off 
\delta^{\dff m\dff -\dff 1}$\dnsp.\oss
If\qss $m\qff \leq\qff n\qff -\qff 1$\nnsp,\oss
then\dss there\dss is\dss exactly\sss one $(\fff m\dff +\dff 1\fff)$\dnsp-simplex\sss $\rho$\sss
of\trs $T^{\dff m\dff +\dff 1}$\sss such\dss that\sss $\sigma$ is\dss a\sss face of\dss
$\rho$\nnsp.\oss 
Clearly\halfff,\pss
$\varphi\dff(\dff \rho\dff)
\off \supset\off
\varphi\dff(\dff \sigma\dff)
\off =\off
\delta^{\fff m}$\dnsp,\oss
and since
$\varphi$ is\dss pseudo-identical\dss map,\pss
$\varphi\dff(\dff \rho\dff)
\off \subset\off
\delta^{\fff m\dff +\dff 1}$\dnsp.\oss
It\dss follows\dss that\sss $\rho$ is\dss a\sss vertex\sss of\trs $\mathbb{G}$\nnsp.\oss
Clearly\halfff,\pss $\sigma$ is\dss connected\dss by\sss an edge only\dss with $\tau$
if\qss $m\off =\off n$\nnsp,\oss
only\dss with $\rho$ if\qss $m\off =\off 0$\nnsp,\oss 
and only\dss with $\tau$ and $\rho$ if\qss
$1\qff \leq\qff m\qff \leq\qff n\qff -\qff 1$\nnsp.\oss

Suppose\sss now\dss that\dss $\sigma$ is\dss 
an $m$\dnsp-simplex\sss of\trs $T^{\dff m}$\sss
such\dss that\qss
$\varphi\dff(\dff \sigma\dff)
\off =\off 
\delta^{\fff m\dff -\dff 1}$\dnsp.\oss
By\dss the definition,\oss in\dss this case $\sigma$ is\dss not\sss connected\dss by\sss
an edge with any\sss simplex of\trs $T^{\dff m\dff +\dff 1}$\sss not\trs belonging\dss
to\dss $T^{\dff m}$\dnsp.\oss
The arguments used\dss for $\mathbb{G}_{\dff i}$ with\dss $T^{\dff m}$\sss
in\dss the role of\trs $T$\sss show\dss that\dss in\dss this case $\sigma$
is\dss connected\dss with exactly\dss two vertices,\oss
both of\dss which are $(\fff m\dff -\dff 1\fff)$\dnsp-simplices of\trs $T^{\dff m}$\dnsp.\oss

Finally\halfff,\oss
let\dss us\dss consider\sss an
$(\fff m\dff -\dff 1\fff)$\dnsp-simplex\sss $\tau$
of\trs $T^{\dff m}$\sss such\dss that\qss
$\varphi\dff(\dff \tau\trf)
\off =\off\qff 
\delta^{\dff m\dff -\dff 1}$\dnsp.\oss
If\dss $\tau$ is\dss actually\sss a\sss simplex of\trs $T^{\dff m\dff -\dff 1}$\dnsp,\oss
then $\tau$ is\dss connected\dss by\sss an edge\sss with\sss two vertices\sss
by\dss the first\dss paragraph of\trs the proof\dss applied\dss to $\tau$ and\qss
$m\qff -\qff 1$\qss in\dss the roles of\dss $\sigma$ and $m$ respectively\halfff.\oss
Otherwise\sss the arguments used\dss for $\mathbb{G}_{\dff i}$ with\dss $T^{\dff m}$\sss
in\dss the role of\trs $T$\sss show\dss that\dss in\dss this case $\tau$
is\dss connected\dss with exactly\dss two vertices,\oss
both of\dss which are $m$\dnsp-simplices of\trs $T^{\dff m}$\dnsp.\oss  \eproof

\myuppar{Corollary\halfff.}
\emph{The\dss graph\trs 
$\mathbb{G}$\sss consists\dss of\qss several\dss
disjoint\dss paths\sss and\dss cycles.\oss
With\dss exception of\pss $\delta^{\dff 0}$\nsp\dnsp,\oss the\sss endpoints of\qss these paths are
$n$\dnsp-simplices $\sigma$ such\dss that\qss
$\varphi\dff(\dff \sigma\dff)
\off =\off\dff 
\delta$\nnsp.\oss}  \eproof

\myuppar{Corollary\halfff.}
\emph{The number\sss of\trs $n$\dnsp-simplices $\sigma$ such\dss that\qss
$\varphi\dff(\dff \sigma\dff)
\off =\off\dff 
\delta$\qss is\trs odd.\oss} \eproof

\myuppar{Path-following\sss algorithms.}
Following\dss the
unique path of\trs $\mathbb{G}$\sss
starting\sss at\sss $\delta^{\dff 0}$\sss
leads\sss to an $n$\dnsp-simplex $\sigma$ sich\dss that\qss
$\varphi\dff(\dff \sigma\dff)\off =\off \delta$\nnsp.\qff\oss
Following\dss this path can be easily\dss turned\dss into an algorithm,\oss
which\dss turns out\dss to be equivalent\dss to\dss
one of\pss Scarf's\qss algorithms\qss \cite{sc3}.\oss
Cf.\qss \cite{sc3},\oss Lemma\qss 3.4.\oss
In\dss the context\sss of\oss Brouwer's\qss fixed-point\dss theorem\dss
such simplices $\sigma$ may\dss be interpreted as approximate
fixed\dss points of\dss continuous maps\qss
$\delta\qff \ttoo\qff \delta$\nnsp,\oss
and\dss the path-following\sss algorithms were actually\dss used\dss for
computing\sss approximations\sss to fixed\dss points.\oss
Cf.\qss \cite{sc2},\oss \cite{sc3}.\oss

\myuppar{Historical\dss remarks.}
A\dss proof\dss of\pss Sperner's\trs lemma\dss
based on\dss path-following arguments was published\dss in\qss 1967\qss
by\pss D.I.A.\qss Cohen\qss \cite{c}.\oss
His\sss proof\dss amounts\sss to using\dss the standard\dss induction\dss
by\sss $n$ and\dss the graph $\mathbb{G}_{\dff i}$ for\dss
the step of\trs induction.\oss
Cohen\dss did\dss not\dss relate\sss his proof\trs
to any\sss of\trs the classical\dss proofs.\oss
In\qss 1979\pss A.{\halfff}W.\qss Tucker\qss \cite{t}\qss wrote about\qss
Sperner's\trs lemma and\dss Cohen's\dss proof\dff:\vspace{-6pt}

\begin{quoting}
This\sss lemma,\qss proved\trs by\sss a\sss simple existential\sss argument\dss through\dss
induction on $n$\nnsp,\oss \ldots\oss
Now,\oss however\halfff,\oss we have an algorithmic proof\dss of\pss Sperner's\qss lemma,\pss
thanks\sss to an\dss idea of\qss Cohen\qss \cite{c}.\oss 
\end{quoting}

\vspace{-6pt}
Also\dss in\qss 1967\pss H.\qss Scarf\qss published\qss \cite{sc2}\qss a\sss proof\dss of\qss
Brouwer's\qss fixed-point\dss theorem\sss based on a combinatorial\trs theorem\dss
resembling\qss Sperner's\trs lemma.\oss
Scarf\qss proved\dss this combinatorial\qss theorem\dss by\sss
a\sss path-following\dss algorithm\dss realizing\dss the whole inductive argument\halfff.\oss
Later on\qss Scarf\qss proved\qss
Sperner's\trs lemma\sss in\sss a\sss similar\dss manner\halfff.\oss
See\qss \cite{sc3}.\oss
In contrast\dss with\qss D.I.A.\qss Cohen,\oss Scarf\qss was quite forthcoming\sss
and explained\dss his sources of\dss inspiration.\oss
His\sss proof\dss of\pss Brouwer's\qss fixed-point\trs theorem\dss
was a\sss byproduct\sss of\qss his\sss fundamental\dss work\dss in game\sss theory\qss
\cite{sc1}\qss and strongly\dss influenced\dss by\trs linear\dss programming
and a\sss paper\dss by\qss C.E.\qss Lemke\qss \cite{le}.\oss
Scarf\dss wrote\dss in\qss \cite{sc2}\dff:\vspace{-6pt}

\begin{quoting}
\ldots\qff,\oss
Sperner's\trs lemma\sss suggests no procedure for\dss the determination
of\dss an approximate\sss fixed\dss point\sss other\dss than an exhaustive search
of\dss all\sss subsimplices until one\dss is\dss found\dss with\sss all\sss vertices\sss
labeled differently\halfff.\oss 
\ldots\qff,\oss
the algorithm\dss is\dss intimately\dss related\dss to\sss the procedure described\dss by\qss
Lemke\qss \cite{le}\qss for\dss the determination of\pss Nash\qss equilibrium\dss points
of\trs two-person nonzero-sum games.\oss 
\end{quoting}

\vspace{-6pt}
The above path-following\dss proof\trs shows\sss that\dss the classical\dss proofs
of\qss Sperner's\trs lemma\dss naturally\dss lead\dss to a proof\dss sharing\dss the main\dss
features of\pss Scarf's\trs proof\dss
and\sss to\sss the same\sss path-following algorithm.\oss
The main difference\dss is\dss in\dss 
the ways used\dss to piece\sss together all\sss
steps of\dss induction.\oss
Scarf\qss used\dss the so-called\qss
\emph{slack\dss vectors},\oss an\dss idea coming\dss from\dss linear\dss programming\halfff.\oss
This\sss idea\sss works in our context\sss also,\oss
but\trs is\dss harder\dss to motivate from\sss a\sss topological\dss point\sss of\dss view.\oss

\myuppar{The graphs\dss $\mathbb{G}_{\dff i}$\sss and\qss Alexander's\qss lemma.}
Admittedly\halfff,\oss one has\sss to be\sss more inventive in order\dss to see\sss
these graphs in\dss the\sss proof\dss of\pss
Alexander's\qss lemma\qss (see\qss Section\qss \ref{alexander-lemma}).\oss
Of\dss course,\pss
the proof\dss
should\dss be specialized\dss to\sss the case\qss $S'\off =\off T$\dnsp,\qss
$S\off =\off \Delta$\nnsp.\oss
But\dss in\dss this case\sss
no $(\fff n\dff -\dff 1\fff)$\dnsp-face of\dss $\delta$\sss
plays any\sss special\dss role.\oss
In\dss fact\halfff,\oss this proof\dss more naturally\dss leads\sss to\sss the union\dss
$\mathbb{U}$\sss of\trs all\dss graphs\dss $\mathbb{G}_{\dff i}$\nsp.

The proof\dss of\qss Alexander's\trs lemma\sss
is\dss based on a study\sss of\trs the action of\trs the map $\varphi$ on simplices
of\qss $T$\sss of\qss dimensions $n$ and\dss $n\qff -\qff 1$\nnsp.\oss
Clearly\halfff,\oss only\dss the simplices $\sigma$ such\dss that\dss the
dimension of\dss $\varphi\dff(\dff \sigma\dff)$\dss is\dss equal\trs to $n$
or\dss $n\qff -\qff 1$\dss matter\halfff.\oss
One may\dss think\dss that\dss $n$\dnsp-simplices $\sigma$ such\dss that\qss
$\varphi\dff(\dff \sigma\dff)$\dss is\dss an  $(\fff n\dff -\dff 1\fff)$\dnsp-simplex
are\sss irrelevant\dss because for\dss them\qss
$\varphi_{\dff *}\dff(\dff \sigma\dff)\off =\off 0$\nnsp.\oss
In\dss fact\halfff,\oss they\sss are highly\dss relevant\trs because
a crucial\sss step in\dss the proof\trs is\dss the application of\qss Theorem\qss 1.\oss
The only\dss nontrivial\dss part\sss of\trs the proof\dss of\qss Theorem\qss 1\qss
is\dss exactly\dss the part\sss dealing\dss with such simplices.\oss 

So,\oss the proof\dss of\qss Alexander's\trs lemma\sss
suggests\sss to consider\dss the graph\dss having as vertices all\sss
simplices $\sigma$ of\trs $T$
such\dss that\dss the image\qss
$\varphi\dff(\dff \sigma\dff)$\dss is\dss equal\dss either\dss to $\delta$ or\dss
to one of\dss the $(\fff n\dff -\dff 1\fff)$\dnsp-faces of\dss $\delta$\nnsp.\oss
In order\dss to encode\sss the relation\qss 
``\dnsp\hnsp$\tau$ is\dss a\sss face of\dss $\sigma$\nnsp''\qss
we connect\dss two vertices by\sss an edge when
one of\dss them\dss is\dss a\sss proper\dss face of\trs the other\halfff.\oss
The resulting\dss graph\dss is\dss the union\dss
$\mathbb{U}$\sss of\trs graphs\dss $\mathbb{G}_{\dff i}$\nsp.\oss

The above\sss graph-theoretical\sss arguments do not\sss apply\dss to\sss the graph $\mathbb{U}$\sss
by\dss itself\halfff.\oss 
While some vertices of\trs $\mathbb{U}$ are connected\dss by\sss an edge with\qss
$n\qff +\qff 1$\qss vertices,\oss they\sss are not\dss the source of\dss
difficulties.\oss
Indeed,\oss they\sss are exactly\dss the vertices $\sigma$ such\dss that\qss
$\varphi\dff(\dff \sigma\dff)\off =\off \delta$\nnsp.\oss
But\dss the parity\sss argument\dss would\dss be destroyed\qss
if\qss there are paths in\dss $\mathbb{U}$\dss connecting\dss
two $(\fff n\dff -\dff 1\fff)$\dnsp-simplices\sss in\dss
two different\dss faces of\dss $\delta$\sss without\dss passing\dss through a vertex $\sigma$
such\dss that\qss
$\varphi\dff(\dff \sigma\dff)\off =\off \delta$\nnsp.\oss
The proof\dss of\qss Theorem\qss 1\qss shows\sss that\dss this never\sss
happens.\oss
Indeed,\oss Case\qss 3\qss of\trs this proof\dss shows\sss that\trs if\trs
$\tau\fff,\pff \sigma$\dss are\sss two vertices of\dss $\mathbb{U}$ and $\tau$ is\dss a\sss face
of\dss $\sigma$\nnsp,\oss
then either\qss
$\varphi\dff(\dff \sigma\dff)\off =\off \delta$\nnsp,\oss 
or\dss the images $\varphi\dff(\dff \sigma\dff)$ and $\varphi\dff(\dff \tau\trf)$
are equal\pss ({\fff}to\sss the same face\dss $\delta_{\dff i}$\nsp,\qss $i\qff \in\pff I$\nsp).\oss
Therefore a path\dss in\dss $\mathbb{U}$\dss starting\sss
at\sss an $(\fff n\dff -\dff 1\fff)$\dnsp-simplex\sss contained\dss in $\delta_{\dff i}$
actually\sss stays\dss in\dss $\mathbb{G}_{\dff i}$\dss until\trs it\dss reaches $\sigma$
with\qss $\varphi\dff(\dff \sigma\dff)\off =\off \delta$\nnsp.\oss

\mysection{Cohomology\qss groups}{cohomology}

\myuppar{Discarding\sss coboundaries.}
The pairing\dss with\dss
$\fclass{\delta\fff}$\dss
allows\dss to discard\dss the coboundaries\dss
$\partial^*\dff (\dff \tau \dff)$\dss
such\dss that\dss $\tau$\dss is a simplex of\dss $T$\nnsp,\oss
but\dss not\sss of\dss $\bd\dff T$\nnsp,\oss 
without any further analysis.\oss
It is only natural\dss to discard\dss these coboundaries in a systematic way
by taking the quotient space of\dss the vector space of\dss $n$\dnsp-cochains 
by a suitable subspace.\oss
In\dss the\qss \emph{cohomology\dss theory}\qss all\sss coboundaries are discarded\dss in\dss this way.\oss
But\dss we need\dss 
to keep\sss the coboundaries\dss
$\partial^*\dff (\dff \tau \dff)$\dss 
such\dss that\dss $\tau$\dss is a simplex of\dss $\bd\dff T$\nnsp.\oss
A natural\sss way\dss to do this is provided\dss by\dss the\qss
\emph{relative cohomology\dss theory},\oss
a version\dss based on\qss
\emph{relative\dss cochains}.\oss
The following\dss theorem\dss is\dss the starting\dss point\halfff.\oss

\myuppar{\nnsp$\partial^*\hnsp\partial^*$\dnsp\dnsp-theorem.}\oss
$\partial^*\dff \circ\pff \partial^*
\off =\off
0$\dnsp.\oss

\proof
This\sss immediately\dss follows from\dss the\dss $\partial\partial$\dnsp-theorem\dss by dualizing\halfff.\oss
Here is a direct\dss proof\halfff.\oss
It\dss is sufficient\dss to show\dss that\qss
$\partial^*\dff \circ\pff \partial^*\dff(\dff \sigma\dff)
\off =\off
0$\qss
for every simplex $\sigma$\dnsp.\oss
If\dss $\sigma$\dss is an $m$\dnsp-simplex,\oss then\vspace{1.5pt}
\[
\quad
\partial^*\dff \circ\pff \partial^*\dff(\dff \sigma\dff)
\off =\off\qff
\sum\qff \partial^*\dff(\dff \tau\dff)
\qff,
\]

\vspace{-10.5pt}
where $\tau$ runs over\dss $(\fff m\dff +\dff 1 \fff)$\dnsp-simplices
having $\sigma$ as a face.\oss
It\dss follows\dss that\qss
$\partial^*\dff \circ\pff \partial^*\dff(\dff \sigma\dff)$\qss
is equal\dss to a sum of\dss $(\fff m\dff +\dff 2 \fff)$\dnsp-simplices
having $\sigma$ as a face.\oss
If\dss $\rho$\dss is such a simplex,\oss
then $\rho$ has exactly\dss two $(\fff m\dff +\dff 1 \fff)$\dnsp-faces
containing $\sigma$\nnsp.\oss
If\qss $\tau'\fff,\pff \tau''$\qss are\sss these\sss two faces,\oss
then $\rho$ is a summand of\qss
$\partial^*\dff(\dff \tau'\dff)$\dss and\qss
$\partial^*\dff(\dff \tau''\dff)$\dss
and of\dss no other\dss
$\partial^*\dff(\dff \tau\dff)$\nnsp.\oss
It\dss follows\dss that\sss every\dss $(\fff m\dff +\dff 2 \fff)$\dnsp-simplex $\rho$ 
having $\sigma$ as a face occurs in our sum exactly\dss two times
and no other\dss $(\fff m\dff +\dff 2 \fff)$\dnsp-simplices do.\oss
Since\qss $2\off =\off 0$\qss in\dss $\ftwo$\nsp,\oss
the whole sum vanishes and\dss hence\qss
$\partial^*\dff \circ\pff \partial^*\dff(\dff \sigma\dff)
\off =\off
0$\dnsp.\oss  \eproof

\myuppar{Cocycles,\oss coboundaries,\oss and cohomology.}
Let\dss $S$\dss be a simplicial complex and\dss $m$\dss be a non-negative integer\halfff.\oss
An $m$\dnsp-cochain $\alpha$ of\dss $S$\dss is called a\qss \emph{cocycle}\pss
if\oss
$\partial^*\dff(\dff \alpha\dff)
\off =\off
0$\pss
or $m$ is equal\dss to\sss the dimension of\dss $S$\dnsp,\oss
and\dss a\qss \emph{coboundary}\pss
if\pss
$\alpha
\off =\off
\partial^*\dff(\trf \beta\dff)$\qss
for\dss some\dss $(\fff m\dff -\dff 1 \fff)$\dnsp-cochain\dss $\beta$\nnsp.\qff\oss
Let\vspace{1.5pt}
\[
\quad
Z^{\fff m}\fff(\dff S \dff)
\hspace*{1.0em}\mbox{and}\hspace*{1.0em}
B^{\fff m}\fff(\dff S \dff)
\]

\vspace{-10.5pt}
be\sss the spaces of\sss $m$\dnsp-cochains which are,\oss
respectively\halfff,\oss cocycles and\dss coboundaries.\qss 
By\qss the\qss $\partial^*\hnsp\partial^*$\dnsp\dnsp-the\-o\-rem\dss 
every \nsp$m$\dnsp-coboundary\sss is\sss an \nsp$m$\dnsp-cocycle,\qss
i.e.\qss
$B^{\fff m}\fff(\dff S \dff)
\dff \subset\qff
Z^{\fff m}\fff(\dff S \dff)$\nnsp.\qss
The quotient\sss space\vspace{1.5pt}
\[
\quad
H^{\fff m}\fff(\dff S \dff)
\off =\off\qff
Z^{\fff m}\fff(\dff S \dff)\dff
\bigl/\bigr.\qff
B^{\fff m}\fff(\dff S \dff)
\]

\vspace{-10.5pt}
is\dss called\dss the
\emph{$m$\dnsp-dimensional cohomology\dss group of}\qss
$S$\dnsp.\oss
The\sss term\qss \emph{``cohomology\dss group''}\pss
is standard even when\dss it\dss is a vector space.\oss
The image of\dss a cocycle\qss
$\alpha\qff \in\qff Z^{\fff m}\fff(\dff S \dff)$\qss
in\dss the cohomology\dss group\qss
$H^{\fff m}\fff(\dff S \dff)$\qss
is called\dss the\qss \emph{cohomology\dss class}\qss of\dss $\alpha$\dss
and\dss is denoted\dss by\dss $\hclass{\alpha}$\nnsp.\oss

\myuppar{Relative cohomology.}
Suppose now\dss that\dss $Q$\dss is a\qss 
\emph{subcomplex}\qss of\dss $S$\nnsp,\oss
i.e.\qss that every vertex of\trs $Q$\dss is a vertex of\trs $S$\dss
and every simplex of\dss $Q$\dss is a simplex of\dss $S$\nnsp.\oss
A\qss \emph{relative\dss $m$\dnsp-cochain of\qss the pair}\qss
$(\dff S\fff,\pff Q \dff)$\dss
is an\dss $m$\dnsp-cochain of\trs $S$\dss
vanishing on every\dss $m$\dnsp-simplex\sss of\trs $Q$\dss
if\trs considered as a linear\dss functional\qss
$C_{\fff m}\fff(\dff S \dff) \ttoo \ftwo$\nnsp.\oss
In\dss the language of\dss formal sums of\dss simplices\sss
this means\sss that\dss no simplices of\trs $Q$\dss are allowed\dss
to enter\dss the sum.\oss

Obviously,\oss the coboundary of\dss a relative\dss
$(\fff m\dff -\dff 1 \fff)$\dnsp-cochain
is a relative $m$\dnsp-cochain.\oss
A relative $m$\dnsp-cochain $\alpha$ is called a\qss \emph{relative cocycle}\qss
if\qss it\qss is\dss a\sss cocycle as a chain,\oss
and a\qss \emph{relative coboundary}\pss
if\pss
$\alpha
\off =\off
\partial^*\dff(\trf \beta\dff)$\qss
for\dss some relative\dss $(\fff m\dff -\dff 1 \fff)$\dnsp-cochain\dss $\beta$\nnsp.\qff\oss
Let\vspace{3pt}
\[
\quad
Z^{\fff m}\fff(\dff S\fff,\pff Q \dff)
\hspace*{1.0em}\mbox{and}\hspace*{1.0em}
B^{\fff m}\fff(\dff S\fff,\pff Q \dff)
\]

\vspace{-9pt}
be\sss the spaces of\dss relative $m$\dnsp-cochains which are 
relative cocycles and\dss relative coboundaries respectively.\oss
Like before,\oss the\dss $\partial^*\hnsp\partial^*$\dnsp\dnsp-the\-o\-rem\dss implies\sss
that\sss every relative coboundary\dss is a relative cocycle,\oss
i.e.\qss
$B^{\fff m}\fff(\dff S\fff,\pff Q \dff)
\qff \subset\off
Z^{\fff m}\fff(\dff S\fff,\pff Q \dff)$\nnsp.\oss
The quotient\sss space\vspace{3pt}
\[
\quad
H^{\fff m}\fff(\dff S\fff,\pff Q \dff)
\off =\off\qff
Z^{\fff n}\fff(\dff S\fff,\pff Q \dff)\dff
\bigl/\bigr.\qff
B^{\fff m}\fff(\dff S\fff,\pff Q \dff)
\]

\vspace{-9pt}
is\dss called\dss the
\emph{$m$\dnsp-dimensional cohomology\dss group of}\pss the pair
$(\dff S\fff,\pff Q \dff)$\dnsp.\oss
The image of\dss a cocycle\qss
$\alpha\qff \in\qff Z^{\fff m}\fff(\dff S\fff,\pff Q \dff)$\qss
in\qss
$H^{\fff m}\fff(\dff S\fff,\pff Q \dff)$\qss
is called\dss the\qss \emph{cohomology\dss class}\qss of\dss $\alpha$\dss
and\dss is denoted\dss by\dss $\hclass{\alpha}$\nnsp.\oss

\myuppar{Pseudo-manifolds.}
We are especially\dss interested\dss
in\dss the\qss \emph{top-dimensional}\qss cohomology\dss groups,\oss
i.e.\qss in\dss the $n$\dnsp-dimensional\sss cohomology\dss groups
of\dss complexes\dss $S$\dss and\dss pairs $(\dff S\fff,\pff Q \dff)$
with $n$ equal\dss to\sss the dimension of\dss $S$\dnsp.\oss
There is a class of\dss complexes for which\dss the\sss
top-dimensional cohomology\dss groups
are easy\dss to determine.\oss
See\dss Theorem\qss 6\qss below.\oss
Let\dss $S$\dss be a simplicial complex of\dss dimension\dss $n$\nnsp.

Guided\dss by\dss the non-branching\sss property of\dss triangulations of\dss a simplex\halfff,\oss
we will say that\dss $S$\dss 
is\qss \emph{non-branching}\oss
if\dss each $(\fff n\dff -\dff 1 \fff)$\dnsp-simplex of\dss $S$\dss
is a face of\dss either one or two $n$\dnsp-simplices.\oss
If\dss $S$\dss is non-branching\halfff,\oss
then the\qss \emph{boundary}\qss $\partial\dff S$\dss  
is defined as simplicial complex having as the simplices\sss all\trs faces of\dss
$(\fff n\dff -\dff 1 \fff)$\dnsp-simplices of\trs $S$\dss which are faces of\dss
exactly one $n$\dnsp-simplex\halfff.\oss

We will say\dss that\dss $S$\dss is\qss 
\emph{strongly\dss connected}\oss
if\qss 
for every\dss two $n$\dnsp-simplices\qss 
$\sigma,\pff \sigma\halfff'$\qss
of\trs $S$\dss 
there is a sequence\oss 
$\sigma
\off =\off 
\sigma_{\fff 0}\fff,\off\off 
\sigma_{\fff 1}\fff,\off\off 
\ldots\fff,\off\off
\sigma_{\fff k\dff -\dff 1}\fff,\off\off
\sigma_{\fff k}
\off =\off 
\sigma\halfff'$\oss 
of\dss $n$\dnsp-simplices
such that\dss 
$\sigma_{\fff l}$\dss and\dss $\sigma_{\fff l\dff +\dff 1}$\dss 
have a common\dss
$(\fff n\dff -\dff 1 \fff)$\dnsp-face\dss
for every\qss $l\qff \leq\qff k\qff -\qff 1$\nnsp.\oss

Usually\dss these\sss two conditions are imposed\dss together with\dss another one.\oss
The complex\dss $S$\dss is said\dss to be\qss
\emph{dimensionally\dss homogenous}\qss
if\dss every\sss its simplex is a face of\dss an $n$\dnsp-dimensional\sss simplex.\oss
This condition is natural,\oss but\dss is hardly\dss relevant\dss for us.\oss
The complex\dss $S$\dss is called\qss
a\qss \emph{pseudo-manifold}\oss 
if\qss it\trs is non-branching\halfff,\oss
strongly connected,\oss
and\dss dimensionally\dss homogeneous.\oss

\myuppar{Triangulations of\dss a simplex.}
Examples of\dss pseudo-manifolds are provided\dss by\dss triangulations\sss of\dss
simplices.\oss
Every\dss triangulation\dss $T$\dss of\dss $\delta$ is a pseudo-manifold of\dss
dimension\dss $n$\nnsp.\oss
We already\dss
implicitly accepted\dss the non-branching\dss property as
geometrically obvious.\oss
The fact\dss that\dss $T$\dss is dimensionally\dss homogenous is also geometrically obvious.\oss

The fact\dss that\dss $T$\dss is strongly connected\dss
is a little\sss less obvious.\oss
Let\qss 
$\sigma\fff,\pff \sigma\halfff'$\qss
be\sss two $n$\dnsp-simplices of\trs $T$\nnsp.\oss
They can\dss be connected\dss by a path\sss in\dss $\delta$\nnsp.\oss
If\trs this path intersects only $n$\dnsp-simplices and
$(\fff n\dff -\dff 1 \fff)$\dnsp-simplices of\dss $T$\nnsp,\oss
then one can read off\trs the required sequence by\dss following\dss this path.\oss
If\dss this path intersects an $m$\dnsp-simplex with\qss
$m\qff \leq\qff n\qff -\qff 2$\nnsp,\oss
one can replace a segment\sss of\dss this path crossing\dss this $m$\dnsp-simplex
by a segment\dss bypassing\sss it\halfff.\oss
We leave the details\sss to the reader\halfff.\oss

By\dss the non-branching\dss property\sss of\dss triangulations of\trs $\delta$\dss
the boundary\dss $\partial\dff T$\dss consists of\dss all simplices of\trs $T$\dss
contained\sss in\dss the boundary\dss $\bd\fff \delta$\nnsp.\oss
In other words,\oss the boundary\dss $\partial\dff T$\dss
is nothing else but\dss the complex which was denoted\dss by\dss $\bd\dff T$\dss above.\oss
Recall\dss that\dss the complex $\Delta$ is\dss the\sss tautological\dss
triangulation of\dss $\delta$\dss consisting of\dss $\delta$\dss itself\dss and\dss its faces.\oss
The discussion at\dss the beginning of\qss Section\qss \ref{alexander-lemma}\qss implies\sss that\qss
$\partial\dff T\off =\off \bd\dff T$\qss 
is a subdivision of\dss $\partial\dff \Delta$\dss
and\dss hence\sss the following\dss theorem\qss
(applied\dss to\dss
$S\off =\off \partial\dff \Delta$\dss 
and\dss $S'\off =\off \partial\dff T$\nsp)\qss
implies\sss that\dss $\partial\dff T$\dss is also a pseudo-manifold.\oss\vspace{4.25pt}

\myuppar{Theorem\qss 5.}
\emph{If\qss a complex\dss $S$\dss is a pseudo-manifold\dss
of\qss dimension\dss $n$\dss and\dss $S'$\dss is a subdivision of\dss $S$\dnsp,\oss
then\dss $S'$\dss is also a pseudo-manifold\dss
of\qss dimension\dss $n$\nnsp.\oss}\vspace{4.25pt}

\proof
Recall\dss that\dss for every\sss simplex $\sigma$ of\dss $S$\dss
the simplices of\dss $S'$\dss contained\dss in $\sigma$
form a\sss triangulation\dss $S'\fff(\dff \sigma\dff)$\dss
of\dss $\sigma$\dnsp.\oss
As we\dss just\dss saw,\oss
$S'\fff(\dff \sigma\dff)$\dss is a pseudo-manifold
and its boundary\dss $\partial\trf S'\fff(\dff \sigma\dff)$\dss
consists of\dss simplices of\dss $S'\fff(\dff \sigma\dff)$\dss
contained\dss in\dss $\bd\dff \sigma$\nnsp.\oss
Let $\tau\fff'$ be an $(\fff n\dff -\dff 1 \fff)$\dnsp-simplex of\dss $S'$\dnsp.\oss
Since\dss $S'$\dss is a subdivision of\dss $S$\dnsp,\oss
the simplex $\tau\fff'$ is contained\dss is some simplex of\dss $S$\nnsp.\oss

If\dss $\tau\fff'$ is contained\dss in an
$(\fff n\dff -\dff 1 \fff)$\dnsp-simplex $\tau$ of\dss $S$\dnsp,\oss
then\dss the simplices $\sigma\fff'$ of\dss $S'$\dss having $\tau\fff'$ as a face 
are in $1$\dnsp-to-\dnsp$1$ correspondence with simplices $\sigma$ of\dss $S$\dss
having $\tau$ as a face,\oss
and\dss hence\sss there are $1$ or $2$ of\dss such simplices $\sigma\fff'$\dnsp.\oss
If\sss $\tau\fff'$ is not\sss contained\dss in any
$(\fff n\dff -\dff 1 \fff)$\dnsp-simplex of\dss $S$\dnsp,\oss
then $\tau$ is contained\dss in a unique $n$\dnsp-simplex $\sigma$ of\dss $S$\dss
and\sss does not\dss belong\dss to\dss
$\partial\trf S'\fff(\dff \sigma\dff)$\nnsp.\oss
In\dss this case
$\tau\fff'$ is a face of\dss exactly\dss two $n$\dnsp-simplices 
of\dss $S'\fff(\dff \sigma\dff)$\dss
and\dss hence of\dss exactly\dss two
$n$\dnsp-simplices of\dss $S'$\dnsp.\oss

This proves\sss that\dss $S'$\dss is non-branching\halfff.\oss
The fact\dss that\dss $S'$\dss is strongly connected\dss follows
from\dss the strong connectedness of\dss $S$\dss and of\dss complexes\dss
$S'\fff(\dff \sigma\dff)$\nnsp.\oss 
Similarly,\oss $S'$\dss is dimensionally\dss homogeneous\dss because\dss
$S$\dss and\dss $S'\fff(\dff \sigma\dff)$\dss are.\oss
We\dss leave\sss the details\sss to\sss
the interested\dss readers.\oss  \eproof\vspace{4.25pt}

\myuppar{Theorem\qss 6.}
\emph{Let\qss 
$S$\dss be a strongly connected\dss non-branching complex of\dss dimension $n$\nnsp.\oss
Then\dss
$H^{\dff n}\dff(\dff S\fff,\pff \partial\dff S \dff)$\dss
is\dss a\sss vector\sss space\dss of\qss dimension\dss $1$\nnsp.\oss
Every\dss $n$\dnsp-simplex\sss $\sigma$\sss of\pss $S$\dss is a cocycle and\dss
its cohomology\sss class\qss $\hclass{\sigma}$\dss is\dss non-zero\sss and\dss hence\sss
is\dss a\sss basis\sss of\oss
$H^{\dff n}\dff(\dff S\fff,\pff \partial\dff S \dff)$\nnsp.\oss}\vspace{4.25pt}

\proof
Since $n$ is\sss the dimension of\dss $S$\dnsp,\oss
every $n$\dnsp-cochain is a cocycle and\dss hence\vspace{4.5pt}
\[
\quad
H^{\fff n}\fff(\dff S\fff,\pff \partial\dff S \dff)
\off =\off\qff
C^{\fff n}\fff(\dff S\fff,\pff \partial\dff S \dff)\dff
\bigl/\bigr.\qff
B^{\fff n}\fff(\dff S\fff,\pff \partial\dff S \dff)
\qff.
\]

\vspace{-7.5pt}
It\dss follows\dss that\qss
$H^{\dff n}\dff(\dff S\fff,\qff \partial\dff S \dff)$\qss
is generated\dss by\dss the cohomology classes
$\hclass{\sigma}$ of\dss $n$\dnsp-sim\-plices $\sigma$ of\trs $S$\dnsp.\oss
By\dss the definition of\dss $\partial\dff S$\dnsp,\oss
if\qss an\dss $(\fff n\dff -\dff 1\fff)$\dnsp-simplex\sss $\tau$\sss is
not\sss a simplex of\dss $\partial\dff S$\nnsp,\oss
then\sss 
$\tau$\sss is a face of\dss exactly\dss two $n$\dnsp-simplices\dss
$\sigma,\pff \sigma\halfff'$\dss of\qss $S$\nnsp.\oss
In\dss this case \vspace{0pt}
\begin{equation*}
\quad
\partial^*\dff (\dff \tau \dff)
\off =\off
\sigma\qff +\qff \sigma\halfff'
\off =\off
\sigma\qff -\qff \sigma\halfff'
\end{equation*}

\vspace{-12pt}
and\dss hence\dss
$\hclass{\sigma}
\off =\off
\hclass{\sigma\halfff'}$\nnsp.\oss
Since $S$ is strongly connected\halfff,\pss
this implies that all\sss cohomology classes
$\hclass{\sigma}$
of\dss $n$\dnsp-simplices $\sigma$
are equal\dss and\dss hence\sss the dimension of\dss
$H^{\dff n}\dff(\dff S\fff,\pff \partial\dff S \dff)$\dss
is\oss $\leq\qff 1$\nnsp.\oss

It\dss remains\sss to prove\sss that\dss the cohomology classes
$\hclass{\sigma}$
are non-zero.\oss
Let\dss us consider cochains as formal sums of\dss simplices
and assign\dss to every cochain $\alpha$ 
the number of\dss the simplices in\dss the sum $\alpha$ taken\dss modulo\dss $2$\nnsp.\oss
This defines a homomorphism\qss \vspace{0pt}
\[
\quad
\varepsilon\dff \colon\dff
C^{\fff n}\fff(\dff S\fff,\qff \partial\dff S \dff)
\qff \ttoo\off
\ftwo
\qff.
\]

\vspace{-12pt}
Since\dss $S$\dss is non-branching\halfff,\oss
$\varepsilon$\sss vanishes on\qss 
$\partial^*\dff \tau$\dss if\qss $\tau$\dss is\dss not\dss a simplex of\qss $\partial\dff S$\dss
and\dss hence defines a homomorphism\oss
$H^{\dff n}\dff(\dff S\fff,\qff \partial\dff S \dff)
\qff \ttoo\qff
\ftwo$\nsp.\oss
Obviously,\oss this homomorphism\dss maps every cohomology class\sss $\hclass{\sigma}$\sss to\sss $1$\nnsp.\oss
It\dss follows\dss that\dss the cohomology classes\sss $\hclass{\sigma}$\sss 
are non-zero.\oss  \eproof

\myuppar{Connecting\dss homomorphisms.}
Let\dss $S$\dss be a simplicial complex 
and\dss $Q$\dss be a subcomplex of\trs $S$\dnsp.\oss
Let\dss $m$\sss be a non-negative integer\halfff.\oss
Then there is a canonical\dss map\vspace{0pt}
\[
\quad
\partial^{*\nsp *}
\qff \colon\qff
H^{\fff m\dff -\dff 1}\fff(\dff Q \dff)
\off \ttoo\off
H^{\fff m}\fff(\dff S\fff,\qff Q \dff)
\qff,
\]

\vspace{-12pt}
called\dss the\qss \emph{connecting\dss homomorphism}\qss 
and defined as follows.\oss

To begin with,\oss
let\dss us consider a cochain\dss 
$\alpha\qff \in\qff C^{\fff m\dff -\dff 1}\fff(\dff Q \dff)$\dnsp.\oss
An\qss \emph{extension}\qss of\dss $\alpha$\dss
is any cochain\qss 
$\widetilde{\alpha}\pff \in\pff C^{\fff m\dff -\dff 1}\fff(\dff S \dff)$
resulting\dss from adding\dss to $\alpha$
several\dss $(\fff m\dff -\dff 1 \fff)$\dnsp-simplices of\dss $S$\dss
not\dss belonging\dss to\dss $Q$\nnsp.\oss
Every cochain $\alpha$ admits a\sss tautological\qss \emph{extension\dss by\dss zero},\oss
resulting\dss from adding\dss no simplices.\oss
But\halfff,\oss as we will see in a moment\halfff,\oss
the freedom\dss to use other extensions is essential.\oss
If\dss $\alpha$\dss is considered as a\sss linear\dss functional\qss
$C_{\fff m\dff -\dff 1}\fff(\dff Q \dff)\dff \ttoo\qff \ftwo$\nsp,\oss
then an\qss \emph{extension}\qss of\dss $\alpha$\dss
can\dss be defined as an extension of\dss $\alpha$\dss
to a linear\dss functional\qss
$C_{\fff m\dff -\dff 1}\fff(\dff S \dff)\dff \ttoo\qff \ftwo$\nsp.\oss

Let\qss
$a\qff \in\qff
H^{\fff m\dff -\dff 1}\fff(\dff Q \dff)$\nnsp.\oss
Then\qss
$a\off =\off \hclass{\alpha}$\qss
for some cocycle\qss
$\alpha\qff \in\qff
Z^{\fff m\dff -\dff 1}\fff(\dff Q \dff)$\dnsp.\oss
Let\qss 
$\widetilde{\alpha}\pff \in\pff C^{\fff m\dff -\dff 1}\fff(\dff S \dff)$\dss 
be an extension of\dss $\alpha$\qss
(usually\dss $\widetilde{\alpha}$\dss is not\sss a cocycle).\oss
The\dss $\partial^*\hnsp\partial^*$\dnsp\dnsp-the\-o\-rem\dss implies\sss that\dss the coboundary\dss
$\partial^*(\qff \widetilde{\alpha} \qff)$\dss
is a cocycle.\oss
In addition,\pss 
since\dss $\alpha$\dss is a cocycle,\qff\oss i.e.\qss 
$\partial^*(\dff \alpha \dff)\off =\off 0$\nnsp,\oss
the coboundary\dss
$\partial^*(\qff \widetilde{\alpha} \qff)$\dss
is a relative $m$\dnsp-cochain of\trs the pair\dss
$(\dff S\fff,\pff Q \dff)$\nnsp.\oss
Therefore\dss
$\partial^*(\qff \widetilde{\alpha} \qff)$\dss
is\sss a\sss relative cocycle.\oss
Let\vspace{2.75pt}
\[
\quad
\partial^{*\nsp *}(\dff a \trf)
\off =\off
\left[{\trf \partial^*(\qff \widetilde{\alpha} \qff) \dff}\right]
\off \in\off
H^{\fff m}\fff(\dff S\fff,\qff Q \dff)
\qff.
\]

\vspace{-10pt}
We need\dss to check\dss that\dss this definition is correct\halfff,\oss
i.e.\qss does not\dss depend on\dss the choices of\qss $\alpha\fff,\pff \widetilde{\alpha}$\nnsp.

\prooftitle{Proof\qss of\qss the correctness}
Any\dss two extensions of\dss $\alpha$\dss differ\dss by a relative cochain of\dss
$(\dff S\fff,\pff Q \dff)$\sss
and\dss hence\sss the coboundaries of\dss any\dss two extensions
differ\dss by a relative coboundary.\oss
This\dss implies\sss the\dss independence\dss on\dss the\dss choice\dss of\qss extension.\qff\oss
Let\pss $\alpha_{\dff 1}\dff,\off \alpha_{\dff 2}$\pss
be\sss two cocycles such\dss that\oss
$a
\off =~\off
\hclass{\fff\alpha_{\dff 1}}
\off =\off
\hclass{\fff\alpha_{\dff 2}}$\nnsp.\oss
Then\qss
$\alpha_{\dff 2}\qff -\qff \alpha_{\dff 1}
\off =\off
\partial^*(\fff \omega \dff)$\qss
for some\qss
$\omega\qff \in\qff
C^{\fff m\dff -\dff 2}\fff(\dff Q \dff)$\nnsp.\oss

Let\dss
$\widetilde{\alpha}_{\fff 1}$\dss
and\dss
$\widetilde{\omega}$\dss be arbitrary extensions of\trs
$\alpha_{\dff 1}$\dss and\dss $\omega$\dss respectively,\oss
and\dss let\vspace{3pt}
\[
\quad 
\widetilde{\alpha}_{\dff 2}
\off\qff =\off\qff
\widetilde{\alpha}_{\dff 1}
\off +\off
\partial^*(\dff \widetilde{\omega} \trf)
\qff.
\]

\vspace{-9pt}
Then\dss $\widetilde{\alpha}_{\dff 2}$\dss
is an extension of\dss $\alpha_{\dff 2}$\qss
(note\sss that\dss even\dss if\qss $\widetilde{\alpha}_{\dff 1}$\dss
and\dss $\widetilde{\omega}$\dss are extensions by zero,\oss
the extension\dss $\widetilde{\alpha}_{\dff 2}$\dss is usually\dss not).\oss 
Therefore\oss\vspace{3pt}
\[
\quad 
\partial^*(\qff \widetilde{\alpha}_{\fff 2} \qff)
\off =\off\qff
\partial^*\dff(\qff \widetilde{\alpha}_{\dff 1} \qff)
\off +\off\qff
\partial^* \circ\pff \partial^*\dff (\dff \widetilde{\omega} \dff) 
\off =\off\qff
\partial^*\dff(\qff \widetilde{\alpha}_{\dff 1} \qff)
\qff,
\]

\vspace{-9pt}
where at\dss the last\sss step we used\dss the\dss $\partial^*\hnsp\partial^*$\dnsp\dnsp-the\-o\-rem.\oss
Therefore\qss
$\dis
\left[{\trf \partial^*(\qff \widetilde{\alpha}_{\dff 2} \qff) \dff}\right]
\off =\off
\left[{\trf \partial^*(\qff \widetilde{\alpha}_{\dff 1} \qff) \dff}\right]$\nnsp.\oss
The independence on\dss the choice of\dss the cocycle\dss $\alpha$\dss follows.\oss  \eproof

\myuppar{Theorem\qss 7.}
\emph{Let\qss $S$\dss be a non-branching strongly connected simplicial\dss complex
of\qss dimension\dss $n$\nnsp.\oss 
Suppose\dss that\qss $\partial\dff S$\dss is
a non-branching\sss strongly\sss connected complex of\qss
dimension\qss $n\qff -\qff 1$\nnsp.\oss
Then\dss the connecting\dss homomorphism}\vspace{3pt}
\[
\quad
\partial^{*\nsp *}
\qff \colon\qff
H^{\fff n\dff -\dff 1}\fff(\dff \partial\dff S \dff)
\off \ttoo\off
H^{\fff n}\fff(\dff S\fff,\qff \partial\dff S \dff)
\qff
\]

\vspace{-9pt}
\emph{is an isomorphism.}

\proof
Let\dss $\tau$\dss be some $(\fff n\dff -\dff 1 \fff)$\dnsp-simplex of\dss $\partial\dff S$\dnsp,\oss
and\dss let\dss $\sigma$\dss be\sss the unique $n$\dnsp-simplex of\dss $S$\dss
such\dss that\dss $\tau$\dss is a face of\dss $\sigma$\dnsp.\oss
Let\dss $\widetilde{\tau}$\dss be\sss the extension of\dss the cochain\dss $\tau$\dss by\dss zero,\oss
i.e.\qss the same\dss $\tau$\nnsp,\oss
but\dss considered as a cochain of\dss $S$\dnsp.\oss
Then\qss
$\partial^*(\qff \widetilde{\tau} \qff)
\off =\off
\partial^*(\dff \tau\dff)
\off =\off
\sigma$\qss
and\dss hence\vspace{3pt}
\[
\quad
\partial^{*\nsp *}\bigl(\qff \hclass{\tau\fff} \qff\bigr)
\off =\off
\hclass{\sigma}
\qff.
\]

\vspace{-9pt}
But\dss by\qss Theorem\qss 6\qss
the cohomology classes\dss 
$\hclass{\tau\fff}$\dss and\dss $\hclass{\sigma}$\dss form\dss bases of\dss
the cohomology\dss groups\qss (vector spaces)\qss
$H^{\fff n\dff -\dff 1}\fff(\dff \partial\dff S \dff)$\dss and\dss
$H^{\fff n}\fff(\dff S\fff,\qff \partial\dff S \dff)$\dss
respectively.\oss
The\dss theorem\dss follows.\oss  \eproof

\myuppar{Induced\dss maps.}
Let\pss
$\varphi\dff \colon\dff S\ttoo S'$\pss 
be a simplicial\dss map.\oss
Theorem\qss 1*\qss implies\sss that\dss $\varphi^*$\dss
maps\dss $Z^{\fff m}\fff(\dff S' \dff)$\dss to\dss $Z^{\fff m}\fff(\dff S \dff)$\dss
and\dss maps\dss $B^{\fff m}\fff(\dff S' \dff)$\dss to\dss $B^{\fff m}\fff(\dff S \dff)$\nnsp.\oss
Therefore\dss $\varphi^*$\dss leads\sss to maps\vspace{3pt}
\[
\quad
\varphi^{*\nsp *}
\dff \colon\dff
H^{\fff m}\fff(\dff S' \dff)
\qff \ttoo\qff
H^{\fff m}\fff(\dff S \dff)
\]

\vspace{-9pt}
of\dss cohomology groups,\oss called\dss the\qss \emph{induced\dss maps}\qss in cohomology.\oss

Let\qss $Q\fff,\pff Q'$\qss be subcomplexes of\qss $S\fff,\pff S'$\qss
respectively.\oss
Suppose\sss that\dss $\varphi$\dss
is\qss \emph{simplicial\dss map of\dss pairs}\qss
$(\dff S\fff,\qff Q \dff)\qff \ttoo\qff (\dff S'\fff,\qff Q' \dff)$\nnsp,\oss
i.e.\qss that\dss $\varphi$\dss takes every simplex of\dss $Q$\dss to a simplex of\dss $Q'$\dnsp.\oss
Then\dss $\varphi$\dss defines a simplicial\dss map\qss
$\varphi_{\dff Q}\dff \colon\dff 
Q\qff \ttoo\qff Q'$\nnsp,\oss
and\dss the induced\dss map\dss $\varphi^*$\dss 
maps relative cochains of\dss
$(\dff S'\fff,\qff Q' \dff)$\dss
to relative cochains of\dss
$(\dff S\fff,\qff Q \dff)$\nnsp.\oss
Hence\dss $\varphi^*$\dss defines a maps\vspace{3pt}
\[
\quad
\varphi^*
\dff \colon\dff
C^{\fff m}\fff(\dff S'\fff,\qff Q' \dff)
\qff \ttoo\qff
C^{\fff m}\fff(\dff S\fff,\qff Q \dff)
\]

\vspace{-9pt}
of\dss relative cochains.\oss
Again,\oss Theorem\qss 1*\qss implies\sss that\dss $\varphi^*$\dss
maps relative cocycles\sss to relative cocycles and\dss
relative coboundaries\sss to relative coboundaries and\dss hence leads\sss to maps\vspace{1.5pt}
\[
\quad
\varphi^{*\nsp *}
\dff \colon\dff
H^{\fff m}\fff(\dff S'\fff,\qff Q' \dff)
\qff \ttoo\qff
H^{\fff m}\fff(\dff S\fff,\qff Q \dff)
\]

\vspace{-10.5pt}
of\dss relative cohomology groups.\oss
They are also called\dss the\qss \emph{induced\dss maps}.\oss
The maps induced\dss by\dss $\varphi$\dss and\dss $\varphi_{\dff Q}$\dss
s together\dss with connecting homomorphisms\sss form\dss 
the following diagram.\vspace{3pt}
\begin{equation*}
\quad
\begin{tikzcd}[column sep=spec, row sep=spe]\dis
H^{\fff m\dff -\dff 1}\fff(\dff Q' \dff)\dff\off
\arrow[r, 
"\dis \partial^{*\nsp *}"]
\arrow[d, 
"\dis \varphi^{*\nsp *}_{\dff Q}"']
& 
\off H^{\fff m}\fff(\dff S'\fff,\pff Q' \dff)
\arrow[d,   
"\dis \varphi^{*\nsp *}"']
\\
H^{\fff m\dff -\dff 1}\fff(\dff Q \dff)\off
\arrow[r,  
"\dis \partial^{*\nsp *}"]
& 
\off H^{\fff m}\fff(\dff S\fff,\pff Q \dff)
\end{tikzcd}
\end{equation*}

\vspace{-6.5pt}
\myuppar{Lemma.}
\emph{The above diagram is commutative,\oss i.e.}\oss
$\varphi^{*\nsp *}\fff \circ\pff \partial^{*\nsp *}
\off =\off 
\partial^{*\nsp *}\fff \circ\pff \varphi^{*\nsp *}_{\dff Q}$\nnsp.\oss

\proof
If\qss
$a\qff \in\qff
H^{\fff m\dff -\dff 1}\fff(\dff Q' \dff)$\nnsp,\oss
then\qss
$a\off =\off \hclass{\alpha}$\qss
for some\qss
$\alpha\qff \in\qff
Z^{\fff m\dff -\dff 1}\fff(\dff Q' \dff)$\dnsp.\oss
Let\dss\vspace{1.5pt}
\[
\quad
\widetilde{\alpha}\pff \in\pff C^{\fff m\dff -\dff 1}\fff(\dff S' \dff)
\]

\vspace{-10.5pt}
be an extension\sss of\dss $\alpha$\nnsp.\oss
Then\oss 
$\dis
\partial^{*\nsp *}(\dff a \trf)
\off =\off\dff
\left[{\trf \partial^*(\qff \widetilde{\alpha} \qff) \dff}\right]
$\oss 
and\dss hence\oss \vspace{4.5pt}
\[
\quad
\varphi^{*\nsp *}\fff \circ\pff \partial^{*\nsp *}(\dff a \trf)
\off =\off\dff
\left[{\trf \varphi^*\fff \circ\pff \partial^*(\qff \widetilde{\alpha} \qff) \dff}\right]
\qff.
\]

\vspace{-7.5pt}
On\dss the other\dss hand,\oss
$\dis
\varphi^{*\nsp *}_{\dff Q}\dff(\dff a \trf)
\off =\off
\left[{\trf \varphi^*_{\dff Q}\dff(\dff \alpha \dff) \dff}\right]$\qss
and\dss
$\varphi^*\dff(\qff \widetilde{\alpha}\qff)$\qss
is an extension of\dss
$\varphi^*_{\dff Q}\dff(\dff \alpha \dff)$\nnsp.\oss
Therefore\vspace{4.5pt}
\[
\quad
\partial^{*\nsp *}\fff \circ\pff \varphi^{*\nsp *}_{\dff Q}\dff(\dff a \trf)
\off =\off\dff
\left[{\trf \partial^*\fff \circ\pff \varphi^*\dff(\qff \widetilde{\alpha} \qff) \dff}\right]
\qff.
\]

\vspace{-7.5pt}
Theorem\qss 1*\pss implies\sss that\qss
$\dis
\varphi^*\fff \circ\pff \partial^*(\qff \widetilde{\alpha} \qff)
\off =\off
\partial^*\fff \circ\pff \varphi^*\dff(\qff \widetilde{\alpha} \qff)$\qss
and\dss hence\vspace{3pt}
\[
\quad
\varphi^{*\nsp *}\fff \circ\pff \partial^{*\nsp *}(\dff a \trf)
\off =\off
\partial^{*\nsp *}\fff \circ\pff \varphi^{*\nsp *}(\dff a \trf)
\qff.
\]

\vspace{-9pt}
The\sss lemma\dss follows.\oss  \eproof

\myuppar{Remark.}
Suppose\sss that\dss we\sss took as\dss $\widetilde{\alpha}$\dss
the extension\dss by\dss zero.\oss
If\dss $\varphi$\dss maps\dss to\dss $Q'$\dss
some simplices of\dss $S$\dss not\dss belonging\dss to\dss $Q$\nnsp,\oss
then\dss $\varphi^*\dff(\qff \widetilde{\alpha}\qff)$\qss usually\dss
will\dss not\dss be an extension\dss by\dss zero of\dss
$\varphi^*_{\dff Q}\dff(\dff \alpha \dff)$\nnsp.\oss
This is another\dss illustration of\dss the usefulness of\dss the freedom\dss
in\dss the choice of\dss extensions,\oss
and\dss we will encounter\dss this situation\dss in\dss the following cohomological\dss
proof\dss of\qss Sperner's\dss lemma.\oss
The coboundaries\dss $\partial^*(\dff \tau\dff)$\dss to be discarded
are discarded\dss here.\oss

\prooftitle{A cohomological\dss proof\qss  
of\pss Sperner's\trs lemma\fff}
Let\qss
$\varphi\dff \colon\dff
T\qff \ttoo\qff \Delta$\qss
be a pseudo-i\-den\-ti\-cal simplicial\dss map.\oss
Since\dss $\varphi$\dss is a pseudo-identical\dss map,\oss
$\varphi$\dss leads\sss to\sss the map\qss 
$\varphi_{\dff \partial\fff T}
\dff \colon\dff
\partial\dff T
\qff \ttoo\qff
\partial\dff \Delta$\nnsp,\oss
which we will\sss denote now\dss by\dss $\partial\fff \varphi$\nnsp.\oss
Let\dss us consider\dss the following diagram.\oss\vspace{2.5pt}
\begin{equation*}
\quad
\begin{tikzcd}[column sep=spec, row sep=spe]\dis
H^{\fff n\dff -\dff 1}\fff(\dff \partial\dff \Delta \dff)\dff\off
\arrow[r, 
"\dis \partial^{*\nsp *}"]
\arrow[d, 
"\dis (\dff\partial\fff \varphi\dff)^{\nsp *\nsp *}"']
& 
\off H^{\fff n}\fff(\dff \Delta\fff,\pff \partial\dff \Delta \dff)
\arrow[d,   
"\dis \varphi^{*\nsp *}"']
\\
H^{\fff n\dff -\dff 1}\fff(\dff \partial\dff T \dff)\off
\arrow[r,  
"\dis \partial^{*\nsp *}"]
& 
\off H^{\fff n}\fff(\dff T\fff,\pff \partial\dff T \dff)
\end{tikzcd}
\end{equation*}

\vspace{-6.5pt}
By\dss the above lemma it\dss is commutative.\oss
By\qss Theorem\qss 6\qss every cohomology\dss group in\dss
this diagram is a vector space of\dss dimension one
with\dss the cohomology class of\dss any\dss top-dimensional\sss simplex
forming a basis\qss
(or\halfff,\oss what\dss is\dss the same,\oss being\dss the only\dss
non-zero element).\oss
In\dss par\-tic\-u\-lar\halfff,\oss every cohomology group in\dss this diagram\dss
is isomorphic\dss to\dss $\ftwo$\nsp.\oss
If\dss a vector space over\dss $\ftwo$\dss is isomorphic\sss to\dss 
$\ftwo$\nsp,\oss
then\dss it\dss is canonically\dss isomorphic\sss to\dss $\ftwo$\nsp.\oss
Hence we can\dss replace all cohomology\dss groups in our diagram\dss by\dss
$\ftwo$\dss and\dss get\dss the following\sss diagram.\oss\vspace{2.25pt}
\begin{equation*}
\quad
\begin{tikzcd}[column sep=specc, row sep=spe]\dis
\ftwo 
\arrow[r, 
"\dis \partial^{*\nsp *}"]
\arrow[d, 
"\dis (\dff\partial\dff \varphi\dff)^{\nsp *\nsp *}"']
& 
\ftwo
\arrow[d,   
"\dis \varphi^{*\nsp *}"']
\\
\ftwo
\arrow[r,  
"\dis \partial^{*\nsp *}"]
& 
\ftwo
\end{tikzcd}
\end{equation*}

\vspace{-6.75pt}
Let\dss us denote by\dss $\mathbb{1}$\dss the non-zero element\sss of\dss $\ftwo$\nsp.\oss 
Theorem\qss 7\qss implies\sss that\dss both connecting\dss homomorphisms\dss
$\partial^{*\nsp *}$\dss are isomorphisms,\oss
i.e.\qss that\qss
$\partial^{*\nsp *}(\dff \mathbb{1}\qff)
\off =\off
\mathbb{1}$\qss 
for both\dss maps\dss $\partial^{*\nsp *}$\nnsp.\oss
By\dss the definition,\oss
$\varphi^*(\dff \delta\trf)$\dss
is equal\dss to\dss the sum of\dss all $n$\dnsp-simplices $\sigma$ of\trs $T$\dss
such\dss that\qss
$\varphi\dff(\dff \sigma \dff)\off =\off \delta$\nnsp.\oss
It\dss follows\dss that\qss
$\varphi^{*\nsp *}(\dff \mathbb{1}\qff)\off =\off e\trf \mathbb{1}$\nnsp,\oss
where $e$\sss is\dss the number of\dss such simplices\dss $\sigma$\dnsp,\oss
and\dss hence\vspace{3pt}
\[
\quad
\varphi^{*\nsp *}\fff \circ\pff \partial^{*\nsp *}\fff (\dff \mathbb{1} \qff)
\off =\off
e\dff \mathbb{1}
\qff.
\]

\vspace{-9pt}
Similarly,\oss if\qss $i\qff \in\qff I$\nnsp,\oss
then\dss $(\dff\partial\fff \varphi\dff)^{*}\dff (\dff \delta_{\dff i}\dff)$\dss
is equal\dss to\dss the sum of\dss all $(\fff n\dff -\dff 1 \fff)$\dnsp-simplices 
$\tau$ of\trs $\partial\dff T$\dss
such\dss that\qss
$\partial\dff\varphi\dff(\dff \tau \dff)
\off =\off
\varphi\dff(\dff \tau \dff)
\off =\off 
\delta_{\dff i}$\nsp.\oss
Since\dss $\varphi$\dss is a pseudo-identical\dss map,\oss
every such simplex $\tau$ belongs\sss to\dss $T_{\dff i}$\nnsp.\oss
It\dss follows\dss that\qss
$(\dff\partial\fff \varphi\dff)^{\nsp *\nsp *}(\dff \mathbb{1}\qff)
\pff =\off 
h\trf \mathbb{1}$\nnsp,\oss
where $h$\sss is\dss the number of\dss simplices of\trs $T_{\dff i}$\dss
such\dss that\qss
$\varphi\dff(\dff \tau \dff)
\off =\off 
\delta_{\dff i}$\nsp,\oss
and\dss hence\vspace{3pt}
\[
\quad
\partial^{*\nsp *}\fff \circ\pff (\dff\partial\fff \varphi\dff)^{\nsp *\nsp *}\fff (\dff \mathbb{1} \qff)
\off =\off
h\trf \mathbb{1}
\qff.
\]

\vspace{-9pt}
Now\dss the commutativity of\dss the last\sss diagram\dss implies\sss that\qss
$e\trf \mathbb{1}\off =\off h\trf \mathbb{1}$\qss
and\dss hence\qss
$e\off \equiv\off h$\qss modulo\dss $2$\nnsp.\oss
As usual,\oss
an\dss induction\dss by $n$
completes\sss the proof\halfff.\oss  \eproof

\myappend{Barycentric\qss subdivisions}{barycent}

\myuppar{Cones.}
Let\dss $X$\dss be a subset\sss of\trs $\rrr^{\dff d}$\dss
and\dss let\sss $z$\sss be a point\dss in\dss  $\rrr^{\dff d}$\nnsp.\oss
Suppose\sss that\dss segments connecting\dss $z$\dss with different\dss 
points of\dss $X$\dss always intersect\sss only\sss at\sss $z$\nnsp.\oss
Then\dss the union\dss $z\dff *\dff X$\dss of\dss all segments connecting\dss $z$\dss with\dss
points of\dss $X$\sss is called\dss the\qss
\emph{cone}\qss over\dss $X$\dss with\dss the\qss \emph{apex}\qss $z$\nnsp.\oss

If\dss $X$\dss is a simplex and\dss $z$\sss is affinely\dss
independent\dss from\dss the vertices of\dss $X$\nnsp,\oss
then\dss $z\dff *\dff X$\dss is also a simplex\halfff.\oss
Its vertices are\sss the vertices of\dss $X$\dss together with\dss the point\sss $z$\nnsp.\oss 
The faces of\dss $z\dff *\dff X$\dss
are\sss the vertex\sss $z$\nnsp,\oss 
the faces\dss $Y$\dss of\dss $X$\dss and\dss the cones\dss
$z\dff *\dff Y$\dss over\dss the faces\dss $Y$\dss of\dss $X$\nnsp.\oss

\myuppar{Geometric simplices as cones over\dss their boundaries.}
For a\sss geometric simplex\sss $\sigma$\sss we will\sss denote by\dss
$\sco{\fff\sigma}$\dss the simplicial\sss complex consisting\sss
of\dss $\sigma$\dss and all\dss its faces.\oss
The boundary\dss $\partial\fff\sco{\fff\sigma}$\dss is\dss the\sss
complex consisting of\dss all\dss proper\dss faces of\dss $\sigma$\dnsp.\oss
Clearly,\pss $\norm{\sco{\fff\sigma}}\off =\off \sigma$\qss and\qss
$\norm{\partial\fff\sco{\fff\sigma}\nsp}\off =\off \bd\dff \sigma$\dnsp.\oss

Let\dss $\sigma$\dss be a geometric simplex and
let\qss $z\qff \in\qff \sigma\qff \smallsetminus\qff \bd\dff \sigma$\qss
be a point\dss in\dss the interior of\dss $\sigma$\dnsp.\oss
Then segments connecting\sss $z$\sss with different\dss points of\dss the 
boundary\dss $\bd\dff \sigma$\dss intersect\sss only\sss at\sss $z$\sss
and\dss hence\dss $\sigma$\dss is a cone over\dss $\bd\dff \sigma$\dss
with\dss the apex\sss $z$\nnsp,\oss
i.e.\qss
$\sigma\off =\off z\dff *\trf \bd\dff \sigma$\nnsp.\oss

Let\dss $S$\dss be a subdivision of\dss the complex\dss $\partial\fff\sco{\fff\sigma}$\dnsp.\oss
Then\qss $\norm{S}\off =\off \norm{\partial\fff\sco{\fff\sigma}\nsp}\off =\off \bd\dff \sigma$\nnsp,\oss
and every simplex of\trs $S$\dss is contained\dss in a proper\dss face of\dss $\sigma$\dnsp.\oss
It\dss follows\sss for every simplex\dss $\tau$\dss of\dss $S$\dss
the point\sss $z$\sss is affinely\dss independent\dss from\dss the vertices of\dss $\tau$\dss
and\dss hence\sss the cone\dss 
$z\dff *\dff \tau$\dss is defined\sss and\dss is a geometric simplex.\oss
Let\sss us define\sss the\qss \emph{cone}\dss $z\dff *\dff S$\dss
with\dss the\qss \emph{apex}\dss $z$\dss as\sss the collection consisting of\dss $z$\nnsp,\oss
the simplices of\dss $S$\nnsp,\oss and\dss the cones\dss
$z\dff *\dff \tau$\dss 
over\dss the simplices\dss $\tau$\dss of\dss $S$\nnsp.\oss
Clearly,\qss $z\dff *\dff S$\dss
is a simplicial\sss complex and\dss is a subdivision\sss of\trs $\sco{\fff\sigma}$\dnsp.\oss
The complex\sss $S$\sss is a sub\-complex of\dss  $z\dff *\dff S$\dss
and\sss every\sss vertex of\dss $z\dff *\dff S$\dss
is\sss either $z$ or\dss is a vertex of\sss $S$\dnsp.

\myuppar{Centers-generated subdivisions.} 
Let\dss $S$\dss be a geometric simplicial\sss complex\halfff.\oss
Suppose\sss that\dss for every simplex\dss $\sigma$\dss of\qss $S$\dss
a point\qss $z\dff(\dff \sigma\dff)\qff \in\qff \sigma\qff \smallsetminus\qff \bd\dff \sigma$\qss
is chosen.\oss
One may\dss think\dss that\dss $z\dff(\dff \sigma\dff)$\dss is\dss a\sss sort\sss of\dss
a\qss \emph{center}\qss of\trs $\sigma$\dnsp.\oss
Any choice of\dss such centers generates a subdivision\dss $\cs S$\dss of\qss $S$\dss as\sss follows.\oss

For every\dss integer\qss $n\qff \geq\qff 0$\qss let\dss $S_{\fff n}$\dss
be\sss the complex consisting of\dss simplices of\qss $S$\dss
having dimension\qss $\leq\qff n$\nnsp.\oss
Then\dss $S_{\fff 0}$\dss is a finite set\sss and\qss
$S_{\fff n}\off =\off S$\qss for all\sss sufficiently\dss big\dss $n$\nnsp.\oss
Let\dss us consecutively construct\dss the subdivisions\qss
$\cs S_{\fff 0}\dff,\off \cs S_{\fff 1}\dff,\off \cs S_{\fff 2}\dff,\off \ldots$\qss
of\oss
$S_{\fff 0}\dff,\off S_{\fff 1}\dff,\off S_{\fff 2}\dff,\off \ldots$\qss
respectively.\oss
Let\qss $\cs S_{\fff 0}\off =\off S_{\fff 0}$\nnsp.\oss
Suppose\sss that\dss the subdivision\dss $\cs S_{\fff m}$\dss of\dss $S_{\fff m}$\dss
is already constructed.\oss
Let\dss $\sigma$\dss be an $(\fff m\dff +\dff 1\fff)$\dnsp-simplex\dss $\sigma$\dss
of\qss $S$\dnsp.\oss
Its\trs boundary\dss $\bd\dff \sigma$\dss is contained\dss in\qss 
$\norm{\cs S_{\fff m}}\off =\off \norm{S_{\fff n}}$\nnsp,\oss
and\sss simplices of\qss $\cs S_{\fff m}$\dss contained\dss in\dss
$\bd\dff \sigma$\dss form\sss a\sss subdivision\sss of\qss
$\partial\fff\sco{\fff\sigma}$\dnsp.\oss
Let\dss us\sss denote\sss this subdivision\dss by\dss
$\cs \partial\fff\sco{\fff\sigma}$\dnsp.\oss
The subdivision\dss $\cs S_{\fff m\dff +\dff 1}$\dss of\pss $S_{\fff m\dff +\dff 1}$\dss
is\dss the result\sss of\trs adding\dss to\dss
$\cs S_{\fff m}$\dss the cone\dss\vspace{1.5pt}
\[
\quad
z\dff(\dff \sigma\dff)\trf *\qff \cs \partial\fff\sco{\fff\sigma}
\]

\vspace{-10.5pt}
for each\dss
$(\fff m\dff +\dff 1\fff)$\dnsp-simplex\dss $\sigma$\dss
of\qss $S$\dnsp.\oss
Clearly,\oss $\cs S_{\fff m\dff +\dff 1}$\dss is a simplicial\sss complex\sss
and\dss is\sss a subdivision of\pss $S_{\fff m\dff +\dff 1}$\nnsp.\oss
Finally,\oss $\cs S$\dss is defined as\dss
$\cs S_{\fff m}$\dss for any\dss $m$\dss such\dss that\qss
$S_{\fff m}\off =\off S$\dnsp.\oss

By\dss the construction,\oss every\sss simplex\dss $\omega$\dss of\qss $\cs S_{\fff m\dff +\dff 1}$\dss
not\sss contained\dss in\dss $\cs S_{\fff m}$\dss
has\sss the center\dss $z\dff(\dff \sigma\dff)$\dss of\trs
some\dss $(\fff m\dff +\dff 1\fff)$\dnsp-simplex\dss $\sigma$\dss as a vertex\halfff,\oss
and\dss the other vertices of\dss $\omega$\dss are all\sss contained\dss
in a face of\dss $\sigma$\dnsp.\oss
It\dss follows\sss that\sss every\sss $n$\dnsp-simplex\dss $\omega$\dss of\dss $\cs S$\dss 
has\sss as\sss its\sss vertices\sss the centers\qss\vspace{0pt}
\[
\quad
z\dff(\dff \sigma_{\dff 0}\dff)\fff,\off
z\dff(\dff \sigma_{\dff 1}\dff)\fff,\off
\ldots\fff,\off
z\dff(\dff \sigma_{\dff n}\dff)
\]

\vspace{-12pt}
of\dss simplices\qss
$\sigma_{\dff 0}\fff,\off
\sigma_{\dff 1}\fff,\off
\ldots\fff,\off
\sigma_{\dff n}$\qss
such\dss that\qss
$\sigma_{\dff i\dff -\dff 1}$\dss
is a face of\dss $\sigma_{\dff i}$\qss
for all\qss
$i\off =\off 1\fff,\off 2\fff,\off \ldots\fff,\off n$\qss
and\dss hence\dss $\sigma_{ j}$\dss
is\dss a\dss face of\dss $\sigma_{\dff i}$\dss if\qss $j\qff \leq\qff i$\nnsp.\oss
In\dss particular\halfff,\oss
the combinatorial\sss structure of\dss $\cs S$\dss does not\sss depend\sss
on\dss the choice of\dss centers\dss $z\dff(\dff \sigma\dff)$\dnsp.
But\dss the size of\dss simplices of\dss $\cs S$\dss
depends on\dss the choice of\dss the centers.\oss
Let\dss us\dss turn\dss to\sss an\sss efficient\dss in\dss this
respect\sss choice of\dss centers.\oss

\myuppar{Barycentric subdivisions.}
Recall\qss (see \qss Section\qss \ref{dimension})\qss
that\dss the\qss \emph{barycenter}\qss of\dss a geometric simplex\dss $\sigma$\dss is\dss the only\dss point\sss
of\dss $\sigma$\dss with all\dss barycentric coordinates equal.\oss
So,\oss if\qss
$w_{\fff 0}\fff,\pff w_{\fff 1}\fff,\pff \ldots\fff,\pff w_{\fff m}$\qss
are\dss the\sss vertices of\dss $\sigma$\dnsp,\oss then\dss
the\qss \emph{barycenter}\qss of\dss $\sigma$\dss is\sss equal\dss to\sss the point\vspace{1.5pt}
\[
\quad
c\dff(\dff \sigma\dff)
\off\off =\off\off
\frac{1}{m\qff +\qff 1}\off \sum\nolimits_{i\qff =\qff 0}^{m}\qff w_{\fff i}
\qff.
\]

\vspace{-10.5pt}
Let\dss $S$\dss be a geometric simplicial\sss complex\halfff,\oss
and\dss let\dss us choose\sss the barycenters as\sss the centers of\dss simplices.\oss
In\sss other\dss words,\oss
let\qss $z\dff(\dff \sigma\dff)\off =\off c\dff(\dff \sigma\dff)$\nnsp.\oss
The corresponding\sss subdivision\dss $\cs S$\dss is called\dss
the\qss \emph{barycentric\dss subdivision}\qss of\qss $S$\dss
and\dss is\sss denoted\dss by\dss $\bs S$\dnsp.\oss
As we will see now,\oss
passing\dss from\dss $S$\dss to\dss $\bs S$\dss
decreases\dss the maximal\sss diameter of\dss simplices by\sss
a definite factor\halfff.\oss
For\qss $x\qff \in\qff \rrr^{\dff d}$\dnsp,\oss
let\dss us denote\sss by\dss $\num{x}$\dss be\sss the norm of\dss the vector\dss $x$\nnsp.\oss
So,\pss $\num{x\qff -\qff y}$\dss is\sss the distance between\qss
$x\fff,\pff y\qff \in\qff \rrr^{\dff d}$\dnsp.

\myuppar{Lemma.}
\emph{Let\dss $\sigma$\dss be an\dss $n$\dnsp-simplex\dss in\qss $\rrr^{\dff d}$\dnsp,\oss
and\dss let\qss
$w_{\fff 0}\fff,\pff w_{\fff 1}\fff,\pff \ldots\fff,\pff w_{\fff n}$\qss
be its vertices.\oss
Then}\vspace{0pt}
\[
\quad
\num{x\qff -\qff a}
\off \leq\off
\max_{\dff i}\qff \num{x\qff -\qff w_{\fff i}}
\]

\vspace{-15pt}
\emph{for every\qss $x\qff \in\qff \rrr^{\dff d}$\qss
and\qss $a\qff \in\qff \sigma$\nnsp.\oss
The diameter of\dss $\sigma$\dss is equal\qss to\qss
$\max_{\dff i\fff,\trf j}\qff \num{w_{\fff i}\qff -\qff w_{\fff j}}$\nsp.\oss}

\proof
By\dss the definition of\dss a simplex\halfff,\oss\vspace{1.5pt}
\[
\quad
a
\off =\off
\sum\nolimits_{i\qff =\qff 0}^{n}\qff a_{\fff i}\dff w_{\fff i}
\qff,
\]

\vspace{-10.5pt}
where\qss $a_{\fff i}\qff \geq\qff 0$\qss for all\dss $i$\dss
and\oss
$\dis
\sum\nolimits_{i\qff =\qff 0}^{n}\qff a_{\fff i}
\off =\off
1$\nnsp.\oss
It\dss follows\dss that\vspace{3pt}
\[
\quad
\num{x\qff -\qff a}
\off =\off
\num{x\off -\off
\sum\nolimits_{i\qff =\qff 0}^{n}\qff a_{\fff i}\dff w_{\fff i}}
\]

\vspace{-33pt}
\[
\quad
\phantom{\num{x\qff -\qff a}
\off }
=\off
\num{\sum\nolimits_{i\qff =\qff 0}^{n}\qff a_{\fff i}\dff (\dff x\qff -\qff w_{\fff i}\dff)}
\off \leq\off\dff
\sum\nolimits_{i\qff =\qff 0}^{n}\qff a_{\fff i}\dff
\num{x\qff -\qff w_{\fff i}}
\]

\vspace{-9.5pt}
and\dss hence\qss
$\num{x\qff -\qff a}
\off \leq\off 
\max_{\dff i}\qff \num{x\qff -\qff w_{\fff i}}
\off =\off
\max_{\dff i}\qff \num{w_{\fff i}\qff -\qff x}$\nnsp.\oss
This proves\sss the first\sss statement\sss of\dss the lemma.\oss
If\qss $x\qff \in\qff \sigma$\nnsp,\oss
then\dss the first\sss statement\dss implies\sss that\qss
$\num{w_{\fff i}\qff -\qff x}
\off \leq\off
\max_{\dff j}\qff \num{w_{\fff i}\qff -\qff w_{\fff j}}$\qss
for every\dss $i$\dss
and\dss hence\qss
$\num{x\qff -\qff a}
\off \leq\off
\max_{\dff i\fff,\trf j}\qff \num{w_{\fff i}\qff -\qff w_{\fff j}}$\nnsp.\oss
The second\sss statement\dss follows.\oss  \eproof

\myuppar{Lemma.}
\emph{Let\dss $\sigma$\dss be an\dss $n$\dnsp-simplex\halfff.\oss 
Then\dss the diameter of\dss every\dss simplex of\dss the 
barycentric subdivision\qss $\bs \sco{\fff\sigma}$\qss is\qss
$\leq$\qss than\qss $n/(\dff n\qff +\qff 1 \dff)$\qss
times\sss the diameter of\qss $\sigma$\dnsp.\oss}\oss

\proof
Let\qss
$w_{\fff 0}\fff,\pff w_{\fff 1}\fff,\pff \ldots\fff,\pff w_{\fff n}$\qss
be\dss the vertices of\dss $\sigma$\dss
and\dss let\dss $\tau$\dss be a simplex of\qss $\bs \sco{\fff\sigma}$\dnsp.\oss
Every\dss vertex of\dss $\tau$\dss is\dss the barycenter of\dss
some face of\dss $\sigma$\dnsp.\oss
Moreover\halfff,\oss
if\qss $v_{\fff 1}\fff,\pff v_{\fff 2}$\qss are\sss two vertices of\dss $\tau$\dss
and are\sss the\sss barycenters of\dss the faces\qss $\sigma_{\dff 1}\fff,\pff \sigma_{\dff 2}$\qss
respectively,\oss
then one of\dss the simplices\qss $\sigma_{\dff 1}\fff,\pff \sigma_{\dff 2}$\qss
is\dss a face of\dss the other\halfff.\oss
Without\sss any\dss loss of\dss generality\dss we may assume\sss that\dss
$\sigma_{\dff 1}$\dss is a face of\dss $\sigma_{\dff 2}$\nsp.\oss
After\dss renumbering\dss the vertices,\pss if\dss necessary,\pss
we now\sss may assume\sss that\qss
$w_{\fff 0}\fff,\pff w_{\fff 1}\fff,\pff \ldots\fff,\pff w_{\fff p}$\qss
are\sss the vertices of\dss $\sigma_{\dff 1}$\dss
and\dss
$w_{\fff 0}\fff,\pff w_{\fff 1}\fff,\pff \ldots\fff,\pff w_{\fff q}$\qss
are\sss the vertices of\dss $\sigma_{\dff 2}$\dss
for some\qss $p\qff \leq\qff q$\nnsp.\oss
Then\vspace{3pt}
\[
\quad
v_{\fff 1}
\off\qff =\off\qff
c\dff(\dff \sigma_{\dff 1}\dff)
\off\qff =\off\qff
\frac{1}{p\qff +\qff 1}\off \sum\nolimits_{i\qff =\qff 0}^{p}\qff w_{\fff i}
\hspace{1.5em}\mbox{and}\hspace{1.5em}
v_{\fff 2}
\off\qff =\off\qff
c\dff(\dff \sigma_{\dff 2}\dff)
\off\qff =\off\qff
\frac{1}{q\qff +\qff 1}\off \sum\nolimits_{i\qff =\qff 0}^{q}\qff w_{\fff i}
\qff.
\]

\vspace{-12.5pt}
By\dss the\sss last\dss lemma,\vspace{2.5pt}
\[
\quad
\num{c\dff(\dff \sigma_{\dff 2}\dff)
\qff -\qff
c\dff(\dff \sigma_{\dff 1}\dff)}
\off \leq\off
\max_{\dff 0\qff \leq\qff i\qff \leq\qff p}\qff \num{c\dff(\dff \sigma_{\dff 2}\dff)\qff -\qff w_{\fff i}}
\qff.
\]

\vspace{-12.5pt}
At\dss the same\sss time,\vspace{-0.25pt}
\[
\quad
\num{w_{\fff i}\qff -\qff c\dff(\dff \sigma_{\dff 2}\dff)}
\off =\off
\num{w_{\fff i}
\qff -\qff
\frac{1}{q\qff +\qff 1}\off \sum\nolimits_{j\qff =\qff 0}^{q}\qff w_{\fff j}}
\off \leq\off
\frac{1}{q\qff +\qff 1}\off \sum\nolimits_{j\qff =\qff 0}^{q}\qff
\num{w_{\fff i}\qff -\qff w_{\fff j}}
\qff. 
\]

\vspace{-9pt}
If\qss $i\qff \leq\qff p$\nnsp,\oss
then\dss the last\dss sum\dss includes\sss the summand\qss
$\num{w_{\fff i}\qff -\qff w_{\fff i}}
\off =\off
0$\qss
and\sss hence\vspace{4.5pt}
\[
\quad
\num{w_{\fff i}\qff -\qff c\dff(\dff \sigma_{\dff 2}\dff)}
\off \leq\off
\frac{q}{q\qff +\qff 1}\off
\max_{\dff i\fff,\trf j}\qff \num{w_{\fff i}\qff -\qff w_{\fff j}}
\off =\off
\frac{q}{q\qff +\qff 1}\off r
\qff,
\]

\vspace{-9pt}
where\qss 
$r
\off =\off
\max_{\dff i\fff,\trf j}\qff \num{w_{\fff i}\qff -\qff w_{\fff j}}$\qss 
is\sss the diameter of\dss $\sigma$\dnsp.\oss
It\dss follows\dss that\qss\vspace{3pt}
\[
\quad
\num{c\dff(\dff \sigma_{\dff 2}\dff)
\qff -\qff
c\dff(\dff \sigma_{\dff 1}\dff)}
\off \leq\off
\frac{q}{q\qff +\qff 1}\off r
\qff.
\]

\vspace{-9pt}
Finally,\pss $q\qff \leq\qff n$\qss implies\dss that\oss
$\dis
q/(\dff q\qff +\qff 1\dff)
\off \leq\off
n/(\dff n\qff +\qff 1\dff)$\oss
and\dss hence\vspace{3pt}
\[
\quad
\num{v_{\fff 2}\qff -\qff v_{\fff 1}}
\off =\off
\num{c\dff(\dff \sigma_{\dff 2}\dff)
\qff -\qff
c\dff(\dff \sigma_{\dff 1}\dff)}
\off \leq\off
\frac{n}{n\qff +\qff 1}\off r
\]

\vspace{-9pt}
Since\qss $v_{\fff 1}\fff,\pff v_{\fff 2}$\qss are\sss two arbitrary\sss
vertices of\dss $\tau$\nnsp,\oss
it\dss remains\sss to apply\dss the last\dss lemma.\oss  \eproof

\myuppar{Iterated\dss barycentric subdivisions.}
Let\sss $S$\sss be a geometric simplicial\sss complex\halfff.\oss
The\qss \emph{iterated\trs barycentric\dss subdivisions}\pss $\ibs{i}\fff S$\dss of\qss $S$\dss
is\dss defined\dss by\dss the rules\qss $\ibs{0} S\off =\off S$\qss
and\qss $\ibs{i}\fff S\off =\off \bs (\trf \ibs{i\dff -\dff 1} S\dff)$\qss
for every\qss $i\qff \geq\qff 1$\nnsp.\oss
Clearly,\oss the complexes\dss $\ibs{i}\fff S$\dss are indeed\sss subdivisions of\qss $S$\dnsp.\oss
Now,\oss let\qss $\varepsilon\qff >\qff 0$\nnsp.\oss
Since\qss $n/(\dff n\qff +\qff 1\dff)\qff <\qff 1$\nnsp,\oss
the\sss last\dss lemma implies\sss that\dss  
the diameters of\dss simplises\dss of\qss $\ibs{i}\fff S$\dss
are\qss $<\qff \varepsilon$\qss for all\sss sufficiently\dss big\dss $i$\nnsp.\oss

\myappend{Stellar\qss subdivisions}{stellar}

\myuppar{Stellar\sss subdivisions.}
Let\dss $S$\dss be a simplicial complex and\dss $\sigma$\dss be\dss a simplex of\trs $S$\nnsp.\oss
Let\dss $K$\dss be\sss the subcomplex of\trs $S$\sss consisting of\dss simplices of\trs $S$\sss
having $\sigma$ as a face\sss together\dss with\dss their other\dss faces.\oss
Let\dss $L$\dss be\sss the subcomplex of\trs $K$\dss consisting of\dss simplices of\trs $K$\dss
not\dss having $\sigma$ as a face.\oss

Let\dss us\sss choose a\sss point\dss $w$\dss in\dss the interior\dss
$\sigma\qff \smallsetminus\qff \bd\fff \sigma$\dss of\qss $\sigma$\dss
and\sss define\sss the\qss \emph{cone}\qss 
$K'\off =\off w\dff *\dff L$\qss
as\sss the collection of\dss simplices consisting of\dss $w$\nnsp,\oss
the simplices of\dss $L$\nnsp,\oss and\dss the cones\dss
$w\dff *\dff \tau$\dss 
with\sss $\tau$\sss running over\dss the simplices of\qss $L$\nnsp.\oss
Clearly\halfff,\pss $K'$\dss is\dss a\sss simplicial\sss complex.\oss
As we will\sss see in a moment\halfff,\oss
$\norm{K'}\off =\off \norm{K}$\nnsp,\oss
i.e.\qss $K'$\dss is\dss a subdivision of\trs $K$\nnsp.\oss
The readers feeling\dss that\dss this\dss is\dss obvious\dss are advised\dss to
skip\dss the proof\dss of\trs the next\dss lemma.\oss

Let\dss us\dss
replace\dss in\dss $S$\dss the
simplices of\trs $K$\dss by\dss the simplices of\trs $K'$\nnsp.\oss
Since\dss $K'$\dss is\dss a\sss subdivision of\trs $K$\nnsp,\oss
the result\dss $S'$\dss
is\dss a\sss subdivision of\trs $S$\nnsp,\oss
called\dss the\qss \emph{stellar\dss subdivision}\qss of\trs $S$\dss
with\dss the\qss \emph{center}\qss $w$\dnsp.\vspace{1.25pt}

\myuppar{Lemma.} 
$\norm{K'}\off =\off \norm{K}$\nnsp.\oss\vspace{1.25pt}

\proof
By\dss the construction,\oss
$\norm{K'}\off \subset\off \norm{K}$\nnsp.\oss
Let\dss us\dss prove\sss the opposite inclusion.\oss 
Given\sss
a\dss point\qss $x\qff \in\qff \norm{K}$\nnsp,\oss
let\dss us\dss consider\dss the ray\sss starting at\dss $w$\dss and going\dss in\dss
the direction of\dss $x$\nnsp.\oss
The intersection of\trs this ray\dss with\dss $\norm{K}$\dss is\dss a\sss segment\dss $J$\dss
containing\dss $x$\dss and\dss having\dss $w$\dss as one of\dss its endpoints.\oss
Let\dss $z$\dss be\sss the other endpoint\sss of\trs $J$\nnsp,\oss
and\dss let\dss $\tau$\dss be\sss the smallest\sss simplex of\dss $K$\dss containing\dss $z$\nnsp.\oss
Then\dss $z$\dss belongs\sss to\sss the interior of\dss $\tau$\nnsp.\oss
If\dss $\tau$\dss has\dss $\sigma$\dss as a face,\oss
then\qss $w\qff \in\qff \tau$\qss and\dss the whole segment\dss $J$\dss is\dss
contained\dss in\dss $\tau$\nnsp.\oss
Moreover\halfff,\oss since\dss $z$\dss is\dss in\dss the interior\dss of\dss $\tau$\nnsp,\oss
the segment\dss $J$\dss can\dss be extended\dss without\dss leaving\dss $\norm{K}$\dss
and\dss hence cannot\dss be\sss the intersection of\dss our\dss ray\dss with\dss $\norm{K}$\nnsp.\oss
The contradiction shows\sss that\dss $\sigma$\dss is\dss not\sss a face of\dss $\tau$\dss
and\dss hence\dss $\tau$\dss is\sss a simplex of\trs $L$\nnsp.\oss
It\dss follows\dss that\qss
$x
\qff \in\qff 
w\dff *\dff \tau
\off \subset\off 
w\dff *\dff L$\nnsp.\oss
Since\qss $x\qff \in\qff \norm{K}$\qss was arbitrary\halfff,\oss
this proves\sss that\qss
$\norm{K'}\off =\off \norm{K}$\nnsp.\oss  \eproof

\myuppar{Stellar\sss subdivisions and chains.}
Let\dss us\dss keep\dss the above assumptions and choose a vertex\dss $v$\dss
of\dss $\sigma$\nnsp.\oss
Let\qss
$\varphi\dff(\dff w\dff)\off =\off v$\qss
and\qss 
$\varphi\dff(\dff z\dff)\off =\off z$\qss
if\dss $z$\dss is\dss a vertex of\trs $S'$\dss different\dss from\dss $w$\nnsp.\oss
Then\dss $\varphi$\dss is\dss a\sss simplicial\dss map\qss
$S'\qff \ttoo\qff S$\nnsp.\oss
Recall\dss that\dss $\fclass{\alpha}$\dss
denotes\sss the
subdivision of\dss a\sss chain\dss $\alpha$\dss of\trs $K$\dss
or\dss $S$\dss with respect\dss to\dss $K'$\dss or\dss $S'$\dss respectively\halfff.\oss
Since\sss the map\dss $\varphi$\dss is\dss obviously\dss pseudo-identical,\oss
Alexander's\dss lemma implies\sss that\qss
$\varphi_{\dff *}\dff(\qff \fclass{\alpha}\qff)
\off =\off
\alpha$\qss
for every chain\dss $\alpha$\dss in\dss $S$\nnsp.\oss
It\dss turns out\dss that\dss the map\qss
$\beta\qff \longmapsto\qff \fclass{\fff\varphi_{\dff *}\dff(\trf \beta\dff)}$\qss
is\dss fairly close\sss to\sss the identity\halfff.\oss
This\dss is\dss the key\dss element\sss of\trs the proof\dss of\trs
the invariance of\trs homology\dss groups under stellar subdivisions.\oss

We need\dss to extend\dss the operation of\trs
taking cones\sss to a simple situation when\dss the geometric cone is\dss not\sss defined.\oss
Let\dss $\tau$\dss be an $n$\dnsp-simplex of\trs $L$\dss for some\dss $n$\nnsp.\oss
If\dss $v$\dss is not\sss a vertex of\dss $\tau$\nnsp,\oss
then\dss $v\dff *\dff \tau$\dss is\dss already\dss well\sss defined\qss
({\fff}because\dss $v$\dss is\dss a vertex and\dss $\tau$\dss is\dss a face
of\dss some simplex of\trs $K$\nsp).\oss
If\dss $v$\dss is\dss a vertex of\dss $\tau$\nnsp,\oss
then\dss we interpret\dss $v\dff *\dff \tau$\dss
as\sss the zero $(\fff n\dff +\dff 1\fff)$\dnsp-chain.\oss
We extend\dss the maps\qss  
$\tau\qff \longmapsto\qff w\dff *\dff \tau$\qss
and\qss
$\tau\qff \longmapsto\qff v\dff *\dff \tau$\qss
from\sss simplices\sss to chains\sss by\dss linearity\halfff.\oss

\myuppar{Lemma.}
\emph{If\pss $\delta$\dss is\dss a\sss face and\dss $z$\dss
is\dss a vertex of\qss a\sss simplex,\oss
then}\qss
$\partial\qff (\dff z\dff *\dff \delta\trf)
\off =\off
\delta
\qff +\qff
z\dff *\dff \partial\dff \delta$\nnsp.\oss

\proof
If\dss $z$\dss is\dss not\sss a\sss vertex of\dss $\delta$\nnsp,\oss
this\dss is\sss obvious.\oss
If\dss $\delta$\dss is\dss an $n$\dnsp-simplex and\dss $z$\dss 
is\sss a\sss vertex of\dss $\delta$\nnsp,\oss
then\dss $z$\dss is\sss a\sss vertex of\dss all\sss $(\fff n\dff -\dff 1\fff)$\dnsp-faces
of\dss $\delta$\dss except\sss of\trs the $(\fff n\dff -\dff 1\fff)$\dnsp-face\dss 
$\delta_{\dff z}$\dss opposite\sss to\dss $z$\nnsp.\oss
In\dss this case\qss
$z\dff *\dff \delta
\off =\off
0$\qss
and\qss
$z\dff *\dff \partial\dff \delta\
\off =\off
z\dff *\dff \delta_{\dff z}
\off =\off
\delta$\nnsp.\oss
Therefore\sss in\dss this case\sss the identity\sss of\trs the lemma\dss reduces\sss
to\qss 
$\partial\dff 0\off =\off \delta\qff +\qff \delta$\nnsp,\oss
which\dss is\dss obviously\dss true over\dss $\ftwo$\nsp.\oss  \eproof

\myuppar{Double cones.}
Given a simplex\dss $\tau$\dss of\trs $L$\nnsp,\oss
let\qss\vspace{1.5pt}
\[
\quad
wv\dff *\dff \tau
\off =\off
w\dff *\dff (\dff v\dff *\dff \tau\dff) 
\qff
\]

\vspace{-10.5pt}
considered as a chain\dss in\dss $K'$\nnsp.\oss
Here\dss $v\dff *\dff \tau$\dss is\dss interpreted as\sss zero\dss
if\dss $v$\dss is\dss a\sss vertex of\dss $\tau$\dss
and\dss
$w\dff *\dff (\dff v\dff *\dff \tau\dff)$\dss
is\dss interpreted as zero\dss
if\dss $w$\dss belongs\dss to\sss the simplex\dss $v\dff *\dff \tau$\nnsp.\oss
Let\dss us extend\dss the map\qss\vspace{1.5pt}
\[
\quad
\tau\qff \longmapsto\qff wv\dff *\dff \tau
\]

\vspace{-10.5pt}
to\sss a\sss linear\dss map\dss from chains of\trs $L$\dss to\sss chains of\trs $K'$\nnsp.\oss

\myuppar{Lemma.} 
\emph{Let\pss $\alpha$\dss be\dss a chain\dss in\qss $L$\nnsp.\oss
Then}\vspace{2.5pt}
\[
\quad
\partial\dff (\dff wv\dff *\dff \alpha\dff)
\off =\off
\fclass{v\dff *\dff \alpha\dff}
\pff +\pff
w\dff *\dff \alpha
\pff +\pff
wv\dff *\dff (\dff \partial\dff \alpha \dff)
\qff.
\]

\vspace{-9pt}
\proof
It\dss is\dss sufficient\dss to consider\dss the case when\dss $\alpha$\dss
is\dss equal\dss to a simplex\dss $\tau$\dss of\trs $L$\nnsp.\oss
Suppose\sss first\dss that\sss $v$\dss is\dss not\sss a vertex of\dss $\tau$\dss
and\dss
$w$\dss does not\dss belong\sss to\dss $v\dff *\dff \tau$\nnsp.\oss
Then\dss
$v\dff *\dff \tau$\dss
is\dss a simplex\dss in\dss $L$\dss and\dss hence\qss
$\fclass{v\dff *\dff \tau\fff}
\off =\off
v\dff *\dff \tau$\nnsp.\oss
In\dss this case\sss the previous\sss lemma\sss  
implies\sss that\vspace{3pt}
\[
\quad
\partial\dff (\dff wv\dff *\dff \tau\dff)
\off =\off
\partial\dff (\dff w\dff *\dff (\dff v\dff *\dff \tau\dff)\trf)
\off =\off
v\dff *\dff \tau
\qff +\qff
w\dff *\dff (\trf \partial\dff (\dff v\dff *\dff \tau\dff)\trf)
\]

\vspace{-35.25pt}
\[
\quad
\phantom{\partial\dff (\dff wv\dff *\dff \tau\dff)
\off }
=\off
v\dff *\dff \tau
\qff +\qff
w\dff *\dff (\trf \tau\qff +\qff v\dff *\dff \partial\dff \tau \trf)
\off =\off
v\dff *\dff \tau
\qff +\qff
w\dff *\dff \tau
\qff +\qff
w\dff *\dff (\dff v\dff *\dff \partial\dff \tau\dff)
\qff.
\]

\vspace{-35.25pt}
\[
\quad
\phantom{\partial\dff (\dff wv\dff *\dff \tau\dff)
\off }
=\off
\fclass{v\dff *\dff \tau\fff}
\qff +\qff
w\dff *\dff \tau
\qff +\qff
wv\dff *\dff (\dff \partial\dff \tau \dff)
\qff.
\]

\vspace{-9pt}
Suppose now\dss that\dss $v$\dss is\dss not\sss a vertex of\dss $\tau$\nnsp,\oss
but\dss $w$\dss belongs\sss to\dss
$v\dff *\dff \tau$\nnsp.\oss
In\dss this case\dss $\sigma$\dss is\dss 
not\sss a face of\dss $\tau$\dss but\dss is\dss
a face of\dss $v\dff *\dff \tau$\nnsp.\oss
This may\dss happen only\dss when\dss $\tau$\dss
contains\sss the face\dss $\sigma_{\fff v}$\dss opposite\sss to\dss $v$\dss
in\dss $\sigma$\dss
but\sss does not\sss contains\dss $v$\nnsp.\oss
Let\dss $n$\dss be\sss the dimension of\dss $\tau$\nnsp.\oss
An $n$\dnsp-face of\dss $v\dff *\dff \tau$\dss
is\dss either equal\dss to\dss $\tau$\nnsp,\oss or\trs
has\sss the form\dss
$v\dff *\dff \lambda$\dss
for some $(\fff n\dff -\dff 1\fff)$\dnsp-face\dss $\lambda$\dss of\dss $\tau$\nnsp.\oss
The simplex\dss $v\dff *\dff \lambda$\trs 
is\sss a simplex of\trs $L$\dss if\trs and\dss only\trs if\dss $\lambda$\dss
does not\sss contain\dss the simplex\dss $\sigma_{\fff v}$\nsp.\oss
Let\dss $\Lambda$\dss be\sss the sum of\dss such $(\fff n\dff -\dff 1\fff)$\dnsp-faces\dss 
$\lambda$\dss of\dss $\tau$\nnsp.\oss
Then\qss 
$\fclass{v\dff *\dff \tau\fff}
\off =\off
w\dff *\dff (\dff v\dff *\dff \Lambda\dff)
\qff +\qff
w\dff *\dff \tau$\nnsp.\qff\oss
On\dss the other\dss hand,\oss\vspace{1.5pt}
\[
\quad
wv\dff *\dff (\dff \partial\dff \tau \dff)
\off =\off
w\dff *\dff (\dff v\dff *\dff \partial\dff \tau\dff)
\off =\off
w\dff *\dff (\dff v\dff *\dff \Lambda\dff)
\qff
\]

\vspace{-10.5pt}
because\dss if\trs an $(\fff n\dff -\dff 1\fff)$\dnsp-face\dss $\mu$\dss of\dss $\tau$\dss
contains\sss the simplex\dss $\sigma_{\fff v}$\nsp,\oss
then\dss $v\dff *\dff \mu$\dss contains\dss $w$\dss
and\dss hence\qss
$w\dff *\dff (\dff v\dff *\dff \mu\dff)
\off =\off
0$\qss
by\dss the definition.\oss
It\dss follows\dss that\vspace{1.75pt}
\[
\quad
\fclass{v\dff *\dff \tau\fff}
\qff +\qff
w\dff *\dff \tau
\qff +\qff
wv\dff *\dff (\dff \partial\dff \tau \dff)
\]

\vspace{-37.25pt}
\[
\quad
=\off
w\dff *\dff (\dff v\dff *\dff \Lambda\dff)
\qff +\qff
w\dff *\dff \tau
\qff +\qff
w\dff *\dff \tau
\qff +\qff
w\dff *\dff (\dff v\dff *\dff \Lambda\dff)
\off =\dff\off
0
\qff.
\]

\vspace{-10.25pt}
In\dss this case also\qss
$wv\dff *\dff \tau\off =\off 0$\qss
and\dss hence\sss the identity\sss of\trs lemma holds.\oss

It\dss remains\sss to consider\dss the case when\dss
$v$\dss is\dss a\sss vertex of\dss $\tau$\nnsp.\oss
Then\qss
$wv\dff *\dff \tau
\off =\off
w\dff *\dff 0
\off =\off
0$\qss
and\qss
$v\dff *\dff \tau
\off =\off
0$\nnsp.\oss
By\dss the previous\sss lemma\qss
$\tau
\off =\off
v\dff *\dff \partial\dff \tau$\qss
and\dss hence\qss
$wv\dff *\dff (\dff \partial\dff \tau \dff)
\off =\off
w\dff *\dff \tau$\nnsp.\oss
Since\qss
$w\dff *\dff \tau
\qff +\qff
w\dff *\dff \tau
\off =\off
0$\nnsp,\oss
the\dss lemma\dss holds in\dss this case also.\oss  \eproof

\myuppar{Lemma.}
\emph{Let\dss $\alpha$\dss be a chain of\pss $K'$\dss
such\dss that\dss its boundary\dss is a chain\dss of\pss $L$\nnsp.\oss
Then\qss
$\fclass{\varphi_{\dff *}\dff(\dff \alpha\dff)}
\qff -\qff
\alpha$\qss
is a cycle and,\oss moreover\halfff,\oss a boundary\halfff.\oss}

\proof
Let\dss $\beta$\dss be\sss the sum of\dss simplices of\dss $\alpha$\dss not\dss
having\dss $w$\dss as a vertex.\oss
Then\dss $\beta$\dss is\dss a chain\dss in\dss $L$\dss and\dss hence\qss
$\varphi_{\dff *}\dff(\trf \beta\dff)
\off =\off
\beta$\qss
and\qss
$\fclass{\varphi_{\dff *}\dff(\trf\beta\dff)}
\off =\off
\fclass{\beta}
\off =\off
\beta$\nnsp.\oss
It\dss follows\dss that\dss
the boundary\sss of\dss $\beta$\dss is\dss a chain\dss in\dss $L$\dss
and\qss
$\fclass{\varphi_{\dff *}\dff(\trf \beta\dff)}
\qff -\qff
\beta
\off =\off
0$\nnsp.\oss
Therefore,\oss after\dss replacing\dss $\alpha$\dss by\qss
$\alpha\qff -\qff \beta$\nnsp,\oss
if\dss necessary\halfff,\oss
we may\sss assume\dss that\sss every\sss simplex of\dss $\alpha$\dss
has\dss $w$\dss as a vertex.\oss
In\dss this case\qss
$\alpha
\off =\off
w\dff *\trf \rho$\qss
for some chain\dss $\rho$\dss of\trs $L$\dss
and\dss hence\qss
$\varphi_{\dff *}\dff(\dff \alpha\dff)
\off =\off
v\dff *\dff \rho$\nnsp.\oss
Also,\oss
$\partial\dff \alpha$\dss is\dss a chain\dss in\dss $L$\dss
and\dss hence
$\partial\dff \alpha\off =\off \rho$\nnsp.\oss
It\dss follows\dss that\qss
$\partial\dff \rho
\off =\off 
\partial\dff \partial\dff \alpha
\off =\off
0$\nnsp.\oss 
By\dss the previous\sss lemma\qss\vspace{1.75pt}
\[
\quad
\partial\dff (\dff wv\dff *\dff \rho\dff)
\off =\off
\fclass{v\dff *\dff \rho\fff}
\pff +\pff
w\dff *\dff \rho
\pff +\pff
wv\dff *\dff (\dff \partial\dff \rho \dff)
\]

\vspace{-37.25pt}
\[
\quad
\phantom{\partial\dff (\dff wv\dff *\dff \tau\dff)
\off }
=\off
\fclass{v\dff *\dff \rho\fff}
\pff +\pff
w\dff *\dff \rho\fff
\off =\off
\fclass{\varphi_{\dff *}\dff(\dff \alpha\dff)}
\pff +\pff
\alpha
\qff
\]

\vspace{-10.25pt}
and\dss hence\qss $\fclass{\varphi_{\dff *}\dff(\dff \alpha\dff)}
\pff -\pff
\alpha
\off =\off
\fclass{\varphi_{\dff *}\dff(\dff \alpha\dff)}
\pff +\pff
\alpha
$\qss is\dss a\sss boundary\halfff.\oss  \eproof

\myuppar{Lemma.}
\emph{If\qss $\gamma$\dss is a cycle in\dss $S'$\nnsp,\oss
then\qss
$\fclass{\varphi_{\dff *}\dff(\dff \gamma\dff)}
\qff -\qff
\gamma$\qss is\dss a\dss boundary\halfff.\oss}

\proof
Let\dss $\alpha$\dss be\sss the sum of\dss simplices of\dss $\gamma$\dss
having\dss $w$\sss as a vertex,\oss
and\dss let\dss $\beta$\dss be\sss the sum of\trs other simplices of\dss $\gamma$\nnsp.\oss
Then\qss
$\gamma\off =\off \alpha\qff +\qff \beta$\nnsp,\oss
the chain\dss $\alpha$\dss is\dss a chain\dss in\dss $K'$\nnsp,\oss
and\dss each simplex of\dss $\beta$\dss belonging\dss to\dss $K'$\dss
actually\dss belongs\sss to\dss $L$\nnsp.\oss
It\dss follows\dss that\qss
$\fclass{\varphi_{\dff *}\dff(\trf\beta\dff)}
\off =\off
\fclass{\beta}
\off =\off
\beta$\qss
and\dss hence\qss 
$\fclass{\varphi_{\dff *}\dff(\dff \gamma\dff)}
\qff -\qff
\gamma
\off =\off
\fclass{\varphi_{\dff *}\dff(\dff \alpha\dff)}
\qff -\qff
\alpha$\nnsp.\oss
Since\dss $\gamma$\dss is\dss a cycle,\oss
$\partial\dff \alpha
\off =\off
-\qff \partial\trf \beta$\nnsp.\oss
It\dss follows\sss that\dss the boundaries\qss $\partial\dff \alpha\fff,\off \partial\trf \beta$\qss
are chains in\dss $L$\nnsp.\oss
It\dss remains\sss to apply\dss the previous lemma.\oss  \eproof

\myuppar{Theorem.}
\emph{For\sss every\qss $m\qff \geq\qff 0$\qss the\dss homology\dss subdivision\dss map\qss
$s
\qff \colon\dff
H_{\fff m}\fff(\dff S \dff)
\qff \ttoo\qff
H_{\fff m}\fff(\dff S' \dff)$\qss
and\dss the\sss induced\dss map\oss
$\varphi_{*\nsp *}
\qff \colon\dff
H_{\fff m}\fff(\dff S' \dff)
\qff \ttoo\qff
H_{\fff m}\fff(\dff S \dff)$\qss
are mutually\dss inverse isomorphisms.\oss}

\proof
By\qss Alexander's\dss lemma\dss 
$\varphi_{\dff *}\dff(\qff \fclass{\alpha}\qff)
\off =\off
\alpha$\qss
for every chain\dss $\alpha$\dss in\dss $S$\dss and\dss hence\dss 
$\varphi_{*\nsp *}\dff \circ\dff s$\dss
is\dss the identity\dss map.\oss  
By\dss the last\dss lemma for every\sss cycle\dss $\gamma$\dss in\dss $S'$\dss
the cycles\dss $\fclass{\varphi_{\dff *}\dff(\dff \gamma\dff)}$\dss
and\dss $\gamma$\dss belong\dss to\sss the same homology\sss class
and\dss hence\qss
$s\dff \circ\dff \varphi_{*\nsp *}$\qss
is\dss the identity\dss map.\oss  \eproof

\myuppar{Stellar\dss subdivisions\sss and centers-generated\dss subdivisions.}
Let\dss $S$\dss be a simplicial complex.\oss
Suppose\sss that\halfff,\oss as\dss in\dss Appendix\qss \ref{barycent},\oss
for every simplex\dss $\sigma$\dss of\qss $S$\dss
a point\qss $z\dff(\dff \sigma\dff)\qff \in\qff \sigma\qff \smallsetminus\qff \bd\dff \sigma$\nnsp,\oss
called\dss the\qss \emph{center}\qss of\dss $\sigma$\nnsp,\oss
is chosen.\oss
Recall\dss that\sss such a choice generates a subdivision\dss $\cs S$\dss of\qss $S$\nnsp.\oss
The simplices of\dss $\cs S$\dss are in one-to-one correspondence with sequences\qss
$\sigma_{\dff 0}\dff,\off
\sigma_{\dff 1}\dff,\off
\ldots\dff,\off
\sigma_{\dff n}$\qss
of\dss simplices of\trs $S$\dss such\dss that\dss $\sigma_{ j}$\dss
is\dss a\dss face of\dss $\sigma_{\dff i}$\qss if\pss $j\qff <\qff i$\nnsp,\oss
and\dss the simplex corresponding\dss to\sss the sequence\qss
$\sigma_{\dff 0}\dff,\off
\sigma_{\dff 1}\dff,\off
\ldots\dff,\off
\sigma_{\dff n}$\qss
has\sss the centers\qss\qss\vspace{0pt}
\[
\quad
z\dff(\dff \sigma_{\dff 0}\dff)\fff,\off
z\dff(\dff \sigma_{\dff 1}\dff)\fff,\off
\ldots\fff,\off
z\dff(\dff \sigma_{\dff n}\dff)
\]

\vspace{-12pt}
as its vertices.\oss
Alexander\dss observed\dss that\sss such a subdivision\dss $\cs S$\dss
can\dss be obtained as\sss the result\sss of\dss a sequence of\dss stellar\sss subdivisions.\oss
Let\dss us\dss arrange all\sss simplices of\trs $S$\dss
into a sequence of\dss simplices\qss
$\sigma_{\dff 1}\dff,\off
\sigma_{\dff 2}\dff,\off
\ldots\dff,\off
\sigma_{\dff N}$\qss
with\dss non-increasing\sss dimensions.\oss
In other\dss words,\oss
the dimension of\dss $\sigma_{\dff i}$\dss
is required\dss to be greater or equal\dss than\dss the dimension of\dss
$\sigma_j$\qss if\pss $i\qff \leq\qff j$\nnsp.\qff\oss
The order of\dss simplices of\trs the same dimension does not\dss matter\halfff.\oss
Let\dss $S\dff(\dff 0\dff)\off =\off S$\qss
and\qss let\dss $S\dff(\dff k\trf)$\dss
be\sss the stellar subdivision of\trs $S\dff(\dff k\qff -\qff 1\dff)$\dss
with\dss the center\dss $z\dff(\dff \sigma_{\fff k}\trf)$\dss
for every\qss $k\qff \leq\qff N$\nnsp.\oss\vspace{-1.5pt}

\myuppar{Lemma.}
\emph{\dnsp$S\dff(\trf N\trf)$\dss is equal\dss to\sss the 
centers-generated\sss subdivision\dss $\cs S$\nnsp.\oss}\vspace{-1.5pt}

\proof
The vertices of\trs $S\dff(\trf N\trf)$\dss are exactly\dss the centers of\dss
simplices of\trs $S$\sss
(one should\dss keep\dss in\dss mind\dss the vertices of\trs $S$\dss are centers
of\trs the corresponding $0$\dnsp-simplices).\oss
Since several\dss vertices are\sss the vertices of\dss a simplex\dss
if\trs and\dss only\trs if\qss they\sss are pairwise connected\dss by edges,\oss
the main\dss part\sss of\trs the proof\dss is\dss to find out\dss when\dss
two centers are connected\dss by\sss an edge in\dss $S\dff(\trf N\trf)$\nnsp.\oss

Let\qss $k\qff \leq\qff N$\nnsp.\oss 
The center\dss $z\dff(\dff \sigma_{\fff k} \trf)$\dss is\dss introduced
as a new\sss vertex\sss in\dss $S\dff(\dff k\trf)$\nnsp.\oss
Since dimensions are non-increasing,\oss the centers\dss $z\dff(\dff \rho\dff)$\dss
of\trs simplices\dss $\rho$\dss having\dss $\sigma_{\fff k}$\dss as a proper\dss face
are already\dss present\dss in\dss $S\dff(\dff k\trf)$\nnsp.\oss
By\dss the definition of\trs stellar\sss subdivisions,\pss
$z\dff(\dff \sigma_{\fff k}\dff)$\dss is\dss connected\dss by\sss an edge in $S\dff(\dff k\trf)$
to\sss these centers\dss $z\dff(\dff \rho\dff)$\dss 
and\dss to no other centers\dss $z\dff(\dff \sigma_{\fff i}\dff)$\dss
with\qss $i\qff <\qff k$\nnsp.\oss
In\dss particular\halfff,\pss $z\dff(\dff \sigma_{\fff k} \trf)$\dss is\dss not\dss connected\dss
in\dss $S\dff(\dff k\trf)$\dss with centers\dss $z\dff(\dff \sigma_{\fff i}\trf)$\dss
of\dss simplices\dss $\sigma_{\fff i}$\dss of\trs the same dimension as\dss $\sigma$\nnsp.\oss

Since every\sss edge connecting\dss two centers is created at\sss one of\dss the steps
of\dss our stellar\sss subdivision\dss process,\oss
we see\sss that\dss two centers\dss $z\dff(\dff \sigma_{\fff k} \dff)$\dss
and\dss $z\dff(\dff \sigma_{\fff i} \dff)$\dss 
are connected\dss by\sss an edge\dss in\dss $S\dff(\trf N\trf)$\dss
if\trs and\dss only\trs if\qss one of\trs the simplices\qss
$\sigma_{\fff k}\dff,\pff \sigma_{\fff i}$\qss
is\dss a\sss proper\dss face of\trs the other\halfff.\oss
It\dss follows\dss that\dss several\dss centers are pairwise connected\dss
by\sss edges of\trs $S\dff(\trf N\trf)$\dss
if\trs and\dss only\trs if\qss
the corresponding simplices can\dss be arranged\dss in\sss a\sss sequence\qss
$\sigma_{\dff 0}\dff,\off
\sigma_{\dff 1}\dff,\off
\ldots\dff,\off
\sigma_{\dff n}$\qss
such\dss that\dss $\sigma_{ j}$\dss
is\dss a\dss face of\dss $\sigma_{\dff i}$\qss if\pss $j\qff <\qff i$\nnsp.\oss
This means\sss that\dss $S\dff(\trf N\trf)$\dss has exactly\dss the same
simplices as\dss $\cs S$\nnsp.\oss  \eproof

\myuppar{Theorem.}
\emph{The\dss homology\dss subdivision\dss map\qss
$s
\qff \colon\dff
H_{\fff m}\dff(\dff S \dff)
\qff \ttoo\qff
H_{\fff m}\dff(\dff \cs S \dff)$\qss
is\dss an\dss isomorphism.\oss}

\proof
Since\sss the composition
of\trs the homology\sss subdivision\dss maps\qss\vspace{1.5pt}
\[
\quad
H_{\dff m}\trf\bigl(\qff S\dff(\dff k\qff -\qff 1\dff) \qff\bigr)
\qff \ttoo\qff
H_{\dff m}\trf\bigl(\qff S\dff(\dff k\dff)  \qff\bigr) 
\]

\vspace{-10.5pt}
is equal\dss to\dss $s\qff \colon\dff
H_{\fff m}\fff(\dff S \dff)
\qff \ttoo\qff
H_{\fff m}\fff(\dff \cs S \dff)$\nnsp,\oss
we need only\dss to apply\dss the previous lemma.\oss  \eproof

\begin{flushright}

September\halfff\qss 1,\oss 2019

https\halfff:/\!/\hspace*{-0.06em}nikolaivivanov.com

E-mail\halfff:\oss nikolai.v.ivanov{\fff}@{\dff}icloud.com

\end{flushright}

\end{document}